%% file: Vanishing1D.tex
\documentclass[10pt]{article}
\usepackage[english,activeacute]{babel}
\usepackage{epsfig}

\usepackage[latin1]{inputenc}
\usepackage[T1]{fontenc}
\usepackage{stmaryrd}
\usepackage{amsmath,amsfonts,amssymb,mathrsfs,amsthm}
\usepackage{xcolor}
\usepackage{dsfont}
\usepackage{graphicx}
\usepackage[mathscr]{eucal}
\usepackage{makeidx}
\usepackage{verbatim}
\usepackage{graphics,graphicx}
\usepackage{textcomp}
\usepackage{float}
\usepackage{tikz}
\usepackage[colorlinks=true, linkcolor=purple, citecolor=purple]{hyperref}

\input mymacros.tex

\DeclareMathOperator{\arcsinh}{arcsinh}

\DeclareMathOperator{\unif}{unif} 

\newcommand\bna{\begin{eqnarray*}}
\newcommand\ena{\end{eqnarray*}}

\def\bal#1\nal{\begin{align*}#1\end{align*}}
\def\baln#1\naln{\begin{align}#1\end{align}}

\newcommand\bnan{\begin{eqnarray}}
\newcommand\enan{\end{eqnarray}}

\newcommand\bnp{\begin{proof}}
\newcommand\enp{\end{proof}}

\newcommand\bneq{\begin{eqnarray*}\left\lbrace \begin{array}{rcl}}
\newcommand\eneq{\end{array} \right.\end{eqnarray*}}
\newcommand\bneqn{\begin{eqnarray}\left\lbrace \begin{array}{rcl}}
\newcommand\eneqn{\end{array} \right.\end{eqnarray}}

 \numberwithin{equation}{section} 
 
\newcommand\grando[1]{\mathcal{O}\left(#1\right)}

\renewcommand\M{[0,L]}

\newcommand\nor[2]{\left\|#1\right\|_{#2}}

\newcommand\e{\varepsilon}

\newcommand{\q}{q}
\newcommand{\qf}{q_{\f}}

\newcommand{\EE}{E}

\newcommand{\TO}{\ensuremath{T}}
\newcommand{\xzero}{\ensuremath{\mathbf{x}_0}}

\begin{document}
\title{On uniform controllability of 1D transport equations in the vanishing viscosity limit}

\author{Camille Laurent\footnote{CNRS UMR 7598 and Sorbonne Universit\'es UPMC Univ Paris 06, Laboratoire Jacques-Louis Lions, F-75005, Paris, France, email: camille.laurent@sorbonne-universite.fr} and Matthieu L\'eautaud\footnote{Laboratoire de Math\'ematiques d'Orsay, UMR 8628, Universit\'e Paris-Saclay, CNRS, B\^atiment 307, 91405 Orsay Cedex France, email: matthieu.leautaud@universite-paris-saclay.fr.}
}

\date{}

\maketitle

\begin{abstract}

We consider a one dimensional transport equation with varying vector field and a small viscosity coefficient, controlled by one endpoint of the interval. We give upper and lower bounds on the minimal time needed to control to zero, uniformly in the vanishing viscosity limit.

We assume that the vector field varies on the whole interval except at one point. The upper/lower estimates we obtain depend on geometric quantities such as an Agmon distance and the spectral gap of an associated semiclassical Schr\"odinger operator. They improve, in this particular situation, the results obtained in the companion paper~\cite{LL:18vanish}.

The proofs rely on a reformulation of the problem as a uniform observability question for the semiclassical heat equation together with a fine analysis of localization of eigenfunctions both in the semiclassically allowed and forbidden regions~\cite{LL:22}, together with estimates on the spectral gap~\cite{HS:84,Allibert:98}.
Along the proofs, we provide with a construction of biorthogonal families with fine explicit bounds, which we believe is of independent interest.  

\end{abstract}


\tableofcontents

\section{Introduction and main results}

We consider the one dimensional diffusive-transport equation, controlled from the left endpoint of the interval:
\begin{equation}
\label{e:1Dcontrol-problemyn}
\left\{
\begin{array}{rl}
(\d_t + \mathfrak{a} \d_x + \mathfrak{b}  - \eps\d_x^2) y= 0 ,& (t,x)\in (0,T) \times (0,L) , \\
y (t,0) = h(t) , \quad y(t,L) = 0 ,&  t \in (0,T) , \\
y(0,x) = y_0 (x),&  x \in (0,L) .
\end{array}
\right.
\end{equation}
Here $L>0$ is the length of the spatial domain, $T>0$ the time horizon, and $\eps>0$ a viscosity parameter. The functions $\mathfrak{a}, \mathfrak{b}  : [0,L] \to \R$ are real-valued and  sufficiently regular. We shall later on rewrite this equation as $$
(\d_t + \f' \d_x +\f'' - \q  - \eps\d_x^2) y= 0 ,
$$
that is to say write $\mathfrak{a} =  \f'$ and $\mathfrak{b}=\f''- \q$ for simplicity of the dual equation and consistency with the companion article~\cite{LL:18vanish}.

For an initial datum $y_0\in L^2(0,L)$ and a control function $h \in L^2(0,T)$, it is known that \eqref{e:1Dcontrol-problemyn} has a unique solution in $C^0(0,T;L^2(0,L))$ in the sense of transposition (see~\cite{FR:71} or~\cite{CG:05}). 
The usual question of null-controllability is whether, the parameters $T,L,\eps$ being fixed, one can drive any initial datum $y_0$ to rest (i.e. the null function) in time $T$ by means of the action on the equation through the function $h(t)$. 

 \begin{definition}[Controllability and cost]
 \label{d:def-0}
  Given $(\eps, T)$, we say that~\eqref{e:1Dcontrol-problemyn} is null-controllable if for any $y_0 \in L^2(0,L)$, there is $h  = h(T,\eps,y_0) \in L^2(0,T)$ such that the associated solution to~\eqref{e:1Dcontrol-problemyn} satisfies $y(T) =0$. 
We define for $y_0 \in L^2(0,L)$ the (possibly empty, closed convex) set $U(y_0)$ of all such controls $h \in L^2(0,T)$, and the cost function 
  $$
  \mathcal{C}_0(T,\eps) := \sup_{y_0 \in L^2(0,L) , \nor{y_0}{L^2(0,L)}\leq 1} \left\{\inf_{h \in U(y_0)} \nor{h}{L^2(0,T)} \right\}  \in [0,+ \infty].
  $$  
  We have $\mathcal{C}_0(T,\eps) <+\infty$ if~\eqref{e:1Dcontrol-problemyn} is null-controllable, and $\mathcal{C}_0(T,\eps) =+\infty$ if not.
   \end{definition}
It is known (see~\cite{FR:71}, or~\cite{FI:96,LR:95,Lea:10} in higher dimension) that for fixed $\eps>0$, the equation~\eqref{e:1Dcontrol-problemyn} is null-controllable in any positive time $T>0$. That is to say, $T,\eps > 0 \implies \mathcal{C}_0(T,\eps) < + \infty$.
 This is linked to the infinite speed of propagation for the heat dissipation.
Here, we address the question of {\em uniform controllability} in the vanishing viscosity limit $\eps \to 0^+$, that is: how does $  \mathcal{C}_0(T,\eps)$ behave for fixed $T>0$ in the limit $\eps \to 0^+$?
This question has first been addressed by Coron and Guerrero~\cite{CG:05} in the case $\mathfrak{a}(x)=M$ (that is to say $\f(x) =Mx$) and $\mathfrak{b}=0$, for $M\in \R^*$, and different behaviors are observed, depending on the sign of $M$. 
In that paper, the authors make a conjecture on the minimal time needed to achieve uniform controllability, i.e.
$$
T_{\unif} (\{0\}) : = \min \{T>0 , \text{ there is }K>0 \text{ such that }   \mathcal{C}_0(T,\eps)\leq K \text{ for all }\eps\in (0,1) \}.
$$
Then, the estimates on this minimal time have been improved in~\cite{Gla:10,Lissy:12,Lissy:14,Lissy:15,Darde:17,YM:19,YM:19bis} with different methods. The result of~\cite{CG:05} was also generalized in several space dimensions and for non-constant transport speed in~\cite{GL:07}. In that paper however, no estimates on the minimal time are given. The first estimates on the minimal time needed for having $ \mathcal{C}_0(T,\eps)$ uniformly bounded as $\eps \to 0^+$ are proved in~\cite{LL:18vanish}, in a setting close to that of the present article. In particular, we exhibited in~\cite{LL:18vanish} higher dimensional situations in which $\frac{T_{\unif}}{T_{\f'}}$ can be as large as desired, where $T_{\f'}$ denotes the minimal time for the controllability of the limit transport equation obtained by formally taking $\eps=0$ in~\eqref{e:1Dcontrol-problemyn} (see Proposition~\ref{p:tpstransp} below for a more precise definition in our 1D context).

 Such uniform control properties in singular limits are also addressed for vanishing dispersion in~\cite{GG:08} and for vanishing dispersion and viscosity in~\cite{GG:09}. 
Controllability problems for nonlinear conservation laws with vanishing viscosity have also been studied in~\cite{GG:07} and~\cite{Lea:11}.
Motivation for studying the vanishing viscosity limit comes from different fields of mathematics:
\begin{itemize}
\item  conservation laws, for which the vanishing viscosity criterium is a selection principle for the physical (called entropy) solution, see~\cite{Kru:70} or~\cite[Chapter 6]{Daf:00}.
\item  control theory, where the study of singular limits sometimes allows to prove controllability properties for the perturbated system itself. See e.g. the papers~\cite{Cor:96,CF:96,ChaNavier:09,CMS:19}, where the authors investigate the Navier-Stokes system with Navier slip boundary conditions, relying on results for the Euler equation.
\item theoretical physics and differential topology, through the Witten-Helffer-Sj\"ostrand theory~\cite{Witten:82,HS:85}.
\item molecular dynamics and statistical physics, via the study of the so-called overdamped Langevin process~\cite{Chan:43,SM:79}.
\end{itemize}
We refer to~\cite[Section~1.2]{LL:18vanish} for more details on motivation.
 Our main results in the present article (namely Theorems~\ref{t:estim-Cobs-expo-1D},~\ref{thmlower+-} and~\ref{thmcontrol1D} below) formulate as explicit (in geometric terms, under some assumptions on the parameters) lower and upper bounds on the cost function $\mathcal{C}_0(T,\eps)$ and the minimal time $T_{\unif}$ of uniform controllability.
  We now give a list of geometric assumptions and related definitions in order to state our main results. 

\subsection{Definitions and assumptions}

\begin{figure}[h!]
  \begin{center}
\begin{tikzpicture}
\draw[->] (-5,0) -- (6,0) node[right] {$x$}; 
\draw[->] (-3,-1) -- (-3,5) ; 
\draw[-] (3,-1) -- (3,5);  

\begin{scope}
\clip (-5,-1) rectangle (4,5);
\draw[color=red,samples=100] plot ({\x},{1/5*(\x-1)*(\x-1)+1});
\end{scope}

\draw[color=blue,dashed] (1,-0.5) -- (1,1.5);
\draw[color=blue,dashed] (-3.5,1) -- (3.5,1);
\node[below left=0.1cm] at (-3,-0) {$0\strut$};
\node[below left=0.1cm] at (3,0) {$L\strut$};
\node[color=red,below right=0.1cm] at (3.8,3) {$V(x)=\frac{\mathfrak{a}(x)^2}{4} =\frac{|  \f'(x)|^2}{4}\strut$};
\node[color=blue,above left=0.1cm] at (-3,1) {$E_0=V(\xzero)\strut$};
\node[color=blue,below left=0.1cm] at (1,0) {$\xzero\strut$};

\end{tikzpicture}
    \caption{Geometric setting of Assumption~\ref{assumptions}}
    \label{f:geom-setting}
 \end{center}
\end{figure}
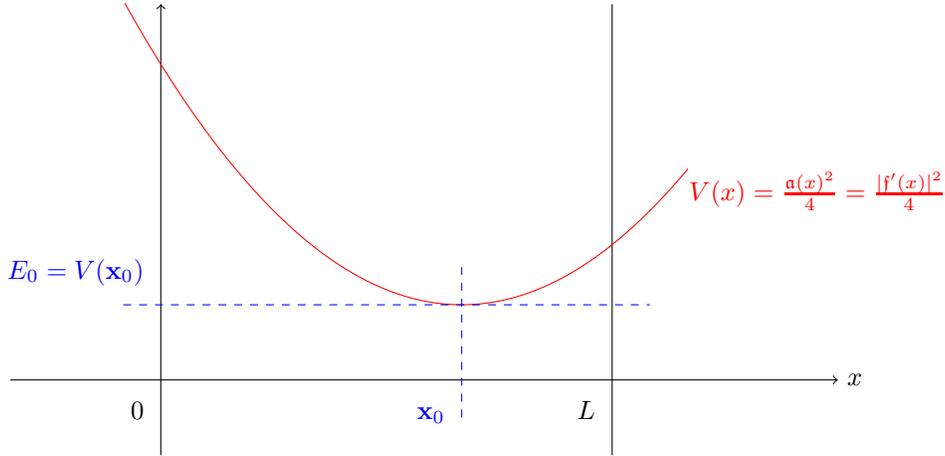

All along the paper, we make intensive use of the effective potential
\bna
V(x)= \frac{\mathfrak{a}(x)^2}{4} =\frac{|  \f'(x)|^2}{4}  .
\ena
In the results presented below, we make (at least part of) the following assumptions, essentially saying that $V$ forms a single non-degenerate well and does not vanish. This assumption is illustrated in Figure~\ref{f:geom-setting}.
\begin{assumption}
\label{assumptions}
With $V(x) =\frac{\mathfrak{a}(x)^2}{4} = \frac{|\f'(x)|^2}{4}$ for $x\in [0,L]$,
\begin{enumerate}
\item \label{A1} $V> 0$ on $[0,L]$;
\item \label{A2} the only $x \in [0,L]$ such that $V'(x) = 0$ is $x= \xzero \in (0,L)$ and $V(\xzero) = \min_{[0,L]}V$;
\item \label{A3} $V(L) \neq V(0)$ ; 
\item \label{A4} $V''(\xzero) > 0$.
\end{enumerate}
\end{assumption}
Assumption~\ref{assumptions}, formulated for simplicity on the potential $V=  \frac{\mathfrak{a}(x)^2}{4} =\frac{|  \f'(x)|^2}{4}$ can be formulated equivalently as (the equivalence is not one to one; however $1$-$i$ above is equivalent to $1$-$i$ below, for all $i \in \{1,\cdots,4\}$):
\begin{enumerate}
\item $\mathfrak{a}(x) \neq 0$ on $[0,L]$ (resp. $\f' \neq 0$ on $[0,L]$);
\item the only $x \in [0,L]$ such that $\mathfrak{a}'(x) = 0$ is $x= \xzero \in (0,L)$ and $|\mathfrak{a}(\xzero)| = \min_{[0,L]}|\mathfrak{a}|$\\
 (resp.  the only $x \in [0,L]$ such that $\f''(x) = 0$ is $x= \xzero \in (0,L)$ and $|\f'(\xzero)| = \min_{[0,L]}|\f'|$);
\item $\mathfrak{a}(0) \neq \mathfrak{a}(L)$ (resp. $\f'(L) \neq \f'(0)$); 
\item $\mathfrak{a}''(\xzero) \neq 0$.
\end{enumerate}

\begin{remark}
Notice that Assumption~\ref{assumptions} does not cover the classical constant speed case, which is largely considered in the literature, see e.g.~\cite{CG:05,Gla:10,Lissy:12,Lissy:14,Lissy:15,Darde:17,YM:19,YM:19bis}. The latter corresponds to $\f(x)=\pm Mx$ and thus to the ``flat potential'' $V=\frac{M^2}{4}$. However, a formal asymptotics will be considered in Section \ref{s:cascst}, starting with a family of potentials satisfying Assumption~\ref{assumptions} converging to the flat potential. 
This allows to compare our situation to the ``flat'' one and shows that our results formally recover (sometimes with slightly less accurate constants) the previously known results for this example.
Our purpose in the present paper (together with~\cite{LL:18vanish}) is not to revisit or consider a perturbation of the usual ``flat'' setting, but rather to reveal effects that are not present in the  ``flat'' case. 

Indeed, the class of vector fields (or functions $\f$) presented here (see Section \ref{s:explicit-comput} for more concrete examples and calculations) will allow to stress the fact that the convexity is responsible for a concentration of some eigenfunctions close to the minimum, which is not the case for the more studied case $\f(x)=\pm Mx$.

Note finally that the result of Guerrero-Lebeau applies in our context, as soon as $\f'$ does not vanish on $[0,L]$, and implies that $T_{\unif}(\{0\})<+\infty$. The goal of the present article is to give a more precise estimate on the quantity $T_{\unif}(\{0\})$ (under the additional Assumption~\ref{assumptions} on $\f'$).
\end{remark}
We also denote by $E_0$ the ground state energy, that is to say
\begin{align}
\label{e:def-E0}
E_0 =\min_{x\in \M} (V(x)) = V(\xzero) =\frac{\mathfrak{a}(\xzero)^2}{4} = \frac{|\f'(\xzero)|^2}{4} .
\end{align}
Let us finally describe geometric and spectral quantities appearing in the statements below. 
The classically allowed region at energy $E$ for the potential $V$ is defined by: 
$$
 K_E = \{x \in \M , V(x) \leq E \} .
$$
We may then define the Agmon distance (see e.g.~\cite[Chapter~3]{Helffer:booksemiclassic}) to the set $K_E$ at the energy level $E$ by
\begin{align*}
d_{A,E}(x) = \inf_{y \in K_E} \left|\int_{y}^{x} \sqrt{\left(V(s) - E \right)_+ } ds\right|,
\end{align*}
that is, the distance to the set $K_E$ for the (pseudo-)metric $(V-E)_+ $ where $\left(V(x) - E \right)_+ = \max \left(V(x) - E , 0\right)$. 
Note that $d_{A,E}$ vanishes identically on $K_E$ (and only on this set).
Under Assumption~\ref{assumptions} Item~\ref{A2}, we have for $E\geq E_0$ 
\bnan
\label{defAgmonbis}
d_{A,E}(x) =  \left|\int_{y}^{x} \sqrt{\left(V(s) - E \right)_+ } ds\right| ,
\enan
where $y$ is any point in $K_E$. 
Another important function in the estimates below is given by
\bnan
\label{e:def-WE}
W_E (x)=d_{A,E}(x)+\frac{\f(x)}{2} .
\enan
The following classical quantities of the Hamiltonian 
\begin{align}
\label{e:def-p}
p(x,\xi)=\xi^2+V(x) =  \xi^2+ \frac{|\f'(x)|^2}{4}.
\end{align}
enter into play in the spectral analysis of the operators involved (and are defined assuming Item~\ref{A2} in Assumption~\ref{assumptions}):
\begin{align}
\label{e:def-phi}
\Phi (E) &: = \int_{x_-(E)}^{x_+(E)} \sqrt{E - V(s)} ds  ,  \quad \text{ for } E\in[\min_{[0,L]}V ,  +\infty) \\
\label{e:def-TT1}
T_1& : =\sup_{E \geq V(\xzero)}T(E)  , \quad \text{with} \quad  T(E) :=  2 \int_{x_-(E)}^{x_+(E)} \frac{\sqrt{E}}{\sqrt{E -V(s)}} ds  .
\end{align}
In these expressions, for $E\geq E_0$, the points $x_\pm(E)$ are such that $K_E = [x_-(E),x_+(E)]$.
Namely, $x_-(E)$ denotes the solution to $V(x_-(E))=E$ which is $\leq \xzero$ for $E \leq V(0)$, and  $x_-(E)=0$ for $E \geq V(0)$.   
Similarly, $x_+(E)$ is the solution to $V(x_+(E))=E$ which is $\geq \xzero$ for $E \leq V(L)$, and  $x_+(E)=L$ for $E \geq V(L)$ (with $\xzero=x_-(E_0)=x_+(E_0)$ if $E=E_{0}$). 
 The geometric content of these quantities, as well as links between them are discussed in Section~\ref{s:section-allibert} below (in particular, $T_1$ is not homogeneous to a time, but we keep this notation issued from~\cite{Allibert:98}).

\subsection{Results}
In the statements below, we recall that $\mathfrak{a}=\f'$ or that $\f(x) = \int_0^x\mathfrak{a}(s)ds$ (all results stated with the function $\f$ are invariant by $\f \mapsto \f+c$ for $c\in \R$). Also, we have $\mathfrak{b} = \f''-\q$ or equivalently $\q=\mathfrak{a}'-\mathfrak{b}$.
A first lower bound is as follows. Recall that the control cost $\mathcal{C}_0(\TO, \eps)$ is introduced in Definition~\ref{d:def-0}.
\begin{proposition}
\label{p:transport-seul}
Assuming $\mathfrak{a} \in C^1([0,L])$ (resp. $\f \in C^2([0,L])$) and setting $T_{\mathfrak{a}} :=\int_0^L \frac{ds}{|\mathfrak{a}(s)|}=\int_0^L \frac{ds}{|\mathfrak{f}'(s)|} \in (0,+\infty]$, we have 
$$
T < T_{\mathfrak{a}}  \quad \implies \quad  \liminf_{\eps \to 0^+}\mathcal{C}_0(\TO, \eps) = +\infty .
$$
\end{proposition}
This result is a direct consequence of the weak convergence of solutions to~\eqref{e:1Dcontrol-problemyn} to those of a limit problem with $\eps=0$ (proof of this weak convergence follows~\cite[Proposition~2.94]{Cor:book} and the limit equation is studied in Section~\ref{s:limit-eq} below)\footnote{Indeed, assuming $\mathcal{C}_0(\TO, \eps_n)\leq C_0$  with $\eps_n\to 0$, then, for any $y_0\in L^2(0,L)$, we can extract a sequence of controls $u_n$ converging in in $\R$ and of solutions $y_n$ converging weakly in $L^2((0,T)\times(0,L))$. By \cite[Proposition~2.94]{Cor:book}, the weak limit of $y_n$ is a solution of the transport equation studied in Section~\ref{s:limit-eq} below, controlled to zero in time $\TO$. Since this holds for any $y_0$, we necessarily deduce $T\geq T_{\mathfrak{a}}$ (see Lemma~\ref{p:tpstransp}).}.

It simply translates the fact that if the ``limit transport equation'' with $\eps=0$ is non-controllable, then there is no hope to obtain uniform controllability as $\eps \to 0^+$.

After this simple non-constructive result, we now provide with two explicit lower bounds on $\mathcal{C}_0(\TO, \eps)$ and an upper bound under stronger assumptions on $\f$ (namely parts of Assumption~\ref{assumptions}).
These three results are commented and compared in Section~\ref{s:rem-comm} below. 

Our first explicit lower bound on  the control cost $\mathcal{C}_0(\TO, \eps)$ is as follows. 
\begin{theorem}
\label{t:estim-Cobs-expo-1D}
Assume that $\mathfrak{a} \in C^1([0,L])$ (i.e. $\f \in C^2([0,L])$), and that Item~\ref{A2}  in Assumption~\ref{assumptions} is satisfied. 
Then, for all $E \in V(\M)$ and all $\delta >0$, there is $\eps_0>0$ such that we have for all $\eps<\eps_0$
$$
\mathcal{C}_0(\TO, \eps) \geq   \exp \frac{1}{\eps} \left(W_E(0) - \min_{\M} W_E  - E \TO - \delta  \right)  , \quad W_E = \frac{\f}{2} + d_{A,E} .
$$
In particular, we have
\bna
T_{\unif}(\{0\})\geq \sup_{E\in V(\M)}\frac{1}{E} \left[W_E (0)- \min_{\M} W_E  \right]  .
\ena
\end{theorem}

Theorem~\ref{t:estim-Cobs-expo-1D} states a result based on single energy levels (called $E$ in the statement). Our next result provides with a lower bound containing a nonlocal quantity defined on the ``limit spectrum'', namely $T_{E,B}$ in~\eqref{TEB} defined on $[V(\xzero),+\infty)$.
As shown by the proof, this term may be interpreted as an interaction term between the different energy levels.
Recall that $\Phi$ is defined in~\eqref{e:def-phi} and set 
$$
\widetilde{W_E}(s)=\frac{\f(s)}{2}-d_{A,E}(s) , \quad \text{ for } s \in [0,L] 
$$
(which is to be compared with $W_E$ defined in~\eqref{e:def-WE}).
\begin{theorem}
\label{thmlower+-}
Assume that $\mathfrak{a} \in C^\infty([0,L])$ (i.e. $\f \in C^\infty([0,L])$),  that Items~\ref{A1}--\ref{A4} in Assumption~\ref{assumptions} are satisfied, and that $\mathfrak{b}=\frac{\mathfrak{a}'}{2}$ (i.e. $\q= \frac{\f''}{2}$).
For any $\delta>0$, $E\in V(\M)$, $B\geq 0$, $C>0$ and for $0<\e<\e_0(\delta,B)$ and $0<T<C$, we have 
 \bna
\mathcal{C}_0(T,\e) \geq \exp \frac{1}{\eps} \left( W_E(0)-\sup_{[0,L]}\widetilde{W_E} -(E+B)T+T_{E,B} - \delta\right)
\ena
with
\begin{align}
\label{TEB}
T_{E,B}  
=\frac{1}{\pi}\int_{V(\xzero)}^{+\infty}\log \left|\frac{x+E+2B}{x-E}\right|\Phi'(x)dx. 
\end{align}
In particular, we have 
\bna
T_{\unif}(\{0\})\geq \sup_{E\in V(\M), B\geq 0} \frac{1}{E+B}\left(W_E(0)-\sup_{[0,L]}\widetilde{W_E}+T_{E,B} \right) .
\ena
\end{theorem}
A few remarks are in order
\begin{itemize}
\item  We prove in Section~\ref{s:section-allibert} that $\Phi$ is globally Lipschitz-continuous, and~\eqref{TEB} makes sense;
\item Note that $T_{E,B}$ increases if $B$ increases while $-BT$ decrease, so there might be an optimal choice of $B$ (hard to determine in general; see Section~\ref{s:explicit-comput} for explicit computations on specific examples).
\end{itemize}

\bigskip
Finally, our last result provides with an upper bound on $\mathcal{C}_0(\TO,\eps)$. Recall that  $T_1$ is defined in~\eqref{e:def-TT1} and $E_0$ in~\eqref{e:def-E0}.

\begin{theorem}
\label{thmcontrol1D}
Assume that $\mathfrak{a} \in C^\infty([0,L])$ (i.e. $\f \in C^\infty([0,L])$),  that Items~\ref{A1}--\ref{A4} in Assumption~\ref{assumptions} are satisfied, and that $\mathfrak{b}=\frac{\mathfrak{a}'}{2}$ (i.e. $\q= \frac{\f''}{2}$).
For all $\TO,\delta>0$, $m\in (0,1)$, there is $\eps_0 = \e_0(\TO,\delta,m)>0$ such that the control cost for the control problem \eqref{e:1Dcontrol-problemyn} satisfies for $0<\e<\e_0$
\begin{equation}
\label{e:cost-born-sup}
\mathcal{C}_0(\TO,\eps) \leq \exp \frac{1}{\eps} \left(\mathsf{G}(T,m,\delta)+\delta \right)  , 
\end{equation}
with 
$$
\mathsf{G}(T,m,\delta) = \frac{2 T_{1}^{2}}{ (1-m)\TO}+\sup_{E\in V([0,L])}\left[W_E(0)-\min_{[0,L]} W_E-m(1-\delta)E  \TO \right] .
$$
Moreover, if 
\begin{equation}
\label{e:T>pos-res}
T > 2\sqrt{2}\frac{T_1}{\sqrt{E_0}}+\frac{1}{E_0}\sup_{E\in V(\M)}\left[W_E(0)-\min_{[0,L]} W_E\right],
\end{equation}
then we have $\min_{m\in[0,1)}\mathsf{G}(T,m,\delta) <0$ for $\delta$ sufficiently small.
In particular, we have
\begin{align}
\label{e:Tunif-born-sup}
T_{\unif}(\{0\})\leq  \frac{1}{E_0} \left(  \sup_{E\in V(\M)}\left[W_E(0)-\min_{[0,L]} W_E\right] + 2\sqrt{2}\sqrt{E_0}T_1 \right). 
\end{align}
\end{theorem}
Notice that~\eqref{e:Tunif-born-sup} is a consequence of~\eqref{e:T>pos-res}, and that, in case~\eqref{e:T>pos-res} holds, the optimal control converges exponentially to zero.

\bigskip
Note that it may seem surprising that the sign of the vector field $\mathfrak{a} = \f'$ does not appear explicitly in the statements of Theorems~\ref{t:estim-Cobs-expo-1D}--\ref{thmlower+-}--\ref{thmcontrol1D}. It does play a role (as stressed by Lemma~\ref{l:WE} below) and some of the quantities involved in the statements above actually simplify in case $\mathfrak{a} >0$ or $\mathfrak{a} <0$.
\begin{corollary}
Under the assumptions of Theorems~\ref{t:estim-Cobs-expo-1D}--\ref{thmlower+-}--\ref{thmcontrol1D}, we have 
\begin{itemize}
\item either $\mathfrak{a}>0$ ($\f$ is increasing) on $[0,L]$, and 
\begin{align*}
T_{\unif}(\{0\})& \geq \sup_{E\in V(\M), B\geq 0} \frac{1}{E+B}\Big(\underbrace{d_{A, E}(0)+d_{A, E}(L)+\frac{\f(0)-\f(L)}{2}}_{\leq 0}+T_{E,B} \Big) , \\
T_{\unif}(\{0\})& \leq  \frac{4\sqrt{2}}{\mathfrak{a}(\xzero)}  T_1 ,
\end{align*}
\item or $\mathfrak{a}<0$ ($\f$ is decreasing) on $[0,L]$, and
\begin{align*}
T_{\unif}(\{0\})& \geq \sup_{E\in V(\M)}\frac{1}{E} \left[ W_E(0)-W_E(L) \right]  , \\
T_{\unif}(\{0\})& \geq \sup_{E\in V(\M), B\geq 0} \frac{1}{E+B}\left(2d_{A, E}(0)+T_{E,B} \right) ,\\ 
T_{\unif}(\{0\})& \leq  \frac{1}{E_0} \left(  \sup_{E\in V(\M)}\left[ W_E(0)-W_E(L)\right] + 2\sqrt{2}\sqrt{E_0}T_1 \right) .
\end{align*}
\end{itemize}
\end{corollary}
This result is a direct consequence of Theorems~\ref{t:estim-Cobs-expo-1D}--\ref{thmlower+-}--\ref{thmcontrol1D} combined with Lemma~\ref{l:WE} below.
We remark on the one hand that in case $\mathfrak{a}>0$ (in which the limit equation is a proper control problem), then the lower bound of Theorem~\ref{t:estim-Cobs-expo-1D} is trivial and the upper bound in Theorem~\ref{thmcontrol1D} only involves the spectral gap quantity $T_1$.
On the other hand, we notice that if $\mathfrak{a}<0$ (in which case the limit equation is not a proper control problem), geometric quantities involving the Agmon distance enter into play. As already remarked in~\cite{LL:18vanish} this allows, in this situation, to have $T_{\unif}$ very large compared to the  minimal flushing time of the limit problem $\eps=0$.
See also Section~\ref{s:explicit-comput} below for explicit computations on an example.
This is consistent with the results in~\cite{CG:05}, in which the sign of the vector field is of key importance.

\subsection{Remarks and comments}
\label{s:rem-comm}

\subsubsection{Remarks about the proofs}

The first step of our proofs consists in conjugating in Section~\ref{sectconjug} the control/observation equations by the weight $e^{\frac{\f}{2\eps}}$. Taking advantage of the fact that, in dimension one, every vector field is a gradient, this reformulates the (seemingly non-selfadjoint) transport equation with vanishing viscosity as a semiclassical heat equation $\eps \d_t w- P_\eps w=0$, involving the following semiclassical Schr\"odinger operator
 \begin{equation}
 \label{e:def-Peps-2}
 P_\eps := - \eps^2 \d_x^2 + \frac{|\f'|^2}{4} + \eps  \qf = - \eps^2 \d_x^2 + V + \eps  \qf .  
 \end{equation}
 see~\eqref{e:def-Peps}.
  All results of the article then rely on fine spectral properties of the operator $P_\eps$, that is to say 
 \begin{itemize}
 \item a precise knowledge of the spatial {\em  localization of eigenfunctions} of $P_\eps$; this is the object of the companion paper~\cite{LL:22} (see also  Section \ref{sect1D} where the results are recalled). 
 Roughly speaking, we use that a solution to $P_\eps \psi = E \psi$ behaves like $|\psi(x)|\sim e^{-\frac{d_{A,E}(x)}{\eps}}$ (up to some loss $e^{\frac{\delta}{\eps}}$) in the sense of $L^2$ density. Here  $d_{A,E}$ is the Agmon distance for the potential $V$ at energy $E$ defined in~\eqref{defAgmonbis}.
 \item a precise knowledge of the {\em distribution of eigenvalues} of $P_\eps$; and in particular the gap between (square roots of) two successive eigenvalues in the limit $\eps\to 0^+$. We extract the results we need from the article of Allibert~\cite{Allibert:98} (itself relying on~\cite{HS:84}) which provides with a precise asymptotics  (as $\eps\to 0^+$, as a function of the energy level $E$) of the distribution and the gap, see Appendix~\ref{s:Allibert}.
 \end{itemize}
 Note that properties of the classical Hamiltonian $p$ defined in~\eqref{e:def-p} (which is the principal symbol of $P_\eps$) naturally arise in the description of spectral properties of the quantum Hamiltonian $P_\eps$ the semiclassical limit $\eps\to 0^+$.
 This is why the functions $\Phi$ and $T$ defined in~\eqref{e:def-phi}--\eqref{e:def-TT1} enter into play (see also Section~\ref{s:section-allibert}).

The proof of Theorem \ref{t:estim-Cobs-expo-1D} is rather direct once the results of localization of eigenfunctions are obtained in Section \ref{sect1D}. Indeed, we only test the observability estimate on solutions of the semiclassical heat equation issued from eigenfunctions of $P_\eps$.

The proof of Theorem~\ref{thmlower+-} follows the spirit of Coron-Guerrero \cite{CG:05} and thus uses interactions between eigenfunctions. It seems more precise than the lower bound of Theorem \ref{t:estim-Cobs-expo-1D} which only considers a single eigenfunction. Indeed, the final estimate contains one part of harmonic analysis related to the spectrum and one part more geometric related to the concentration of eigenfunctions. Unfortunately, the geometric part related to the concentration of eigenfunctions seems less precise than that in Theorem \ref{t:estim-Cobs-expo-1D}. This is why we have chosen to keep both results.
Here, we use the spectral gap estimates of Allibert~\cite{Allibert:98}, which require $q_\f = constant$.

The proof of Theorem \ref{thmcontrol1D} uses the moment method which is classical  for 1D control problems. Yet, in this context, we need precise information about the localization of both the spectrum and the eigenfunctions. We thus rely on both items above. 
Moreover, in order to obtain estimates uniform in $\e$, we need to have a ``quantitative'' moment method, that is with explicit constants, at least uniform in $\e$. This is obtained in Proposition \ref{propgapheat} that provides a result of moment for heat type equation with an assumption on the spectral gap. The main advantage of this construction, which is of independent interest, is that we can follow (almost) explicitly the constants with respect to the parameters, which will be crucial to have estimates uniform in $\e$. The proof relies on Ingham estimates and a transmutation method of \cite{EZ:11}.

\subsubsection{Reformulation of the problem with a constant vectorfield}
We here notice that the control problem~\eqref{e:1Dcontrol-problemyn} can be reformulated as a problem with a constant vectorfield $\pm \d_x$, but with a varying viscosity. This uses Item~\ref{A1} of Assumption~\eqref{assumptions}, stating that the vectorfield is nondegenerate.
To straighten the vectorfield $\f'(x)\d_x$, we introduce its flow $Y_x(s)$ defined by 
$$
\d_s Y_{s_0}(s) = \f'(Y_{s_0}(s)), \quad Y_{s_0}(0) = s_0  .
$$
We also introduce its inverse $J_x(y) = \int_x^y \frac{ds}{\f'(s)}$, so that we have 
$$
J_x(Y_x(s)) = s  . 
$$
Denoting $v(x)=u(Y_0(x))$, we have $\d_xv(x)=\f'(Y_0(x)) (\d_yu)(Y_0(x))$.
Equivalently, we have $v(J_0(y))=u(y)$ and thus
$$
\d_y u(y) = J_0'(y)(\d_xv)(J_0(y)) = J_0'(Y_0(x)) (\d_xv)(x)  = \frac{1}{\f'(Y_0(x))} (\d_xv)(x), \quad x=J_0(y) . 
$$
As a consequence, if $u$ satisfies the equation
$$
(\d_t + \f'(y) \d_y +\f''(y) - \q(y)  - \eps\d_y^2) u(t,y) = 0 ,
$$
then $v(x)=u(Y_0(x))$ satisfies the equation
$$
\left(\d_t + \d_x +\f''(Y_0(x)) - \q(Y_0(x))  - \eps   \frac{1}{\f'(Y_0(x))}\d_x  \frac{1}{\f'(Y_0(x))} \d_x \right) v(t,x) = 0.
$$
The latter is a linear transport equation by the constant coefficient vector field $\d_x$ with a variable coefficient viscous perturbation operator $\eps \frac{1}{\f'(Y_0(x))}\d_x  \frac{1}{\f'(Y_0(x))} \d_x$.
Theorems~\ref{t:estim-Cobs-expo-1D}--\ref{thmlower+-}--\ref{thmcontrol1D} can all be translated as estimates on the minimal  time of uniform null-controllability in this context.

Note that the fact that, given a fixed vectorfield, the choice of the viscous perturbation changes the minimal uniform control time $T_{\unif}$ was already observed in~\cite{LL:18vanish}. 
\subsubsection{Remarks about the assumptions}

The assumptions we make on the vector field $\mathfrak{a}$ (resp. on $\f$ or on $V$) are issued from the analysis of the limit problem and from the two tools we use here as a black box in the analysis (namely localization of eigenfunctions~\cite{LL:22}, and spectral gap estimates~\cite{Allibert:98}).

Item \ref{A1} of Assumption~\ref{assumptions} is necessary for the transport equation with $\eps=0$ to be controlable (see Section \ref{s:limit-eq}) and is therefore quite natural.

Then, the essential assumption in both references~\cite{LL:22,Allibert:98} is Assumption~\ref{assumptions} Item~\ref{A2}, namely that the potential forms a single well.
Removing this would be an interesting problem, but would require a careful study of the interaction between the different wells and the tunneling effect, see~\cite{HS:84}. This is however beyond the scope of the present article.

The remaining assumptions: $\f \in C^\infty([0,L])$, Items~\ref{A3} and \ref{A4} in Assumption~\ref{assumptions} and $\mathfrak{b}=\frac{\mathfrak{a}'}{2}$ (i.e. $\q= \frac{\f''}{2}$, which amounts to $\qf=0$ in~\eqref{e:def-Peps-2}) come from the paper of Allibert \cite{Allibert:98}. The assumption $\mathfrak{b}=\frac{\mathfrak{a}'}{2}$ (i.e. $\q= \frac{\f''}{2}$ is technical and we believe that it can be removed (this would however require to reprove most of the spectral gap estimates in~\cite{Allibert:98} with an additional lower order term).
Item \ref{A4} of Assumption~\ref{assumptions} concerns the non degeneracy of the minimum of the potential. It could probably be weakened (as long as the potential is not too ''flat'' at the minimum), at the cost of several complications in the proofs (because of the associated degeneracy of the spectral gap near the potential minimum).

\subsubsection{Interpretation of the spectral quantities}
\label{s:section-allibert}
Here, we provide with some comments on the classical/spectral quantities $\Phi(E)$ and $T(E)$ (defined in~\eqref{e:def-phi} and~\eqref{e:def-TT1} respectively) entering into play in the above results, under Assumption~\ref{assumptions} Item~\ref{A2}. They are linked with properties of the classical Hamiltonian $p$ defined in~\eqref{e:def-p}.
First notice that $\Phi(E)$ is related to the following phase-space volume of the set $\{p\leq E\}$, that is 
\begin{align*}
\Phi(E) &= \Vol \left( \left\{(x,\xi) \in [0,L]\times \R , \quad V(x)\leq E \textnormal{ and } 0\leq \xi \leq \sqrt{E-V(x)} \right\}  \right) \\
& =  \frac12 \Vol \left( \left\{(x,\xi) \in [0,L]\times \R  ,\quad  p(x,\xi)\leq E \right\}  \right).
\end{align*}
As such, it is linked via the Weyl law to the asymptotic number of eigenvalues of $P_\eps = -\eps^2\d_x^2+V(x)+ \O(\eps)$ in the semiclassical limit $\eps \to 0^+$ as 
$$
\sharp\{\lambda_k^\eps\in \Sp(P_\eps) , \lambda_k^\eps \leq E \} = \frac{1}{2\pi \eps}  \left( \Vol(\{p\leq E\})  + o_E(1)\right)=  \frac{1}{\pi \eps}  \left( \Phi(E) + o_E(1)\right)
$$
(see e.g.~\cite{DS:book} in the boundaryless case or Theorem~\ref{t:Peps-spect-Allib} in the present setting).
Notice also that, for $\EE \geq \max V = \max \{ V(0), V(L) \}$, we have $\Phi(\EE) = \int_{0}^{L} \sqrt{\EE - V(s)} ds$ and in particular $\Phi(\EE) \sim L \sqrt{\EE}$ as $\EE \to + \infty$.
We use a more precise version of this formula (due to~\cite{HR:84,Allibert:98}) stating that 
\begin{align}
\label{e:Phi-lambda-count}
\lambda_k^\eps \approx \Phi^{-1}(\eps \pi k) , \text{ uniformly in both } \eps \to 0^+, k\in \N ,
\end{align}
where the meaning of $\approx$ is made precise in Theorem~\ref{t:Peps-spect-Allib}.

\medskip
Concerning the quantity $T(\EE)$ in~\eqref{e:def-TT1}, we now explain how it is linked to the period of trajectories of the Hamiltonian vector field $H_p$ associated to the Hamiltonian $p$ (defined in~\eqref{e:def-p}), in the energy level $p(x,\xi)=\EE$. More precisely, the Hamiltonian flow of $p(x,\xi)=\xi^2+V(x)$ is defined by $\dot{x}(s)=2\xi(s) , \dot{\xi}(s) = -V'(x(s))$, and the Hamiltonian $p(x,\xi)$ is preserved by the flow. Hence, under Assumption~\ref{assumptions} Items~\ref{A1}--\ref{A2}, if a curves has $p = E$, then in any time interval such that  $\xi(t)> 0$ and $x_-(E) \leq x(0) \leq x(T) \leq x_+(E)$, we have
$$
T = \int_0^T dt = \int_{x(0)}^{x(T)}\frac{dx}{2\xi(x)} =  \int_{x(0)}^{x(T)}\frac{dx}{2\sqrt{E-V(x)}}.
$$
Hence, in the energy level $\{p = E\}$, the Hamiltonian flow of $p$ (consists in two different trajectories and) is periodic with period 
$$
\mathcal{T}(E) = 2\int_{x_-(E)}^{x_+(E)}\frac{dx}{2\sqrt{E-V(x)}} = \int_{x_-(E)}^{x_+(E)}\frac{dx}{\sqrt{E-V(x)}} .
$$
As a consequence, we deduce that the quantity $T(E)$ (defined in~\eqref{e:def-TT1}) verifies
\begin{align}
\label{e:period-TT}
T(E) = 2 \sqrt{E} \mathcal{T}(E).
\end{align}
Note that for large energies, $\mathcal{T}(E) \sim_{E\to +\infty} \frac{L}{\sqrt{E}}$ and $T(E) \sim_{E\to +\infty}2 L$.
\begin{lemma}
\label{l:phiT}
The quantities $\Phi(\EE)$ and $T(\EE)$ defined in~\eqref{e:def-phi} and~\eqref{e:def-TT1} are linked by the following: $\Phi \in W^{1,\infty}([E_0,\infty))$ and 
$$
\Phi'(\EE) = \frac{1}{4\sqrt{\EE}} T(\EE) = \frac12 \mathcal{T}(E)  ,\qquad  \text{i.e.}\quad \big(\Phi(\EE^2) \big)' =  \frac{1}{2} T(\EE^2) .
$$
\end{lemma}
This lemma is proved in Appendix~\ref{A:technical}.
 
\bigskip
We now give a spectral interpretation of $T(E)$ and $\mathcal{T}(E)$, explaining how these classical quantities enter into the description of spectral properties of $P_\eps$. Precise statements and proofs are provided in Appendix~\ref{s:Allibert}, based on results obtained by Allibert in~\cite{Allibert:98} (themselves relying on~\cite{HS:84,HR:84}).

Denoting $N(\beta) :=\Phi(\beta^2)$, we have according to~\eqref{e:Phi-lambda-count} that 
$\Phi(\lambda_k^\eps)\approx \e \pi k$. 
Writing $\beta_k^\eps=\sqrt{\lambda_k^\eps}$ for the square root of the eigenvalues, we thus have $\e \pi k\approx \Phi(\lambda_k^\eps)=\Phi((\beta_k^\eps)^2)=N(\beta_k^\eps)$. Hence, denoting by 
\bna
\mathbf{G}_k^\eps := \beta_{k+1}- \beta_k = \sqrt{\lambda_{k+1}^\eps} - \sqrt{\lambda_k^\eps}
\ena
 the local spectral gap for the square roots of eigenvalues (spacing of square roots of eigenvalues), we obtain 
$$
\e\pi(k+1) - \e \pi k \approx N(\beta_{k+1}^\eps)- N(\beta_k^\eps) \approx \mathbf{G}_k^\eps N'(\beta_k^\eps).
$$
We finally obtain
\bna
\mathbf{G}_k^\eps  \approx \frac{\e\pi}{N'(\beta_k^\eps)} = \frac{2 \e\pi}{T(\lambda_k^\eps)} . 
\ena
As a consequence, the quantity $T(E)$ measures the local spectral gap for the square roots of eigenvalues at energy $E$, and hence $T_1 = \sup_{E\geq E_0}T(E)$ yields a uniform lower bound for the spectral gap for the square roots of eigenvalues.
This is actually stronger than a uniform lower bound for the spectral gap of eigenvalues themselves, for 
$$
\lambda_{k+1}^\eps - \lambda_k^\eps = (\beta_{k+1}^\eps+\beta_k^\eps)(\beta_{k+1}^\eps- \beta_k^\eps)
\approx 2 \beta_k^\eps \mathbf{G}_k^\eps \approx 2 \e\pi  \frac{2\sqrt{\lambda_k^\eps}}{T(\lambda_k^\eps)}  =\frac{2 \e\pi}{\mathcal{T}(\lambda_k^\eps)} ,
$$
where we have used~\eqref{e:period-TT} in the last equality.

\subsubsection{Comparing the minimal times appearing in Theorems~\ref{t:estim-Cobs-expo-1D},~\ref{thmlower+-}, and~\ref{thmcontrol1D}}

In the estimates we obtain on $T_{\unif}$ (we write $T_{\unif} = T_{\unif}(\{0\})$ in this section for short) in Theorems~\ref{t:estim-Cobs-expo-1D},~\ref{thmlower+-} and~\ref{thmcontrol1D}, two parameters enter into play:
\begin{itemize}
\item first a ``Spectral'' parameter, related to the localization of the spectrum, or, more precisely, to the size of the spectral gap at energy $E$;
\item and second a ``Geometrical'' parameter, related to the localization of the eigenfunctions at energy $E$.
\end{itemize} 
In order to compare these parameters, we are led to define the following spectral/geometric constants (the index of the constant refers to the theorem where it appears)
\begin{equation}
\label{e:def-ctes}
\begin{array}{lcl}
S_{\ref{t:estim-Cobs-expo-1D}}=0 , &&  G_{\ref{t:estim-Cobs-expo-1D},E}=W_E(0) - \min_{\M} W_E ,  \\
S_{\ref{thmlower+-},E,B}=T_{E,B} , &  & G_{\ref{thmlower+-},E}=W_E(0)-\sup_{[0,L]}\widetilde{W_E},   \\
S_{\ref{thmcontrol1D}} = 2\sqrt{2}\sqrt{E_0}T_1 , & &G_{\ref{thmcontrol1D},E} =W_E(0)-\min_{[0,L]} W_E =G_{\ref{t:estim-Cobs-expo-1D},E}, 
 \end{array}
 \end{equation}
where, recall, $E_0 = V(\xzero) = \frac{|\f'(\xzero)|^2}{4} = \min_{[0,L]}V >0$.
With these definitions in hand, the critical times appearing in Theorems~\ref{t:estim-Cobs-expo-1D},~\ref{thmlower+-} and~\ref{thmcontrol1D}   respectively are 
 \begin{align}
T_{\ref{t:estim-Cobs-expo-1D}}& =\sup_{E\in V(\M)} \frac{1}{E}  \left( G_{\ref{t:estim-Cobs-expo-1D},E} +S_{\ref{t:estim-Cobs-expo-1D}} \right) =\sup_{E\in V(\M)} \frac{1}{E}G_{\ref{t:estim-Cobs-expo-1D},E}  , \label{e:def-T14} \\  
 T_{\ref{thmlower+-}}&= \sup_{E\in V(\M), B\geq 0} \frac{1}{E+B}\left( G_{\ref{thmlower+-},E}+S_{\ref{thmlower+-},E,B} \right) , \label{e:def-T15} \\
 T_{\ref{thmcontrol1D}} &=\frac{1}{E_{0}}\left(\sup_{E\in V(\M)}G_{\ref{thmcontrol1D},E} + S_{\ref{thmcontrol1D}} \right) 
= \sup_{E\in V(\M)} \frac{1}{E_{0}}\left( G_{\ref{thmcontrol1D},E} + S_{\ref{thmcontrol1D}} \right) , \nonumber 
 \end{align}
 and the associated result formulates (sometimes assuming $\q= \frac{\f''}{2}$) as 
$$
T_{\ref{t:estim-Cobs-expo-1D}}\leq T_{\unif} , \qquad  T_{\ref{thmlower+-}}\leq T_{\unif} ,\qquad  T_{\unif} \leq T_{\ref{thmcontrol1D}} .
$$
We now try to compare the different quantities involved in~\eqref{e:def-ctes}.
We first need to compare $T_{E,B}$ and $T_1$.

\begin{lemma}
\label{l:TEB-T}
The quantities $T_{E,B}$ and $T_1$ are linked by  
\begin{align*}
T_{E,B}  \leq T_1 \frac{\sqrt{E_0}}{2\pi} \Gamma_0\left(\sqrt{\frac{E+2B}{E_0}} , \sqrt{\frac{E}{E_0}}\right),
\end{align*}
where, for $\alpha \geq 1, \beta>1$,
\begin{align*}
 \Gamma_0(\alpha, \beta)   & = \int_{1}^{+\infty}\log \left|\frac{y^2 +\alpha^2}{y^2-\beta^2}\right|  dy \\
& =\pi \alpha - \log(1+\alpha^2)  - 2 \alpha \arctan \left( \frac{1}{\alpha}\right) 
+ \log(\beta^2-1)  + \beta \log \left( \frac{\beta+1}{\beta-1}\right) , 
\end{align*}
extends as a continuous function on $\{(\alpha,\beta) \in \R^2 , 1\leq \beta\leq \alpha\}$.
\end{lemma}
This lemma is proved in Appendix~\ref{A:technical}

\begin{lemma}
\label{l:compare}
The above quantities are linked by 
 \bnan
 \label{e:compG}
 \text{for all }E\geq E_0 , \quad 
G_{\ref{t:estim-Cobs-expo-1D},E} = G_{\ref{thmcontrol1D},E} \geq  G_{\ref{thmlower+-},E} , \\
 \text{for all }E\geq E_0, B\geq 0 , \quad 
\frac{1}{E_0} G_{\ref{thmcontrol1D},E} \geq \frac{1}{E}G_{\ref{t:estim-Cobs-expo-1D},E} \geq \frac{1}{E+B} G_{\ref{thmlower+-},E} , \nonumber
 \enan

 and for all $E\geq E_0$, $B\geq0$,
 \begin{align}
 \label{e:linky}
0= S_{\ref{t:estim-Cobs-expo-1D}} \leq  \frac{S_{\ref{thmlower+-},E,B}}{E+B} \leq \kappa_0  \frac{S_{\ref{thmcontrol1D}}}{E_0} , \quad \text{ with } \quad \kappa_0 = \frac{1}{2\pi \sqrt{2}} \max_{\alpha\geq \beta\geq1}  \frac{\Gamma_0(\alpha,\beta)}{\alpha^2+\beta^2} <+\infty .
\end{align}
\end{lemma}
This lemma is proved in Appendix~\ref{A:technical}.
From this lemma, we draw the following conclusions:
\begin{itemize}
\item The geometric quantity $G_{\ref{thmcontrol1D},E} =G_{\ref{t:estim-Cobs-expo-1D},E}=W_E(0)-\min_{[0,L]} W_E$ seems to be the appropriate one to describe the localization of eigenfunctions. In this particular 1D situation, we indeed know very precisely where eigenfunctions are localized, see Section~\ref{sect1D} below or~\cite{LL:22}. Theorems~\ref{thmcontrol1D} and~\ref{t:estim-Cobs-expo-1D} are thus accurate in this respect, whereas Theorem~\ref{thmlower+-} is not. Note that the quantity $\frac{1}{E_{0}} G_{\ref{thmcontrol1D},E}$ instead of $\frac{1}{E} G_{\ref{thmcontrol1D},E}$ makes however Theorem~\ref{thmcontrol1D} less accurate than Theorem~\ref{t:estim-Cobs-expo-1D}. This can be summarized as 
$$
\text{Theorem~\ref{t:estim-Cobs-expo-1D}} \gg \text{Theorem~\ref{thmcontrol1D}} \gg \text{Theorem~\ref{thmlower+-}} ,
$$
where $\gg$ stands for ``more accurate as far as the geometric quantity is involved''.
Note however that if one rather compares 
$$
\sup_{E\in V(\M), B\geq 0} \frac{1}{E+B} G_{\ref{thmlower+-},E} = \sup_{E\in V(\M)} \frac{1}{E}  G_{\ref{thmlower+-},E} , \quad \text{  with  } \quad 
\frac{1}{E_{0}}  \sup_{E\in V(\M)} G_{\ref{thmcontrol1D},E} ,
$$
then Theorems~\ref{thmcontrol1D} and Theorem~\ref{thmlower+-} are no longer comparable.
\item Now, as far as the spectral quantity is concerned, Theorem~\ref{t:estim-Cobs-expo-1D} does not say anything, Theorem~\ref{thmlower+-} yields a seemingly fine lower bound (comparable to that obtained in~\cite{CG:05}), whereas Theorem~\ref{thmcontrol1D} seems to provide with a relatively rough upper bound. This can be summarized as 
$$
\text{Theorem~\ref{thmlower+-}} \gg \text{Theorem~\ref{thmcontrol1D}} \gg \text{Theorem~\ref{t:estim-Cobs-expo-1D}} ,
$$
where $\gg$ stands for ``more accurate as far as the spectral quantity is involved''.
\end{itemize}
 In particular, the lower bound of Theorem \ref{t:estim-Cobs-expo-1D} is better than the one of Theorem \ref{thmlower+-} from the geometrical point of view, while the latter is better from the spectral point of view.

\medskip
The following lemma allows to better understand the importance of the direction of the vector field $\f'$, i.e. to distinguish properties of $\f'>0$ from $\f'<0$ (recall that the asymmetry comes from the fact that the control acts only on left boundary).
\begin{lemma}
\label{l:WE}
Assume that Items~\ref{A1} and~\ref{A2} in Assumption~\ref{assumptions} are satisfied. Then one of the following two statements hold:
\begin{itemize}
\item either $\f$ is increasing: then for any $E\geq E_{0}$, the functions $x\mapsto W_{E}(x)$ and $\widetilde{W_E}(x)$ are increasing, and the constants defined in~\eqref{e:def-ctes} satisfy
\begin{align*}
G_{\ref{t:estim-Cobs-expo-1D},E}&=0,  \quad \text{ and }\\
 G_{\ref{thmlower+-},E}&=d_{A, E}(0)+d_{A, E}(L)+\frac{\f(0)-\f(L)}{2}\leq 0 ;
 \end{align*}
 \item or $\f$ is decreasing: then for any $E\geq E_{0}$, the functions $x\mapsto W_{E}(x)$ and $\widetilde{W_E}(x)$ are decreasing, and the constants defined in~\eqref{e:def-ctes} satisfy
 \begin{align*}
G_{\ref{t:estim-Cobs-expo-1D},E}=& W_E(0)-W_E(L)\geq 0 ,\\
 G_{\ref{thmlower+-},E}=&2d_{A, E}(0)\geq 0 .
 \end{align*}
 \end{itemize}
In both cases,  $E\mapsto G_{\ref{thmlower+-},E}$ is a nonincreasing function. 

If we assume additionally that $\f$ is an odd function with respect to $\xzero=L/2$, that is to say, $\f(L/2+x)=-\f(L/2-x)$ for all $x \in [0,L/2]$, then we have the following simplifications:
if $\f$ is increasing, then 
\begin{align*}
G_{\ref{t:estim-Cobs-expo-1D},E}=0, \quad \text{ and } \quad 
 G_{\ref{thmlower+-},E} =2d_{A, E}(L)-\f(L)=2d_{A, E}(0)+\f(0) ,
\end{align*}
while if $\f$ is decreasing, we have
\bna
G_{\ref{t:estim-Cobs-expo-1D},E}= \f(0) = -\f(L) , \quad \text{ and } \quad 
 G_{\ref{thmlower+-},E}=2d_{A, E}(L) =2d_{A, E}(0) .
\ena
In particular, $G_{\ref{t:estim-Cobs-expo-1D},E}$ is independent on $E$ in both cases.
\end{lemma}
This lemma is proved in Appendix~\ref{A:technical}. It is also very useful to compute the value of the different constants on explicit examples.

Our results lead us to conjecture that, 
under Assumption~\ref{assumptions} Items~\ref{A1}--\ref{A2} and~\ref{A4}, there is a distribution kernel $\mathsf{K}(x,E)$ such that 
$$
T_{\unif} = \sup_{E\in [E_0,+\infty)} \frac{1}{E}  \left( G_{\ref{t:estim-Cobs-expo-1D},E} +S_{E} \right) , \quad \text{ with } S_E = \int_{E_0}^\infty\mathsf{K}(x,E) \Phi'(x) dx = \frac12\int_{E_0}^\infty\mathsf{K}(x,E)\mathcal{T}(x) dx . 
$$
However, we do not have a precise idea of what the kernel $\mathsf{K}$ should be, but $\mathsf{K}(x,E) =\log \left|\frac{x+E+2B}{x-E}\right|$ would look to be a good candidate for some $B$.

\subsubsection{Explicit computations on an example}
In section~\ref{s:explicit-comput} below, we compute explicitly and further compare all upper/lower bounds for the functions
\bna
 \f_{M,a}^{\pm}(x)=\pm\int_0^x\sqrt{a^2s^2+M^{2}}ds  , \quad \text{ that is to say} \quad  \mathfrak{a}^{\pm}(x) = (\f_{M,a}^{\pm})'(x) = \pm \sqrt{a^2x^2+M^{2}} ,
\ena
defined on the shifted interval $(-L/2,L/2)$ (instead of $(0,L)$).
The latter are associated to the harmonic potential $V(x) = \frac{|\f^{\pm '}_{M,a}(x)|^2}{4}=\frac{a^2x^2+M^{2}}{4}$, and our results apply for $M > 0$ and $a>0$.
For $a=0$ (to which our results do not apply), the vector fields correspond to the case studied in~\cite{CG:05,Gla:10,Lissy:12,Lissy:14,Lissy:15,Darde:17,YM:19,YM:19bis}.
For large values of $a$, the potential is very convex and far from the situation $a=0$.
We draw in particular the following consequences:
\begin{itemize}
\item in case $-$ (that is, for $\f_{M,a}^{-}$), then $T_{\f^{-'}_{M,a}} (\{-L/2\}) \underset{a\to +\infty}{\longrightarrow } 0^+$ (the flushing time associated to the limit equation $\eps=0$, see Section~\ref{s:limit-eq}) whereas $T_{\unif}  (\{-L/2\}) \underset{a\to +\infty}{\longrightarrow } + \infty.$
In particular, we recover~\cite[Section~3.3]{LL:18vanish}, stating that  $\frac{T_{\unif}(\{-L/2\})}{T_{\f'^{- }_{M,a}} (\{-L/2\})}   \underset{a\to +\infty}{\longrightarrow } + \infty$.
We obtain actually the stronger statement that if $a \to 0^+$, the limit problem is controllable in a time $T_{\f^{\pm '}_{M,a}} (\{-L/2\})\to 0$ whereas uniform controllability holds for a time $T_{\unif}(\{-L/2\}) \to +\infty$. This is a refinement of~\cite[Section~3.3]{LL:18vanish}. See Section~\ref{s:a=infty}.

\item In the formal limit $a \to 0^+$, we obtain the lower bounds
\begin{align}
\liminf_{a \to 0^+}T_{\unif,a} \geq \liminf_{a \to 0^+}T_{\ref{thmlower+-},a}\geq &\frac{L}{M} ,\quad \text{ (Case $+$)} , \label{e:ato0+} \\
\liminf_{a \to 0^+}T_{\unif,a} \geq \liminf_{a \to 0^+}T_{\ref{thmlower+-},a}\geq&\frac{2\sqrt{2}L}{M},\quad \text{ (Case $-$)}  \label{e:ato0-}.
\end{align}
As a consequence, the formal limit $a \to 0^+$ coincides with the known lower bounds for the Coron-Guerrero problem $a=0$, appearing in the literature. The first one was obtained by Coron-Guerrero~\cite{CG:05} while the second was obtained by Lissy \cite[Theorem 1.3]{Lissy:15} (using a variant of method of~\cite{CG:05}). See Section~\ref{s:cascst}.

\item In the formal limit $a \to 0^+$, the upper bound of Theorem \ref{thmcontrol1D} degenerates since $T_1 \underset{a\to 0^{+}}{\rightarrow} +\infty$. This suggest that the quantity $T_1$ is not the appropriate one (at least in this regime). 
A variation of our approach however applies to the case $a=0$, but yields slightly less accurate constants than those available in the literature~\cite{CG:05,Gla:10,Lissy:12,Lissy:14,Lissy:15,Darde:17}, see Section~\ref{s:cascst} for a discussion.
\end{itemize}

\subsubsection{Comparison with the results in~\cite{LL:18vanish}}
The result in Theorem~\ref{t:estim-Cobs-expo-1D} is a one-dimensional refinement of~\cite[Theorems~1.5 and 3.1]{LL:18vanish}, which instead states (in a much more general setting of a compact manifold with boundary, with essentially no assumptions on $\f$ or $V$)
\begin{align*}
\mathcal{C}_0(\TO, \eps) \geq  \exp \frac{1}{\eps} \left( W_E(0) - \max_{K_E} W_E - E \TO - \delta \right)  , \quad  \text{for all }E \in V(\M), \delta >0 .
\end{align*}
and in particular
\begin{align*}
T_{\unif}(\{0\}) \geq  \sup_{E\in V(\M)} \frac{1}{E}\left( W_E(0) - \max_{K_E} W_E  \right) = \frac{1}{E_0}\left( W_{E_0}(0) - \frac{\f(\xzero)}{2}  \right)  , \quad E_0 = \min_{[0,L]}V = V(\xzero).
\end{align*}
This last equality is explained in~\cite[remark following Theorem~3.1]{LL:18vanish}.
In Theorem~\ref{t:estim-Cobs-expo-1D}, we are able to replace $-\max_{K_E}W_E$ by $-\min_{\M} W_E$. 
This improvement comes from two additional knowledge we have on the eigenfunctions of a conjugated operator $P_\eps$ (see~\eqref{e:def-Peps} below) in this very particular $1D$ single well problem (see Section~\ref{sect1D} below o~\cite{LL:22}): an eigenfunction $\psi_\eps$ associated to an eigenvalue $E_\eps$ of $P_\eps$, converging to $E$ as $\eps \to 0^+$:
\begin{itemize}
\item spreads over the whole classically allowed region $K_E$ (propagation estimates);
\item vanishes {\em at most} like $e^{-\frac{1}{\eps}d_{A,E}}$ in the classically forbidden region (Allibert estimates).
\end{itemize}
In higher dimension, the first result is false (and the issue of understanding the asymptotic distribution of the distribution $|\psi_\eps|^2dx$ is extremely intricate, even for a given energy $E$) in general; and to our knowledge, the second result does not seem to be well-understood. 

Along the proof of the present paper, we could state an analogue of Theorem~\ref{t:estim-Cobs-expo-1D} for the internal uniform controllability/observability question by an open set $\omega \subset [0,L]$. The latter problem is not considered in the main part of the paper but is the main focus in~\cite{LL:18vanish}. 

\subsubsection{Uniform controllability of the semiclassical heat equation}
As already mentioned, all results proved in Theorems~\ref{t:estim-Cobs-expo-1D},~\ref{thmlower+-} and~\ref{thmcontrol1D} may be reformulated in terms of uniform (resp. non-) observability/controllability results for the semiclassical heat equation
\begin{equation}
\label{e:semiclass-heat}
\left\{
\begin{array}{rl}
\left(\eps \d_t -\eps^2 \d_x^2 + V(x) \right) v = 0 ,& (t,x)\in (0,T) \times (0,L) , \\
v (t,0) = h(t) , \quad v(t,L) = 0 ,&  t \in (0,T) , \\
v(0,x) = v_0 (x),&  x \in (0,L) ,
\end{array}
\right.
\end{equation}
 in the semiclassical limit $\eps \to 0^+$ and in weighted $L^2$-spaces of type $e^{\frac{\f}{2\eps}}L^2(0,L) = L^2((0,L),e^{-\frac{\f}{\eps}}dx)$. Note that in that setting, we do not need that $\f$ and $V$ be linked one to the other (and then have to change the definitions of $W,\widetilde{W}$ accordingly).  We do not state these results for the sake of brevity.

\begin{remark}
The semiclassical heat equation~\eqref{e:semiclass-heat} can be rewritten as 
$$\left( \eps^{-1}\d_t - \d_x^2 + \eps^{-2}V(x) \right) v =0$$ on a fixed time interval $[0,T]$. Rescaling in time, this amounts to study 
\begin{equation}
\label{e:chelou-chaleur-large-pot}
\left( \d_t - \d_x^2 + \eps^{-2}V(x) \right) v =0,
\end{equation}
on a time interval $[0,\eps T]$, that is, the heat equation with a large potential in small time. If we are interested in the controllability of the same equation~\eqref{e:chelou-chaleur-large-pot} in fixed time (independent of $\eps$), the techniques described in the present paper (see in particular Section~\ref{s:proof-upper-bound}) allow to obtain uniform estimates as well, and recover for instance the results of~\cite[Proposition 1.5.]{BDE:20} (proved by different techniques, namely Carleman estimates). In that reference, it is used to control the Grushin equation.
More precisely, the techniques above imply the following Proposition (analogue of~\cite[Proposition 1.5.]{BDE:20}).
\begin{proposition}
\label{propcontrolzerointermbis}
Let $V\in C^{\infty}([0,L])$ satisfy Items~\ref{A1}--\ref{A4} in Assumption~\ref{assumptions}. 
Let $T>0$ and fix $\delta>0$. Then, there exists $\e_{0}$ and $C>0$ so that for any $0<\e<\e_{0}$ and $v_0\in L^{2}(0,L)$, there exists a control $h\in L^{2}(0,T)$ to zero of
\begin{equation*}
\left\{
\begin{array}{rl}
\left(\d_t -\d_x^2 + \frac{1}{\eps^2} V(x) \right) v = 0 ,& (t,x)\in (0,\tau) \times (0,L) , \\
v (t,0) = h(t) , \quad v(t,L) = 0 ,&  t \in (0,\tau) , \\
v(0,x) = v_0 (x),&  x \in (0,L) .
\end{array}
\right.
\end{equation*} with the control cost
\begin{align}
\label{estimz}
\nor{h}{L^2(0,T)}^2 \leq  C e^{\frac{d_{A}(0)+\delta}{\eps}}\nor{v_{0}}{L^{2}}^2 .
\end{align}
\end{proposition}
Note that the equation can also be rewritten as  $(\eps^2\d_t + P_\eps ) v=0$ (compare with the semiclassical heat equation~\eqref{e:semiclass-heat} where we have $\eps\d_t$).
\end{remark}

\section{General facts about transport equation and vanishing viscosity limit}
 
\subsection{Duality between boundary control and observation problems}

In the present one dimensional setting recall the control problem under consideration is~\eqref{e:1Dcontrol-problemyn} (and is written in a ``gradient field'' way, which is always possible in dimension one).
The associated (forward in time) observation problem is 
\begin{equation}
\label{e:1Dobs-problem}
\left\{
\begin{array}{rl}
(\d_t - \f'\d_x - \q  - \eps\d_x^2) u = 0 ,& (t,x)\in (0,T) \times (0,L) , \\
u (t,0)  = u(t,L) = 0 ,&  t \in (0,T) , \\
u(0,x) = u_0 (x),&  x \in (0,L) ,
\end{array}
\right.
\end{equation}
with $\f = \f(x)$ and $q=q(x)$.
The solution $y$ of the controlled equation~\eqref{e:1Dcontrol-problemyn} and the solution $u$ of free equation~\eqref{e:1Dobs-problem} are linked by the following duality equation:
 \begin{align}
 \label{e:duality-boundary-case}
\left( u(T) , y_0 \right)_{L^2(0,L)} - \left( u_0, y(T)\right)_{L^2(0,L)}  + \int_0^T  \eps \partial_x u(t,0) h(T-t) dt =0 .
\end{align}

 The boundary observability problem for~\eqref{e:1Dobs-problem} can be formulated as follows. 
 Does there exist a constant $C>0$ such that 
 \bnan
 \label{e:transport-viscous-obs-boundary}
C^2 \int_0^T  |\eps \d_x u(t,0)|^2  dt \geq \|u (T)\|_{L^2(0,L)}^2 , \quad \text{ for all } u_0 \in L^2(0,L)\text{ and } u \text{ solution of~\eqref{e:1Dobs-problem}}.
 \enan
We define accordingly
 $$
C_0(T,\eps) := \inf\{C\in\R^{+} \text{ such that \eqref{e:transport-viscous-obs-boundary} holds} \}.
  $$  
Classical duality arguments (see~\cite{DR:77} or~\cite[Chapter~2.3]{Cor:book}) yield the following statement.
 
 \begin{lemma}[Observability constant = control cost]
 \label{l:obs-cont}
 Given $(\eps, T)$, Equation~\eqref{e:1Dcontrol-problemyn} is null-controllable if and only if the observability inequality~\eqref{e:transport-viscous-obs-boundary} holds. 
 Moreover, we then have $\mathcal{C}_0(T,\eps)= C_0(T,\eps)$.
 \end{lemma}
As usual, this allows us to mainly focus on the observability inequality~\eqref{e:transport-viscous-obs-boundary}.
 
\subsection{Gradient flows, conjugation and reformulation}
\label{sectconjug}

As in~\cite{LL:18vanish}, we first proceed with the following conjugation:
$$
e^{-\frac{\f}{2\eps}} \d_x^2 e^{\frac{\f}{2\eps}} = \d_x^2 + \frac{1}{\eps} \f' \d_x + \frac{|\f'|^2}{4\eps^2} + \frac{\f''}{2\eps} .
$$
We denote by $\frac{1}{\eps^2}P_\eps := - \d_x^2 + \frac{|\f'|^2}{4\eps^2} + \frac{\f''}{2\eps}-\frac{\q}{\eps}$, that is to say
 \begin{equation}
 \label{e:def-Peps}
 P_\eps := - \eps^2 \d_x^2 + \frac{|\f'|^2}{4} + \eps  \qf  , \quad \text{ with } \quad \qf= \frac{\f''}{2}-\q .
 \end{equation}
The above computation implies that 
\bnan
\label{e:trans-eps}
e^{-\frac{\f}{2\eps}} \big( \frac{1}{\eps^2}P_\eps  \big)e^{\frac{\f}{2\eps}} = - \d_x^2 - \frac{1}{\eps} \f' \d_x-\frac{\q}{\eps} ,
\enan
the last operator being that appearing in the observation/free evolution problem~\eqref{e:1Dobs-problem} multiplied by $\e$.
The operator $P_\eps$ is selfadjoint in $L^2(0,L)$ endowed with domain $D(P_\eps) = H^2(0,L)\cap H^1_0(0,L)$. Hence, the operator $ - \d_x^2 - \frac{1}{\eps} \f' \d_x-\frac{\q}{\eps}$ is also formally selfadjoint, but in the weighted space $L^2((0,L), e^{\frac{\f}{\eps}} dx)$.
We reformulate the uniform observability problem~\eqref{e:1Dobs-problem} in terms of the heat equation involving the operator $P_\eps$ defined in~\eqref{e:def-Peps} (see~\cite[Lemma~2.9]{LL:18vanish}).

\begin{lemma}[Observation problem: equivalent reformulation]
\label{lemequiveunif}
Given $\TO, C_0, \eps >0$, the following statements are equivalent. 
\begin{enumerate}
\item The function $u$ solves 
  \begin{equation}
  \label{e:heat-transp-eps}
 \left\{
 \begin{aligned}
&(\d_t -  \f' \d_{x} u-\q - \eps \d_{x}^{2} ) u = 0 , &(t,x)\in (0,\TO)\times  (0,L) ,  \\
&u (t,0)  = u(t,L) = 0 , & t \in (0,\TO) , \\ 
 \end{aligned}
 \right.
 \end{equation}
\begin{align}
\label{e:obs-heat-transp-epsboundary}
\text{resp. }\quad \nor{u(\TO)}{L^2(0,L)}^2 \leq C_0^2\int_0^{\TO} |\eps\d_x u(t,0)|^2 dt .
\end{align}
\item The function $\zeta(t,x)=e^{\f(x)/2\e} u(t,x)$ solves
  \begin{equation}
  \label{e:heat-transp-eps-Witt}
 \left\{
 \begin{aligned}
&\eps \partial_t \zeta + P_\eps \zeta =0  , & (t,x)\in (0,\TO)\times  (0,L) ,  \\
&\zeta(t,0)  = \zeta(t,L) = 0 , & t \in (0,\TO) , \\ 
 \end{aligned}
 \right.
 \end{equation}
\begin{align}
\label{e:obs-heat-transp-eps-Wittboundary}
\text{resp. }\quad \nor{e^{-\frac{\f}{2\eps}} \zeta(\TO)}{L^2(0,L)}^2 \leq C_0^2\int_0^{\TO}  \left|e^{-\frac{\f(0)}{2\eps}} \eps \d_x \zeta(t,0)\right|^2 dt .
\end{align}
\end{enumerate}
\end{lemma}

\medskip
A similar conjugation result also holds in the controllability side, using conjugation with the opposite sign. More precisely, still with $P_\eps$ defined in~\eqref{e:def-Peps}, we now have (instead of \eqref{e:trans-eps})
\bna
e^{\frac{\f}{2\eps}} \big( \frac{1}{\eps^2}P_\eps  \big)e^{-\frac{\f}{2\eps}} = -\d_{x}^{2} + \frac{1}{\eps} \f' \d_{x}-\frac{\q}{\eps}+\frac{ \f''}{\eps} .
\ena
This time, the conjugation of $P_\eps$, which is selfadjoint on $L^2((0,L),dx)$ with domain $H^1_0\cap H^2(0,L)$, yields the operator $- \d_{x}^{2} + \frac{1}{\eps} \f'\d_{x}-\frac{\q}{\eps}+\frac{ \f''}{\eps}$ (which thus becomes formally selfadjoint in $L^2((0,L), e^{-\frac{\f}{\eps}}dx)$ with Dirichlet boundary conditions). We can then obtain a similar version of Lemma \ref{lemequiveunif} from the control point of view to relate the control problem~\eqref{e:1Dcontrol-problemyn} to a control problem with $P_{\e}$. 

\begin{lemma}[Control problem: equivalent reformulation]
\label{lemequivecontrol}
Given $T, \eps >0$, the following statements are equivalent. 
\begin{enumerate}
\item The function $y(t,x)$ solves the control problem~\eqref{e:1Dcontrol-problemyn} (with initial datum $y_0(x)$ and control $h(t)$)
\item The function $v(t,x)=e^{-\f(x)/2\e}y(t,x)$ solves
  \begin{equation}
  \label{e:heat-control-eps-Witt}
\left\{
\begin{array}{rl}
\eps \partial_t v + P_\eps v =0  , & (t,x)\in (0,T)\times  (0,L) ,  \\
v(t,0) =e^{-\f(0)/2\e} h (t) , \quad v(t,L) = 0 ,&  t \in (0,T) , \\
v(0,x) = e^{-\f(x)/2\e} y_0 (x),&  x \in (0,L) .
 \end{array}
 \right.
 \end{equation}
\end{enumerate}
\end{lemma}
Note that in this lemma, ``solving'' the equation is meant in the classical sense for regular solutions but has to be taken in the transposition sense for ``rough solutions''. The conjugation works the same way in this weak sense according to the duality~\eqref{e:duality-boundary-case}. The latter now rewrites
\begin{align}
\label{e:duality-conjugated}
&\left( \zeta(T) , v(0) \right)_{L^2(0,L)} - \left( \zeta(0), v(T)\right)_{L^2(0,L)}  + \int_0^T  \eps \partial_x \zeta(t,0) v(T-t, 0) dt =0 , \quad \text{ that is to say, } \nonumber\\ 
& \left( \zeta(T) , e^{-\f/2\e} y_0  \right)_{L^2(0,L)} - \left( e^{\f/2\e} u_0, v(T)\right)_{L^2(0,L)}  + \int_0^T  \eps \partial_x \zeta(t,0) e^{-\f(0)/2\e}h(T-t) dt =0 .
\end{align}
with $v$ solving~\eqref{e:heat-control-eps-Witt} and $\zeta$ solving~\eqref{e:heat-transp-eps-Witt}.

\subsection{Controllability of the limit equation $\eps=0$} 
\label{s:limit-eq}
In this section, we consider the observability question for the formal control problem obtained from~\eqref{e:1Dcontrol-problemyn} in the limit $\eps=0$. It is a transport equation of hyperbolic type and the number of boundary conditions to be imposed for well-posedness is different from its parabolic counterpart. The limit equation that we expect on the control side is the following 
\begin{equation}
\label{e:1Dcontrol-transp}
\left\{
\begin{array}{rl}
(\d_t + \f' \d_x +\f'' - \q ) y= 0 ,& (t,x)\in (0,T) \times (0,L) , \\
y(0,x) = y_0 (x),&  x \in (0,L) ,
\end{array}
\right.
\end{equation}
where, again, the transport equation can be written equivalently as $(\d_t + \mathfrak{a} \d_x + \mathfrak{b})y=0$.
 We assume $\f'(0)\neq 0$ and $\f'(L)\neq 0$ for simplicity.
The expected boundary conditions actually depend on the sign of $\f'(0)$ and $\f'(L)$. Namely, the expected relevant boundary conditions in view of the parabolic control problem~\eqref{e:1Dcontrol-problemyn} for $\eps>0$ are, for $t \in (0,T)$, 
\begin{enumerate}
\item \label{case+-}$y (t,0) = h (t)$ and $y(t,L)=0$,  if $\f'(0)>0$ and $\f'(L)<0$,
\item \label{case++}$y (t,0) = h (t)$, if $\f'(0)>0$ and $\f'(L)>0$,
\item \label{case--}$ y(t,L)=0$, if $\f'(0)<0$ and $\f'(L)<0$,
\item \label{case-+} no boundary conditions to be imposed, if $\f'(0)<0$ and $\f'(L)>0$.
\end{enumerate}
Note that in most of the article, we actually assume that the vector field $\f'=\mathfrak{a}$ is $C^1$ and does not vanish on the interval $[0,L]$; the only two relevant boundary conditions are then~\ref{case++} in case $\f'>0$ and~\ref{case--} in case $\f'<0$.
If we denote $\Omega=(0,L)$, $\d \Omega=\{0,L\}$ and $\d_\nu$ the outgoing normal unit field to the boundary (i.e. $\d_\nu|_{x=L} =\d_x$ and $\d_\nu|_{x=0}=-\d_x$), we can write the previous boundary conditions in a more concise (but maybe more complicated) way by $y|_{\d \Omega\cap \{\d_\nu \f <0\} }=H(t,x)|_{\d \Omega\cap \{\d_\nu \f<0\}}$ where $H(t,x)$ is defined on $(0,T)\times \d \Omega$ by $H(t,0)=h(t)$ and $H(t,L)=0$.

Note that only Cases \ref{case+-} and \ref{case++} define a control problem. In Cases \ref{case--} and \ref{case-+}, the only relevant question is whether or not the solutions do vanish at time $T$.

Solutions to~\eqref{e:1Dcontrol-transp} are meant in the weak sense i.e. in the sense of transposition, see e.g.~\cite[Section~2.1.1]{Cor:book}. We say that $y$ is a solution~\eqref{e:1Dcontrol-transp} (with appropriate boundary conditions) in the sense of transposition if for all $\tau\in [0,T]$,
\bna
0=-\int_{0}^{\tau}\int_{0}^{L}y (\d_t \phi+ \f'(x)\d_x \phi+ \q \phi)~dtdx-\f'(0)\int_{0}^{\tau}h(t) \phi(t,0)~dt\\
+\int_{0}^{L} y(\tau,x)\phi(\tau,x)~dx- \int_{0}^{L} y_{0}(x)\phi(0,x)~dx
\ena
for every $\phi\in C^{1}([0,\tau]\times [0,L])$ satisfying $\phi(x,t)=0$ for all $t\in [0,\tau]$ and every $x\in \{0,L\}$ so that $\d_\nu \f(x)>0$, that is  
\begin{enumerate}
\item no assumption if $\f'(0)>0$ and $\f'(L)<0$,
\item  $\phi(t,L)=0$, $\forall t\in [0,\tau]$ , if $\f'(0)>0$ and $\f'(L)>0$,
\item $\phi(t,0)=0$, $\forall t\in [0,\tau]$ , if $\f'(0)<0$ and $\f'(L)<0$,
\item  $\phi(t,0)=0$ and  $\phi(t,L)=0$, $\forall t\in [0,\tau]$, if $\f'(0)<0$ and $\f'(L)>0$.
\end{enumerate}
The arguments of \cite[Proposition 1]{CG:05} can be adapted here to prove that the weak limit of solutions of system \eqref{e:1Dcontrol-problemyn} are solutions to~\eqref{e:1Dcontrol-transp} with the boundary conditions given in Items~\ref{case+-}--\ref{case-+}. In particular, this allows to prove Proposition~\ref{p:transport-seul}.

For the observation problem, we expect the following limit system
\begin{equation}
\label{e:1Dobs-transp}
\left\{
\begin{array}{rl}
(\d_t - \f'(x)\d_x - \q ) u = 0 ,& (t,x)\in (0,T) \times (0,L) , \\
u(0,x) = u_0 (x),&  x \in (0,L) .
\end{array}
\right.
\end{equation}
The boundary conditions are then the same as for $\phi$, namely:  
\begin{enumerate}
\item no boundary conditions to be imposed, if $\f'(0)>0$, $\f'(L)<0$,
\item $u (t,L) = 0$,  $t \in (0,T) , $ if $\f'(0)>0$, $\f'(L)>0$,
\item $ u(t,0)=0$,  $t \in (0,T) , $ if $\f'(0)<0$, $\f'(L)<0$,
\item $u (t,0) = 0$, $u(t,L)=0$,  $t \in (0,T)$, if $\f'(0)<0$, $\f'(L)>0$.
\end{enumerate}
Following closely~\cite[Section 2.1]{Cor:book}, it is possible to prove the following two lemmata.
\begin{lemma}
\label{l:exsie0}
For any $u_{0} \in L^{2}(0,L)$, the Cauchy problem \eqref{e:1Dobs-transp} with above boundary conditions has a unique solution $u\in C((0,T),L^{2}(0,L))$.

Moreover, for any $y_{0}\in L^{2}(0,L)$, $h\in L^{2}(0,T)$, the Cauchy problem \eqref{e:1Dcontrol-transp}  with above boundary conditions has a unique solution $y\in C((0,T),L^{2}(0,L))$ in the sense of transposition.
\end{lemma}
\begin{lemma}
We have the duality relation
\bna
\left<y(T),u_{0}\right>_{L^{2}(0,L)}-\left<y_{0},u(T)\right>_{L^{2}(0,L)}= D ,
\ena
with 
\begin{enumerate}
\item $D=\f'(0) \int_{0}^{T}u(t,0)h (t-T)dt$, if $\f'(0)>0$ and  $\f'(L)<0$,
\item $D=\f'(0) \int_{0}^{T}u(t,0)h (t-T)dt$, if $\f'(0)>0$ and $\f'(L)>0$,
\item $D=0$ if $\f'(0)<0$ and $\f'(L)<0$,
\item $ D=0$, if $\f'(0)<0$ and $\f'(L)>0$.
\end{enumerate}
Moreover, null-controllability (or the problem of having $y(T)=0$) holds true if and only if, for all $u$ solution of \eqref{e:1Dobs-transp} with above boundary conditions, we have
\begin{enumerate}
\item $\int_{0}^{T}|u(t,0)|^{2}dt\geq C \nor{u(T)}{L^{2}(0,L)}^{2}$, if $\f'(0)>0$ and $\f'(L)<0$,
\item $\int_{0}^{T}|u(t,0)|^{2}dt\geq C \nor{u(T)}{L^{2}(0,L)}^{2}$, if $\f'(0)>0$ and $\f'(L)>0$,
\item $u(T)=0$ if $\f'(0)<0$ and $\f'(L)<0$,
\item $ u(T)=0$, if $\f'(0)<0$, $\f'(L)>0$.
\end{enumerate}
\end{lemma}
The next Proposition simply says that null-controllability (or the problem of having $y(T)=0$) holds if and only if all the trajectories exit the interval.
 \begin{proposition}
 \label{p:tpstransp}
The conditions in the previous Lemma hold if and only if $\f' \neq 0$ in $[0,L]$ and $T\geq T_{\f'}=\int_0^{L} \frac{ds}{|\f'(s)|}$.
 \end{proposition}
 Note that the system considered is not the same depending on the sign of $\f'$.
 \bnp
 First, we can check that if there is one point $x_{0}\in (0,L)$ (it cannot be on the boundary with the assumptions) such that $\f'(x_{0})=0$, then the conditions are not fulfilled. Indeed, we can construct some non zero solutions localized arbitrary close to $x_{0}$ that remain zero close to the boundary. Note also that since $\f$ is sufficiently regular and $\f'(x_{0})=0$, then $\int_0^{L} \frac{ds}{\f'(s)}$ is not convergent.
 
 In the other case, we can write explicitly the solution. We first consider the second case $\f'(0)>0$, $\f'(L)>0$. 

For any $x\in [0,L]$, denote $J_{x}(x_{0})=\int_{x}^{x_{0}} \frac{ds}{\f'(s)}$. It is an increasing function from $[x,L]$ to $[0, T_{r,\f'}(x)]$ with $0\leq T_{r,\f'}(x):=\int_{x}^{L} \frac{ds}{\f'(s)}\leq T_{\f'}$  the exit time on the right starting from $x$. Denote $t\mapsto Y_{x}(t)$ its inverse from $[0, T_{r,\f'}(x)]$ to $[x,L]$. Deriving with respect to $t$ the equation $J_{x}(Y_{x}(t))=t$, we see that $\d_{t}Y_{x}(t)= \f'(Y_{x}(t))$ while deriving with respect to $x$ gives $\d_{x}Y_{x}(t)=\frac{\f'(Y_{x}(t))}{\f'(x)}$. Moreover, we have $Y_{x}(0)=x$ and $Y_{x}(T_{r,\f'}(x))=L$. 

We define for $x\in [0,L]$ and $t\geq 0$,
 \begin{equation*}
 \begin{array}{ll}
 u(t,x)=e^{\int_{0}^{t}q(Y_{x}(\tau))d\tau}u_{0}(Y_{x}(t)), &\quad  (t,x)\in \R_{+}\times [0,L] \textnormal{ and }0\leq t\leq T_{r,\f'}(x) , \\
 u(t,x)=0, & \quad  (t,x)\in \R_{+}\times [0,L] \textnormal{ and }t> T_{r,\f'}(x) . 
 \end{array}
 \end{equation*}
 Note first that it is well defined since $0\leq\tau\leq  t\leq T_{r,\f'}(x)$ implies that $Y_{x}(\tau)$ and $Y_{x}(t)$ are well defined. We also notice that if $u_{0}$ is $C^{1}$ with $u_{0}(L)=0$, then $u$ is solution of \eqref{e:1Dobs-transp} in the classical sense with the appropriate boundary conditions $u (t,L) = 0$. 
 
We first check that $u$ is continuous and it is therefore sufficient to verify the equation in each zone. We compute for $0\leq t\leq T_{r,\f'}(x)$,
 \bna
 \d_{t}u(t,x)= q(Y_{x}(t))e^{\int_{0}^{t}q(Y_{x}(\tau))d\tau}u_{0}(Y_{x}(t))+e^{\int_{0}^{t}q(Y_{x}(\tau))d\tau}\f'(Y_{x}(t))u_{0}'(Y_{x}(t))\\
  \d_{x}u(t,x)= \left(\int_{0}^{t}\frac{\f'(Y_{x}(\tau))}{\f'(x)}q'(Y_{x}(\tau))d\tau \right)e^{\int_{0}^{t}q(Y_{x}(\tau))d\tau}u_{0}(Y_{x}(t))+e^{\int_{0}^{t}q(Y_{x}(\tau))d\tau}\frac{\f'(Y_{x}(t))}{\f'(x)}u_{0}'(Y_{x}(t))
 \ena
 We conclude by remarking that $\int_{0}^{t}\f'(Y_{x}(\tau))q'(Y_{x}(\tau))d\tau=\int_{0}^{t}\frac{\d}{\d t}\left(q(Y_{x}(\tau))\right)d\tau=q(Y_{x}(t))-q(x)$. Note that the computations only make sense if $u$ and $q$ are regular enough, but we obtain the same result in general by approximation. For the boundary conditions, $T_{r,\f'}(L)=0$ so that we always have $t> T_{r,\f'}(L)=0$ and $u(t,L)=0$. Also, the assumption $u_{0}(L)=0$ ensures that at time $t= T_{r,\f'}(x)$, $Y_{x}(T_{r,\f'}(x))=L$, so that $u(T_{r,\f'}(x),x)=0$ and the function $u$ is continuous. The formula extends to $L^{2}$ functions and therefore defines the flow map described in Lemma \ref{l:exsie0}.
 
 Now, we have to check if the defined formula fulfills or not the observability estimate. For $x=0$, we have $T_{r,\f'}(0)=T_{\f'}$ and we have seen that $Y_{0}(t)$ is an increasing bijection from $[0, T_{\f'}]$ to $[0,L]$. We can then compute for $0\leq T\leq T_{\f'}$.
 \bna
 \int_{0}^{T}|u(t,0)|^{2}dt\approx \int_{0}^{T}|u_{0}(Y_{0}(t))|^{2}dt=\int_{0}^{Y_{0}(T)}|u_{0}(y)|^{2}\frac{dy}{\f'(y)}\approx \int_{0}^{Y_{0}(T)}|u_{0}(y)|^{2}dy ,
 \ena
 where we have made the substitution $y=Y_{0}(t)$, that is $t=J_{0}(y)$. The symbol $a\approx b $ means that there exists one constant $C$ depending on $T$, $\f$, $L$ and $q$ so that $C^{-1}a\leq b\leq Ca$. For $T\in [0, T_{r,\f'}(x)]$ $\d_{x}Y_{x}(T)=\frac{\f'(Y_{x}(T))}{\f'(x)}>0$, so that $x\mapsto \phi_{T}(x):=Y_{x}(T)$ is a diffeomorphism from $[0,x_{T}]$ to $[Y_{0}(T),L]$ where $x_{T}\in [0,L]$ is so that $Y_{x_{T}}(T)=L$ (which implies $T_{r,\f'}(x_{T})=T$ since $Y_{x}(T_{r,\f'}(x))=L$ by definition).
 
In particular, for $0\leq T\leq T_{\f'}$, we have  $0\leq T\leq T_{r,\f'}(x) \Leftrightarrow x_{T}\leq x\leq L$ and therefore
 \bna
 \nor{u(T)}{L^{2}(0,L)}^{2}=\int_{0}^{L}|u(T,x)|^{2}dx\approx \int_{x_{T}}^{L}|u_{0}(Y_{x}(T))|^{2}dx\approx \int_{Y_{0}(T)}^{L}|u_{0}(y)|^{2}dy.
 \ena
 In particular, the observability inequality holds if and only if $Y_{0}(T)=L$, that is $T=T_{\f'}$. For $T\geq T_{\f'}$, we have $u(T)=0$ so that the observability is trivial. This ends the result in the case $\f'(0)>0$ and $\f'(L)>0$. For the other case $\f'(0)<0$ and $\f'(L)<0$, the change of variable $x \leftrightarrow L-x$ reduces to the same case as before. The only difference now is that we want the solution to be zero instead of the observability. The condition is actually the same.
 \enp
 
\section{Proofs of the main results}
\subsection{Localization of Schr\"odinger eigenfunctions in a one dimensional well}
\label{sect1D}
In this section, we recall results proved in the companion paper~\cite{LL:22} in which we study localization properties for  eigenfunctions of the Schr\"odinger operator $P_\eps$ defined in~\eqref{e:def-Peps}. 
We now state these two results, which take the form of uniform (in terms of both $\eps$ and $E$) upper and lower bounds on eigenfunctions.

\begin{theorem}[Upper bounds for eigenfunctions:~\cite{LL:22} Theorem~1.3]
\label{t:uniform-agmon}
Assume that $V \in C^1([0,L])$ ($\f \in C^2([0,L])$), and that Item~\ref{A2} in Assumption~\ref{assumptions} is satisfied.
Then, for all $\delta>0$ there exist $\eps_0 = \eps_0(\delta) \in (0,1]$ such that for all $E, \psi$ solution to 
\begin{align}
\label{e:eignfct-E-bis}
P_\eps \psi = E \psi , \quad \psi \in H^2(0,L)\cap H^1_0(0,L) , \quad \nor{\psi}{L^2(\M)} =1, 
\end{align}
 we have for all $\eps< \eps_0$
\begin{align}
\label{e:agmon-H1}
 \nor{ e^{\frac{d_{A,E}}{\eps}}\frac{\eps}{\sqrt{|E|+1}}\psi'}{L^2} + \nor{ e^{\frac{d_{A,E}}{\eps}}\psi}{L^2} \leq e^{\frac{\delta}{\eps}}  .\\
\label{e:Agmon-boundary}
\frac{\eps}{\sqrt{|E|+1}}|\psi'(0)| \leq e^{- \frac{d_{A,E}(0)-\delta}{\eps}} , \quad \frac{\eps}{\sqrt{|E|+1}}|\psi'(L)| \leq e^{- \frac{d_{A,E}(L)-\delta}{\eps}}  .
\end{align}
\end{theorem}
Only continuity of the potential $V$ (i.e. $\f \in C^1([0,L])$) is assumed in~\cite{LL:22}, but this refinement is not relevant here since $V \in C^1([0,L])$ (i.e. $\f \in C^2([0,L])$) is needed to use Theorem~\ref{t:allibert-lower-uniform}.

\begin{theorem}[Lower  bounds for eigenfunctions:~\cite{LL:22} Theorem~1.4]
\label{t:allibert-lower-uniform}
Assume that $V \in C^1([0,L])$ ($\f \in C^2([0,L])$), and that Item~\ref{A2} in Assumption~\ref{assumptions} is satisfied.
Then, for any $\mathbf{y}_0\in [0,L], \nu>0$ and any $\delta >0$, there is $\eps_0>0$ such that for all $E,\psi$ satisfying~\eqref{e:eignfct-E-bis},
 we have for all $\eps<\eps_0$, 
\begin{align}
&\nor{\psi}{L^2(U)} \geq e^{-\frac{1}{\eps}(d_{A,E}(U) + \delta)} ,\quad d_{A,E}(U) = \inf_{x\in U}d_{A,E}(x) , \quad  U=(\mathbf{y}_0-\nu, \mathbf{y}_0+\nu) \cap [0,L] , \nonumber  \\
&\frac{\eps}{\sqrt{|E|+1}} |\psi'(0)|  \geq e^{-\frac{1}{\eps}(d_{A,E}(0) + \delta)}, \quad \frac{\eps}{\sqrt{|E|+1}} |\psi'(L)|  \geq e^{-\frac{1}{\eps}(d_{A,E}(L) + \delta)} \label{e:bound-obs}  .
\end{align}
\end{theorem}
Note that this improved lower bound is as precise as the upper bound~\eqref{e:agmon-H1} (except for the $\delta$ loss) and thus essentially optimal. 
Uniformity with respect to the energy level $E$ is necessary for the proof of the cost of controllability in Theorem~\ref{thmcontrol1D}. The latter indeed involves all the spectrum of $P_{\e}$ since we are studying all solutions. 

Remark that Theorems~\ref{t:allibert-lower-uniform} and~\ref{t:uniform-agmon} are counterparts one to the other. They state essentially that, in this very particular one dimensional setting, an eigenfunction $\psi$ associated to the energy $E$ satisfies $|\psi(x)|\sim e^{-\frac{d_{A,E}(x)}{\eps}}$ (and that this is uniform in $E, x,\eps$). 
The symbol $\sim$ is slightly abusive in our setting since we only have equivalence up to multiplicative terms of the form $e^{\frac{\delta}{\eps}}$, which can be very large. Yet, in the present context where only exponentially small quantities are compared, this kind of estimates is sufficient for our purposes and provides with the correct asymptotics. See e.g.~\cite[end of Section~1]{LL:22} for a discussion on possible refinements.

\subsection{Proof of Theorem~\ref{t:estim-Cobs-expo-1D} from Theorems~\ref{t:allibert-lower-uniform} and~\ref{t:uniform-agmon}}
In this section, we give a proof of Theorem~\ref{t:estim-Cobs-expo-1D}. The latter relies on consequences of Theorems~\ref{t:allibert-lower-uniform} and~\ref{t:uniform-agmon} that are not using the uniformity in $E$ and could be deduced from softer versions of these two results.
\begin{proposition}
\label{prop:improved-lower-1D}
Under the assumptions of Theorem~\ref{t:allibert-lower-uniform}, we have
\begin{align*}
\nor{e^{-\frac{\f}{2\eps}} \psi_\eps}{L^2(\M)} \geq  e^{-\frac{1}{\eps}(\min_{\M}W_E + \delta)} .
\end{align*}
\end{proposition}
\bnp[Proof of Proposition~\ref{prop:improved-lower-1D} from Theorem~\ref{t:allibert-lower-uniform}]
We take $x_m \in \M$ such that $W_E(x_m) = \min_{\M} W_E$. Then, by continuity, there is $\rho >0$ such that for all $x \in (x_m-\rho , x_m+\rho)$, we have $\f (x) < \f(x_m) + \delta$ and $d_{A,E}(x)< d_{A,E}(x_m) +\delta$. 
We then estimate
\begin{align*}
\nor{e^{-\frac{\f}{2\eps}} \psi_\eps}{L^2(\M)} & 
\geq  \nor{e^{-\frac{\f}{2\eps}} \psi_\eps}{L^2(x_m-\rho , x_m+\rho)} 
\geq e^{-\frac{1}{2\eps}(\f(x_m)+\delta)} \nor{ \psi_\eps}{L^2(x_m-\rho , x_m+\rho)} \\
& \geq e^{-\frac{1}{2\eps}(\f(x_m)+\delta)}
  e^{-\frac{1}{\eps}(d_{A,E}((x_m - \rho , x_m+\rho))  + \delta)} ,
\end{align*}
after having used Theorem~\ref{t:allibert-lower-uniform} in the last inequality for $\eps < \eps_0(\delta)$, and with $d_{A,E}((x_m - \rho , x_m+\rho)) = \inf_{x\in (x_m - \rho , x_m+\rho)}d_{A,E}(x)< d_{A,E}(x_m) +\delta$. As a consequence, we obtain 
\begin{align*}
\nor{e^{-\frac{\f}{2\eps}} \psi_\eps}{L^2(\M)} & \geq e^{-\frac{1}{\eps} (\frac{\f}{2}(x_m)+\frac{\delta}{2})}
 e^{-\frac{1}{\eps}(d_{A,E}(x_m) + 2 \delta)} \\
 & \geq  e^{-\frac{1}{\eps} (W_E(x_m)+3\delta)} = e^{-\frac{1}{\eps}(\min_{\M} W_E+3\delta)},
\end{align*}
where we have used the definition of $x_m$ in the last inequality.
\enp

We conclude this section with a proof of Theorem~\ref{t:estim-Cobs-expo-1D}, which relies on both Theorems~\ref{t:allibert-lower-uniform} (under the statement of Proposition~\ref{prop:improved-lower-1D}) and~\ref{t:uniform-agmon}.
\bnp[Proof of Theorem~\ref{t:estim-Cobs-expo-1D} from Theorem~\ref{t:uniform-agmon} and Proposition~\ref{prop:improved-lower-1D}]
We follow the proof of \cite[Theorem~3.1]{LL:18vanish}.
We first use Lemma~\ref{lemequiveunif} and Lemma~\ref{l:obs-cont} to see that we have~\eqref{e:obs-heat-transp-eps-Wittboundary} for all solutions to~\eqref{e:heat-transp-eps-Witt}, 
with $C_0 =C_0(T,\eps) = \mathcal{C}_0(T,\eps)$. 
Second, given $E \in V([0,L])$ and $\eps \in (0,1]$, there is $E_\eps = E + \grando{\eps^{2/3}}$ and $\psi_\eps \in H^2(\M) \cap H^1_0(\M)$ such that $P_\eps \psi_\eps = E_\eps \psi_\eps$ (see Lemma \ref{l:exist-eig} below). 
We then test~\eqref{e:obs-heat-transp-eps-Wittboundary} with $\zeta(t) = e^{-tE_\eps / \eps}\psi_\eps$, solution to~\eqref{e:heat-transp-eps-Witt} with $\zeta(0, x) = \psi_\eps(x)$.
This reads 
\begin{equation}
\label{e:tototiti}
e^{-\frac{2TE_\eps}{\eps}} \nor{e^{-\frac{\f}{2\eps}} \psi_\eps}{L^2(0,L)}^2 = \nor{e^{-\frac{\f}{2\eps}} \zeta(\TO)}{L^2(0,L)}^2 \leq C_0^2 \int_0^{\TO}  \left|e^{-\frac{\f(0)}{2\eps}} \eps \d_x \zeta(t,0)\right|^2 dt \leq C_0^2 Te^{-\frac{\f(0)}{\eps}}  \left|\eps \psi_\eps'(0)\right|^2 .
\end{equation}
Using Proposition~\ref{prop:improved-lower-1D}, we have $e^{-\frac{2TE_\eps}{\eps}}\nor{e^{-\frac{\f}{2\eps}}\psi_\eps}{L^2(\M)}^2 \geq e^{-\frac{2T(E+\delta)}{\eps}}e^{-\frac{2}{\eps}(\min_{\M}W_E + \delta)}$ for all $0<\eps<\eps_0(\delta)$ and using \eqref{e:Agmon-boundary} in Theorem~\ref{t:uniform-agmon} we have $\eps |\psi_\eps'(0)| \leq e^{- \frac{d_{A,E}(0)-\delta}{\eps}}$ (recall that $E$ is fixed). Combining together with~\eqref{e:tototiti}, these two inequalities yield, for $0<\eps<\eps_0(\delta, T)$,
$$
e^{-\frac{2T(E+\delta)}{\eps}}e^{-\frac{2}{\eps}(\min_{\M}W_E + \delta)}
\leq C_0^2e^{-\frac{\f(0)}{\eps}} e^{- 2\frac{d_{A,E}(0)-\delta}{\eps}}=C_0^2e^{-2 \frac{W_E (0)-\delta}{\eps}} ,
$$
which is the first statement of the theorem when recalling $C_0  = \mathcal{C}_0(T,\eps)$ and changing the notation for $\delta$.
\enp
In the course of the proof, we have used the following Lemma, taken from \cite{LL:18vanish},
 proving the existence of eigenvalues at any allowed energy level. It can also be deduced from the much more precise Theorem \ref{t:Peps-spect-Allib} adapted from \cite{Allibert:98}.
\begin{lemma}\cite[Lemma~3.2]{LL:18vanish}
\label{l:exist-eig}
Assume $V \in W^{1,\infty}(\M)$ and $\qf \in L^\infty(\M)$  are both real valued. 
For all $E \in V(\M) = [\min_{\M} V, \max_{\M} V]$ and all $\eps \in (0,1]$, there is $E_\eps = E + \grando{\eps^{2/3}}$ and $\psi_\eps \in H^2(\M) \cap H^1_0(\M)$ such that $P_\eps \psi_\eps = E_\eps \psi_\eps$.
\end{lemma}

\subsection{Coron-Guerrero type lower bound: proof of Theorem~\ref{thmlower+-}}

\bnp[Proof of Theorem~\ref{thmlower+-}]
Recall the simpler way of writting the control problem \eqref{e:1Dcontrol-problemyn}, the observation problem \eqref{e:1Dobs-problem} and the duality statement \eqref{e:duality-boundary-case}. Let $(\varphi_k^{\e})_{k\in \N}$ denote the sequence of eigenfunctions of the selfadjoint operator $P_{\e}$, associated with eigenvalues $\lambda_k^\eps$ sorted in increasing order.

Now, we fix $E\in V([0,L])$.  As a consequence of~\cite[Lemma~3.2]{LL:18vanish}, there exists a sequence of eigenvalues $E_\eps$ of  $P_\eps$ with $E_\eps \to E$ as $\eps \to 0$. That is to say, there is $\psi_\eps \in H^2(0,L) \cap H^1_0(0,L)$ such that $P_\eps \psi_\eps = E_\eps \psi_\eps$. Denote $n=n(E,\eps)\in \N$ the index (in the non-decreasing sequence of eigenvalues of $P_\eps$) such that $\lambda_n^\eps=E_\eps$ and $\psi_{\e}=\varphi_{n}^{\e}$.
We choose for initial datum for~\eqref{e:1Dcontrol-problemyn} the function 
$$y_0=y_n = e^{\frac{\f}{2\e}} \psi_\eps = e^{\frac{\f}{2\e}} \varphi_{n}^{\e} .$$ 
We denote by $h_n$ any control driving the initial datum $y_n$ to zero and produce lower bounds for its norm.
According to the Agmon estimate~\eqref{e:agmon-H1}, we have 
\begin{align}
\label{e:norm-yn}
\nor{y_n}{L^2(0,L)}  = \nor{e^{\frac{\f}{2\e}} \psi_\eps}{L^2(0,L)} =\nor{e^{\frac{1}{\eps}(\frac{\f}{2} - d_{A,E})}  e^{\frac{d_{A,E}}{\eps}}\psi_\eps}{L^2(0,L)} \leq e^{\frac{D_E+\delta}{\e}} ,
\end{align}
with 
\begin{align}
\label{e:norm-yn-bis}
D_E=\sup_{[0,L]}\widetilde{W_E} , \quad \widetilde{W_E}(x) : =\frac{\f}{2}(x)-d_{A,E}(x). 
\end{align}

We remark that the function $(t,x) \mapsto e^{-\frac{\lambda_k^{\e} t}{\e}}\varphi_k^{\e}(x)$ is a solution of \eqref{e:heat-transp-eps-Witt}. As a consequence of Lemma~\ref{lemequiveunif}, the function
$u_k(t,x) =e^{-\frac{\f(x)}{2\e}}\varphi_k^{\e}(x)e^{-\frac{\lambda_k^{\e}}{\e} t}$ is solution to~\eqref{e:heat-transp-eps}, that is of \eqref{e:1Dobs-problem} with $u_0(x)=e^{-\frac{\f(x)}{2\e}}\varphi_k^{\e}(x)$. Since we assume $h_n$ is a null-control, we have $y_n(T)=0$ and the duality formula \eqref{e:duality-boundary-case} taken for $u=u_k$ implies
\bnan
\label{formuledualite}
 ( u_k(T) , y_n )_{L^2(0,L)}  + \int_0^T  \eps \partial_x u_k(t,0) h_n (T-t) dt =0 .
\enan
Since $u_k(t,0)=0$ (Dirichlet boundary conditions), we have $\partial_x u_k(t,0)=e^{-\frac{\f(0)}{2\e}}e^{-\frac{\lambda_k^{\e} t}{\e}}(\varphi_k^{\e})'(0)$. Moreover, we have 
 $$ ( u_k(T) , y_n )_{L^2(0,L)} = ( e^{-\frac{\f(x)}{2\e}}\varphi_k^{\e}(x)e^{-\frac{\lambda_k^{\e}}{\e} T} , e^{\frac{\f}{2\e}} \varphi_{n}^{\e})_{L^2(0,L)} =  e^{-\frac{\lambda_k^{\e}}{\e} T}  \delta_{k,n} .$$
 These two identities together with~\eqref{formuledualite} (and the change $T-t \leftarrow t$ in the integral) imply
\bnan
\label{e:delta-d}
\int_0^T h_n(t)e^{\frac{\lambda_k^\e}{\e} t}dt=-\frac{e^{\frac{\f(0)}{2\e}}}{\e(\varphi_k^\eps)'(0)}\delta_{k,n}, \quad \text{ for all }k\in \N ,
\enan
We next set  for $n \in \N$
\bna
v_n(s):= \F (h_n \mathds{1}_{[0,T]})(s) = \int_{0}^{T}h_n(t)e^{-ist}dt  , \quad s\in \C ,
\ena
which defines an entire function $v_n:\C\rightarrow \C$.
Identity~\eqref{e:delta-d} reformulates as 
\bna
v_n\left(i\frac{\lambda_k^\eps}{\e}\right)=-\frac{e^{\frac{\f(0)}{2\e}}}{\e(\varphi_n^\eps)'(0)}\delta_{k,n}.
\ena
Moreover, writing $f(s)_+ = \max \{ f(s) , 0\}$, we have for all $T\geq 0$ 
\begin{align*}
|v_n(s)| \leq e^{T\Im(s)_+}\int_0^T|h_n(t)|dt \leq T^{1/2} e^{T\Im(s)_+}\|h_n\|_{L^2(0,T)}, \quad \text{ for all } s\in \C  .
\end{align*}
We now introduce the parameter $B\geq 0$. We define the entire function 
\bna
g_n:\C\rightarrow \C , \quad 
g_n(s):=v_n\left(\frac{s-iB}{\e}\right) .
\ena
From the above properties of $v_n$ we obtain
\begin{align}
\label{e:b-k-raci}
g_n(i\mathsf{b}_k)=-\frac{e^{\frac{\f(0)}{2\e}}}{\e(\varphi_n^{\e})'(0)}\delta_{k,n} , \quad \mathsf{b}_k=\lambda_k^\e+B ,
\end{align}
together with the general bound  
\begin{align}
\label{generalboundg+}
|g_n(s)| \leq   T^{1/2} e^{\frac{T}{\eps}(\Im(s)-B)_+}\|h_n\|_{L^2(0,T)}, \quad \text{ for all } s\in \C  .
\end{align}
Now, we want to apply the complex analysis Lemma \ref{lmanalysecomplexe} below with the following parameters:
\begin{itemize}
\item $x:=E+B>0$, and $x_\eps := \lambda_n^\e+B \to x$ as $\eps \to 0^+$,
\item $g:=g_n$,
\item $\sigma=\limsup_{y \to +\infty}\frac{\log |g_n(iy)|}{y}\leq \frac{T}{\e}$ according to~\eqref{generalboundg+},
\item  $\log|g(\tau)| \leq C_{\R} := \log(T^{1/2}\|h_n\|_{L^2(0,T)})$ for all  $\tau\in\R$, according to~\eqref{generalboundg+} (using that $B\geq 0$),
\item $\mathsf{b}_k =\lambda_k^\e+B$ and $b_k:=\mathsf{b}_k$ for $k<n$ and $b_k:=\mathsf{b}_{k+1}$ for $k\geq n$ (that is to say $\{b_k,k\in \N \} = \{\mathsf{b}_k,k\in \N ,k\neq n\}$). According to~\eqref{e:b-k-raci} applied with $k\neq n$, the sequence $(b_k)_{k\in \N}$ satisfies $g(ib_k)=0$. Moreover, the assumption~\eqref{e:bkZR} is satisfied with $Z(s):=\Phi^{-1}(\pi s)+B$ according to estimate~\eqref{stimlambda} in Theorem~\ref{t:Peps-spect-Allib}. Note that this uses $\q= \frac{\f''}{2}$, that is to say $\q_\f=0$. We also recall that the function $\Phi$ is defined in~\eqref{e:def-phi}.
 \end{itemize} 
 Application of Lemma \ref{lmanalysecomplexe} implies 
 $$
 \log|g_n(i(\lambda_n^\e+B))| \leq \frac1{\eps} (- T_{E,B} + \delta) + \frac{T}{\eps}(E+B) +  \log(T^{1/2}\|h_n\|_{L^2(0,T)}), 
 $$ 
 where we have set 
 \begin{align*}
T_{E,B} : = I(E+B) = \frac{1}{\pi}\int_{0}^{+\infty}\log \left|\frac{\Phi^{-1}(y)+E+2B}{\Phi^{-1}(y)-E}\right|dy
=\frac{1}{\pi}\int_{V(\xzero)}^{+\infty}\log \left|\frac{x+E+2B}{x-E}\right|\Phi'(x)dx ,
\end{align*}
and, in the last expression used the definition of $\Phi$ and the properties of $V$ to write $\Phi^{-1}(0) = V(\xzero) = \min V$.
According to~\eqref{e:b-k-raci} applied with $k=n$, we also have 
 $$
\log \left| \frac{e^{\frac{\f(0)}{2\e}}}{\e(\varphi_n^{\e})'(0)}\right|  = \log|g_n(i\mathsf{b}_n)| = \log|g_n(i(\lambda_n^\e+B))|  .
 $$
Combining these two lines, we obtain,
\bnan
\label{egal123}
\log\left|\frac{e^{\frac{\f(0)}{2\e}}}{\e(\varphi_n^\eps)'(0)}\right|\leq \frac{1}{\e}\left(-T_{E,B}+T(E+B)+ \delta\right)+\log(T^{1/2}\|h_n\|_{L^2(0,T)}).
\enan
Moreover, thanks to the Agmon estimate \eqref{e:Agmon-boundary}, we have, for $\eps \in (0,\eps_0)$, $ |\eps(\varphi_n^\eps)'(0)|\leq e^{\frac{-d_{A,E}(0)+\delta}{\e}}$.
As a consequence, we have
\bnan
\label{lowerbn}
\log\left|\frac{e^{\frac{\f(0)}{2\e}}}{\e(\varphi_n^\eps)'(0)}\right|\geq \frac{\frac{\f(0)}{2} +d_{A,E}(0)-\delta}{\e}=\frac{W_{E}(0)-\delta}{\e}.
\enan
Combining \eqref{lowerbn} together with~\eqref{egal123}, we finally obtain 
\begin{align}
\label{e:antepen-est}
\log(T^{1/2}\|h_n\|_{L^2(0,T)}) \geq
\frac{1}{\e}\left(W_{E}(0) +T_{E,B} - T(E+B)- 2\delta \right). 
\end{align}
Finally, assuming observability/controllability, Lemma \ref{l:obs-cont} implies that the control cost (observability constant) necessarily satisfies
$$
C_0(T,\e) \geq \frac{\|h_n\|_{L^2(0,T)}}{\nor{y_n}{L^2(0,L)}} , \quad \text{ that is,}\quad   \log C_0(T,\e) \geq \log \|h_n\|_{L^2(0,T)}  - \log \nor{y_n}{L^2(0,L)} .
$$ 
Recalling~\eqref{e:norm-yn}-\eqref{e:norm-yn-bis} together with~\eqref{e:antepen-est}, we have now obtained, for $\eps \in (0,\eps_0)$ small enough
$$
 \log C_0(T,\e) \geq \frac{1}{\e}\left(W_{E}(0) +T_{E,B} - T(E+B)- 3\delta  - D_E - \delta \right) ,
$$
which concludes the proof of the theorem. 
\enp 

The proof of the above result relied on the following complex analysis lemma.
\begin{lemma}
\label{lmanalysecomplexe}
Let $Z:\R^+\mapsto \R^+$ be a continuous strictly increasing function such that $Z^{-1}$ is locally Lipschitz continuous on $[Z(0),+\infty)$ and $\frac{1}{Z}\in L^1([1,+\infty[)$, and set 
$$I : \R\to \R , \quad  I(x) : =\int_{0}^{+\infty}\log \left|\frac{Z(y)+x}{Z(y)-x}\right|dy , $$
(which, under the above assumptions, is well-defined and continuous on $\R^+$). 
Let $R: \R^+ \to \R^+$ be an increasing function tending to zero at zero. Then, for any $x>0$, $\delta>0$ $D>0$, and any family $(x_\eps)_{\eps \in (0,\eps_0)}$ such that $x_\eps \to x$ as $\eps \to 0^+$, there exists $\e_0$ so that for any holomorphic function $g$ on $\C^+$ satisfying
\begin{enumerate}
\item $g$ is of exponential type on  $\C^+$; we write $\sigma :=\limsup_{y \to +\infty}\frac{\log |g(iy)|}{y} <+\infty$,
\item \label{e:real-line} $\tau \mapsto \log|g(\tau)|$ is bounded above on $\R^+$; we write $C_{\R} := \sup_{\tau \in \R }\log|g(\tau)|  <+\infty$,
\item $g(ib_k)=0$ for any $k\in\N^*$ where $b_k = b_k(\eps)>0$ is a sequence such that 
\begin{align}
\label{e:bkZR}
Z(\e k-D\e)-R(\e)\leq b_k \leq   Z(\e k+D\e)+R(\e) , \quad \text{for all } k \geq D , \eps \in (0,\eps_0) ,
\end{align}
\end{enumerate}
we have 
\begin{align}
\label{e:estim-imagg}
\log|g(ix_\eps)|\leq \frac{-I(x)+\delta}{\e}+\sigma x_{\e}+C_{\R},  \quad \text{ for all }\eps \in (0,\eps_0) .
\end{align}
\end{lemma}
We first prove the following lemma, proving in particular that $I(x)$ is well-defined and continuous.
\begin{lemma}
\label{l:continu-I}
Given $Z:\R^+\mapsto \R^+$ a continuous strictly increasing function such that $Z^{-1}$ is locally Lipschitz continuous on $\R^+$ and $\frac{1}{Z}\in L^1([1,+\infty[)$, the functions
$$
F(a,b,c,d) = \int_c^{d}\log \left|\frac{Z(y)+b}{Z(y)-a}\right|dy ,
 \qquad 
 F_\infty(a,b,c) = \int_c^{+\infty}\log \left|\frac{Z(y)+b}{Z(y)-a}\right|dy ,
$$
are well-defined and continuous in $(a,b,c,d) \in \R^4$, resp. $(a,b,c) \in \R^3$.
 \end{lemma}
\bnp
Concerning first the function $F$, it suffices by linearity to check that $\int_c^{d}\log \left|Z(y)-a\right|dy$ is well-defined and continuous in $(a,c,d)$. The change of variable formula for Lipschitz map ($Z^{-1}$ is locally Lipschitz continuous) yields
$$
\mathcal{I} (a,c,d) :=\int_c^{d}\log \left|Z(y)-a\right|dy = \int_{Z(c)}^{Z(d)}\log \left|x-a\right|  \ (Z^{-1})'(x) dx ,
$$
and that the left hand-side is well-defined/continuous if and only if so is the right hand-side. But the right hand-side is well defined since $(Z^{-1})'$ is bounded (a.e.) on every compact interval and $\log|x|$ is integrable on compact sets. Let us now prove by hand that the right hand-side is continuous. Fix $(a,c,d) \in \R^3$ and let $ \eps = (\eps_1,\eps_2,\eps_3) \to 0$. We write 
\begin{align}
& \left| \mathcal{I} (a+\eps_1,c+\eps_2,d+\eps_3)  - \mathcal{I} (a,c,d) \right|  \nonumber \\
& = \left| \int_{Z(c+\eps_2)}^{Z(d+\eps_3)}\log \left|x-a-\eps_1\right|  \ (Z^{-1})'(x) dx - \int_{Z(c)}^{Z(d)}\log \left|x-a\right|  \ (Z^{-1})'(x) dx \right|   \leq I_1(\eps) + I_2(\eps) + I_3(\eps) 
\label{e:continuity-I}
\end{align}
with 
\begin{align*}
& I_1(\eps) = \left| \int_{Z(c)}^{Z(c+\eps_2)}\log \left|x-a\right|  \ (Z^{-1})'(x) dx \right| ,\quad 
 I_2(\eps) = \left| \int_{Z(d)}^{Z(d+\eps_3)}\log \left|x-a\right|  \ (Z^{-1})'(x) dx \right| \\
 & I_3(\eps) = \left| \int_{Z(d+\eps_3)}^{Z(c+\eps_2)}\log \left|\frac{x-a-\eps_1}{x-a}\right|  \ (Z^{-1})'(x) dx \right| .
\end{align*}
We have $I_1(\eps) +I_2(\eps) \to 0$ as $\eps\to 0$ by dominated convergence, and it only remains to study $I_3(\eps)$. Using that $(Z^{-1})'\in L^\infty_{\loc}(\R)$ and choosing $D_1,D_2\in \R$ such that $a, Z(d+\eps_3) , Z(c+\eps_2) \in (D_1,D_2)$ for all $\eps$ sufficiently small, we have 
$$
I_3(\eps) \leq C \int_{D_1}^{D_2} \left|\log \left|\frac{x-a-\eps_1}{x-a}\right| \right|dx = C \int_{D_1}^{a} \left|\log \left|1 - \frac{\eps_1}{x-a}\right| \right|dx +C \int_{a}^{D_2} \left|\log \left|1 - \frac{\eps_1}{x-a}\right| \right|dx .
$$
Assuming that $\eps_1 \geq 0$ (the case $\eps_1\leq 0$ is treated similarly) and changing variables in these two integrals implies 
\begin{align}
\label{e:I_3}
I_3(\eps) \leq  C \eps_1 \int_{-\infty}^{\frac{\eps_1}{D_1-a}} \left|\log \left|1 - s\right| \right|\frac{ds}{s^2} +C \eps_1 \int_{\frac{\eps_1}{D_2-a}}^{+\infty} \left|\log \left|1 -s\right| \right|  \frac{ds}{s^2}.
\end{align}
Then, we write, for $\eps_1$ small
$$
\eps_1 \int_{\frac{\eps_1}{D_2-a}}^{+\infty} \left|\log \left|1 -s\right| \right|  \frac{ds}{s^2} = \eps_1 \int_{\frac{\eps_1}{D_2-a}}^{1/2} \left|\log \left|1 -s\right| \right|  \frac{ds}{s^2} + \eps_1 \int_{1/2}^{+\infty} \left|\log \left|1 -s\right| \right|  \frac{ds}{s^2} .
$$
The function $\frac{1}{s^2}\left|\log \left|1 -s\right| \right|$ is integrable on $[1/2,+ \infty)$ and hence the second term converges to $0$. For the second term, we write $|\log|1-s|| \leq K |s|$ (with $K=2\log2$) on $[0,1/2]$ and thus
$$
\eps_1 \int_{\frac{\eps_1}{D_2-a}}^{1/2} \left|\log \left|1 -s\right| \right|  \frac{ds}{s^2} \leq \eps_1 \int_{\frac{\eps_1}{D_2-a}}^{1/2} Ks  \frac{ds}{s^2} = K \eps_1  \log\left(\frac{D_2-a}{2 \eps_1}\right) \to 0 , \quad \text{as }\eps\to 0 .
$$
This implies that the second term in the right hand-side of~\eqref{e:I_3} converges to zero. The first term is treated similarly, and we deduce that $I_3(\eps)\to 0$ as $\eps \to 0$. In view of~\eqref{e:continuity-I}, this implies that $\mathcal{I} (a,c,d)$ is continuous on $\R^3$ and thus $F$ is continuous on $\R^4$.

We now turn to the study of $F_\infty$, and remark that it suffices now to prove that $F_\infty(a,b,0)$ is well-defined and continuous on $\R^2$.
We have from the assumptions that $Z(y) \to +\infty$ as $y \to +\infty$ and thus 
$$
\left|\log \left|\frac{Z(y)+b}{Z(y)-a}\right|\right| = \left|\log \left|1 + \frac{a+b}{Z(y)-a}\right| \right| =  \left| \log \left(1 + \frac{a+b}{Z(y)-a}\right) \right| \leq  \frac{2|a+b|}{Z(y) - |a|}   , \quad \text{ as } y \to + \infty.
$$
Since $\frac1{Z}$ is decreasing and in $L^1([1,+\infty))$, we deduce that $y \mapsto \log \left|\frac{Z(y)+b}{Z(y)-a}\right|$ is integrable near $+\infty$ (integrability on compacts sets has already been proved for $F$). Moreover, its integral near infinity is continuous in $(a,b)$ by dominated convergence. This concludes the proof of the lemma.
\enp
We now recall a classical representation theorem for the modulus of entire functions of exponential type, which will be crucial in the proof of Lemma~\ref{lmanalysecomplexe} below.
\begin{theorem}\cite[Theorem p56]{Koosis:bookI}
\label{thmKoosis}
Let $f (z)$ be entire and of exponential type and suppose that
\bna
\int_{-\infty}^{+\infty}\frac{\log_{+} |f (x)|}{1+x^{2}}dx < +\infty ,
\ena
where $\log_+(t) = \max\{0, \log(t)\}$. Denote by $\{\lambda_{n}\}$, the set of zeros of $f (z)$ in $\Im(z) > 0$ (repetitions according to
multiplicities), and put
\bna
A = \limsup_{y\to +\infty} \frac{\log |f (iy)|}{y}
\ena
Then, for $\Im(z) > 0$,
\bna
\log|f(z)|=A\Im z+\sum_{n=1}^{+\infty}\log \left|\frac{1-z/\lambda_{n}}{1-z/\overline{\lambda_{n}}}\right|+\frac{1}{\pi}\int_{-\infty}^{+\infty}\frac{\Im z}{|z-\tau|^2}\log|f(\tau)|d\tau.
\ena
\end{theorem}
\bnp[Proof of Lemma~\ref{lmanalysecomplexe}]
We now prove the main statement of the lemma. We apply Theorem \ref{thmKoosis} at the point $ix_\eps$ . Given that $x_\eps \to x$, we may assume that $x_\eps >0$, and have 
\bnan
\label{represent}
\log|g(ix_\eps)|=\sum_{\ell=1}^{+\infty}\log \left|\frac{ix_\eps-a_{\ell}}{ix_\eps-\overline{a}_{\ell}}\right|+\sigma x_\eps+\frac{x_\eps}{\pi}\int_{-\infty}^{+\infty}\frac{\log|g(\tau)|}{|\tau-ix_\eps|^2}d\tau ,
\enan
where $(a_{\ell})_{\ell \in \N}$ is the sequence of zeros of $g$ in $\C_+ := \{z\in \C, \Im(z)>0\}$ (repeated according to multiplicities).

We first estimate the third term in the right handside of~\eqref{represent} using Assumption~\ref{e:real-line}, as 
\bnan
\label{e:est-real}
\frac{x_\eps}{\pi}\int_{-\infty}^{+\infty}\frac{\log|g(\tau)|}{|\tau-ix_\eps|^2}d\tau\leq C_{\R}\frac{x_\eps}{\pi}\int_{-\infty}^{+\infty}\frac{1}{\tau^2+x_\eps^2}d\tau=\frac{C_{\R}}{\pi}\int_{-\infty}^{+\infty}\frac{1}{t^2+1}d\tau=C_{\R} .
\enan
The estimate of the first term in the right handside of~\eqref{represent} is more complicated. First, we notice that since $ix_\eps$ and $a_\ell$ are in $\C_+$, we have $|ix_\eps-a_{\ell}|\leq |ix_\eps-\overline{a}_{\ell}|$ and thus $\log \left|\frac{ix_\eps-a_{\ell}}{ix_\eps-\overline{a}_{\ell}}\right| \leq 0$ for all $\ell \in \N$. Therefore, since $\left\{ib_k, k\in\N\right\}\subset \left\{a_{\ell},\ell\in \N\right\} $, we deduce
\bnan
\label{e:im-rac}
\sum_{\ell\in \N} \log \left|\frac{ix_\eps-a_{\ell}}{ix_\eps-\overline{a}_{\ell}}\right|\leq \sum_{k\in \N}\log \left|\frac{x_\eps-b_k}{x_\eps+b_{k}}\right|.
\enan
We can also assume without loss of generality that $x_\eps\neq b_k$ for all $k \in \N$ (otherwise, the left handside in~\eqref{e:estim-imagg} is $-\infty$ and~\eqref{e:estim-imagg} holds true), and that the sequence $(b_k)_{k\in \N} \in (0,\infty)^\N$ is an increasing sequence. Denote then $N=N(\eps)$ the integer such that 
\begin{align}
\label{e:encad-x}
\cdots \leq b_{N-1} \leq b_N < x_\eps < b_{N+1} \leq b_{N+2}  \leq \cdots .
\end{align}
Notice that since $x>0$,  we have $N(\e)\to +\infty$ as $\eps\to 0$ (see \eqref{e:conv-eN} below for a more precise estimate).
We are thus left to study
\begin{align}
\label{e:Spm}
\sum_{k\in \N}\log \left|\frac{x_\eps-b_k}{x_\eps+b_{k}}\right| = S_\leq +S_>  , \quad \text{ with } \quad S_\leq :=\sum_{k\leq N}\log \left(\frac{x_\eps-b_k}{x_\eps+b_{k}}\right) , \quad S_>  := \sum_{k> N}\log \left(\frac{b_k-x_\eps}{x_\eps+b_{k}}\right) .
\end{align}
Using again that all terms in the sum are nonpositive, together with~\eqref{e:bkZR} and the fact that the functions $s\mapsto \frac{x_\eps-s}{x_\eps+s}$ (resp. $s\mapsto \frac{s-x_\eps}{s+x_\eps}=1-2\frac{x_\eps}{s+x_\eps}$) are decreasing (resp. increasing if $x_\eps>0$), we obtain respectively
\begin{align}
\label{e:Sp}
S_\leq & \leq \sum_{D+1\leq k\leq N}\log \left(\frac{x_\eps-b_k}{x_\eps+b_{k}}\right)\leq\sum_{D+1\leq k\leq N}\log \left(\frac{x_\eps-Z(\e k-D\e)+R(\e)}{x_\eps+Z(\e k-D\e)-R(\e)}\right)  ,\\
\label{e:Sm}
S_> & \leq \sum_{k> N}\log \left(\frac{b_k-x_\eps}{x_\eps+b_{k}}\right)\leq \sum_{k> N}\log \left(\frac{Z(\e k+D\e)+R(\e)-x_\eps}{x_\eps+Z(\e k+D\e)+R(\e)}\right) .
\end{align}
 Note also that $b_k< x_\eps$ implies $Z(\e k-D\e)-R(\e)\leq x_\eps$ and $ x_\eps-R(\e)>0$ (for $\e$ small enough), so the first expression makes sense (the same applies for the other term). We may rewrite these two inequalities as 
\bna
S_\leq\leq \sum_{D+1\leq k\leq N}f_\leq(k) , \qquad 
S_>\leq \sum_{k> N}f_>(k) , 
\ena
with  
$$
f_\leq(s)=\log \left(\frac{x_\eps-Z(\e s-D\e)+R(\e)}{x_\eps+Z(\e s-D\e)-R(\e)}\right) , \quad  f_>(s)=\log \left(\frac{Z(\e s+D\e)+R(\e)-x_\eps}{x_\eps+Z(\e s+D\e)+R(\e)}\right) .
$$
Note then  that the function $f_\leq$ is negative decreasing, whereas $f_>$ is negative increasing to zero.
As a consequence, we have 
\begin{align}
S_\leq & \leq  \sum_{k=D+1}^{N}f_\leq(k)  \leq \int_{D}^{N}f_\leq(s)ds = \int_{0}^{N-D}\log \left(\frac{x_\eps-Z(\e s)+R(\e)}{x_\eps+Z(\e s)-R(\e)}\right) ds \nonumber \\
&=\frac{1}{\e}\int_{0}^{\e (N-D)}\log \left|\frac{x_\eps-Z(y)+R(\e)}{x_\eps+Z(y)-R(\e)}\right| dy \nonumber \\
&=-\frac{1}{\e}I_{N,\e}^{\leq}  , \quad \text{with} \quad I_{N,\e}^{\leq} :=\int_{0}^{\e (N-D)}\log \left|\frac{x_\eps+Z(y)-R(\e)}{x_\eps-Z(y)+R(\e)}\right| dy .
\label{e:Spp}
\end{align}
Similarly, we have
\begin{align}
S_> & \leq \sum_{k=N+1}^\infty f_>(k) \leq \int_{N+1}^{+\infty}f_>(s)ds = \int_{N+1+D}^{+\infty}\log \left(\frac{Z(\e s)+R(\e)-x_\eps}{x_\eps+Z(\e s)+R(\e)}\right)ds \nonumber \\
&= \frac{1}{\e}\int_{\e(N+1+D)}^{+\infty}\log \left|\frac{Z(y)+R(\e)-x_\eps}{x_\eps+Z(y)+R(\e)}\right|dy \nonumber \\
&= -\frac{1}{\e}I_{N,\e}^> , \quad \text{with} \quad I_{N,\e}^> :=\int_{\e(N+1+D)}^{+\infty}\log \left|\frac{x_\eps+Z(y)+R(\e)}{Z(y)+R(\e)-x_\eps}\right|dy .
\label{e:Smm}
\end{align}
Combining~\eqref{represent}-\eqref{e:est-real}-\eqref{e:im-rac}-\eqref{e:Spm}-\eqref{e:Sp}-\eqref{e:Sm}-\eqref{e:Spp}-\eqref{e:Smm}, we have obtained so far that 
\begin{align}
\label{e:comp-alm-fin}
\log|g(ix_\eps)| \leq \sigma x_\eps +C_\R -\frac1{\eps}(I_{N,\e}^{\leq}  + I_{N,\e}^>) ,
\end{align}
and it only remains to study $I_{N,\e}^{\leq} , I_{N,\e}^>$.

\medskip
Note that~\eqref{e:bkZR} applied to $N$ and $N+1$, and the definition of $N$ in~\eqref{e:encad-x} yield
\bna
Z(\e (N-D))-R(\e)\leq b_N < x_{\e}<b_{N+1}\leq Z(\e (N+1+D))+R(\e) .
\ena 
In particular, by continuity of $Z$, this implies that $Z( \e N) = Z( \e N(\eps))$ converges to $x$ as $\eps\to 0^+$, and hence
\begin{align}
\label{e:conv-eN}
\eps N = \eps N(\eps) \to Z^{-1}(x) , \quad \text{ as } \eps \to 0^+.
\end{align}
Next, we define the function $h_{\e} :  \R^+ \to \ovl{\R}$ by 
\bna
h_{\e}(y) :=&\log \left|\frac{x_\eps+Z(y)-R(\e)}{x_\eps-Z(y)+R(\e)}\right|  = \log \left(\frac{x_\eps+Z(y)-R(\e)}{x_\eps-Z(y)+R(\e)}\right) , &\textnormal{ for }Z(y)\leq x_\eps-R(\e) ,\\
h_{\e}(y) :=&\log \left|\frac{x_\eps+Z(y)+R(\e)}{x_\eps-Z(y)-R(\e)}\right| = \log \left(\frac{x_\eps+Z(y)+R(\e)}{-x_\eps+Z(y)+R(\e)}\right) , &\textnormal{ for }Z(y)\geq x_\eps+R(\e) ,\\
h_{\e}(y) :=&0 ,&\textnormal{ otherwise.}
\ena
Recalling the definition of $I_{N,\e}^{\leq } , I_{N,\e}^>$ in~\eqref{e:Spp}-\eqref{e:Smm}, we now have
\begin{align}
\label{e:conv-I-I}
\left| I_{N,\e}^{\leq }+I_{N,\e}^> - \int_0^{+\infty}h_{\e}(y)dy \right| & \leq  \left| \int_{\e (N-D)}^{Z^{-1}(x_\eps-R(\e))}\log \left|\frac{x_\eps+Z(y)-R(\e)}{x_\eps-Z(y)+R(\e)}\right| dy\right| \nonumber \\
 & \quad + \left| \int_{Z^{-1}(x_\eps+R(\e))}^{\e(N+1+D)}\log \left|\frac{x_\eps+Z(y)+R(\e)}{x_\eps-Z(y)-R(\e)} \right|\right| ,
\end{align}
and now examine the convergence of the different terms involved. 
We shall prove that 
\begin{align}
\label{e:heps-to-I}
 \int_0^{+\infty}h_{\e}(y)dy \to I(x), \quad \text{ as } \eps\to 0^+ ,
\end{align}
and that the right handside of~\eqref{e:conv-I-I} converges to zero. This, together with~\eqref{e:comp-alm-fin} will then yield~\eqref{e:estim-imagg}, concluding the proof of the lemma.

Let us now prove~\eqref{e:heps-to-I} by splitting the integral into the two intervals: 
$$
 \int_0^{+\infty}h_{\e}(y)dy = \int_0^{Z^{-1}(x-R(\e))}  \log \left|\frac{x_\eps+Z(y)-R(\e)}{x_\eps-Z(y)+R(\e)}\right|  dy +  \int_{Z^{-1}(x+R(\e))}^{+\infty}  \log \left|\frac{x_\eps+Z(y)-R(\e)}{x_\eps-Z(y)+R(\e)}\right|  dy  .
$$
Lemma~\ref{l:continu-I} implies the following convergence of the two integrals as $\eps\to 0^+$:
$$
 \int_0^{+\infty}h_{\e}(y)dy \to \int_0^{Z^{-1}(x)} \log \left|\frac{x+Z(y)}{x-Z(y)}\right|  dy +  \int_{Z^{-1}(x)}^{+\infty}  \log \left|\frac{x+Z(y)}{x-Z(y)}\right|  dy =  I(x), 
$$
which is~\eqref{e:heps-to-I}.

We finally consider the right handside of~\eqref{e:conv-I-I}. According to~\eqref{e:conv-eN}, both endpoints of these two intervals converge to $Z^{-1}(x)$. Using again Lemma~\ref{l:continu-I}, this implies that the right handside of~\eqref{e:conv-I-I} converges to zero, which concludes the proof of the lemma.
\enp

\subsection{Upper bound: Proof of Theorem~\ref{thmcontrol1D}}
\label{s:proof-upper-bound}
In this section, we give a proof of Theorem~\ref{thmcontrol1D}.  
In particular, we assume that $\f \in C^\infty([0,L])$, that Items~\ref{A1}--\ref{A4} in Assumption~\ref{assumptions} are satisfied, and that $\q= \frac{\f''}{2}$. These assumptions are made so that to apply the spectral results of Theorem~\eqref{t:Peps-spect-Allib}, deduced from~\cite{Allibert:98}.

According to Definition~\ref{d:def-0} and Lemma~\ref{lemequivecontrol}, null controllability of~\eqref{e:1Dcontrol-problemyn} in time $T$ is equivalent to having for any $y_0\in L^2(0,L)$, the existence of $h\in L^2(0,T)$ such the solution $v$ to~\eqref{e:heat-control-eps-Witt} satisfies $v(T,\cdot) = 0$.

With the duality~\eqref{e:duality-conjugated} this can equivalently be formulated as: having for any $y_0\in L^2(0,L)$, the existence of $h \in L^2(0,T)$ such that for all $\zeta_*\in L^2(0,L)$ we have
\begin{align}
\label{e:controllability}
&0 = \left( \zeta(T) , v(0) \right)_{L^2(0,L)}  + \int_0^T  \eps \partial_x \zeta(t,0) v(T-t, 0) dt =0 , \quad \text{ or equivalently, } \nonumber\\ 
& 0 =  \int_0^L e^{-\frac{\f(x)}{2\eps}} y_0 (x) \big( e^{-\frac{T}{\eps}P_\eps}\zeta_*\big)(x) dx  +\eps \int_0^T e^{-\frac{\f(0)}{2\eps}} h(T-t) \d_x (e^{-\frac{t}{\eps}P_\eps}\zeta_*)|_{x=0} dt ,
\end{align}
where $\zeta$ solves~\eqref{e:heat-transp-eps-Witt} with $\zeta(0,x) =\zeta_*(x) = e^{\frac{\f(x)}{2\eps}}u_0(x)$ and $v$ solves~\eqref{e:heat-control-eps-Witt}.

With domain $H^2\cap H^1_0$ (i.e. Dirichlet boundary conditions), the operator $P_\eps$ is selfadjoint on $L^2$, with compact resolvent. We introduce a Hilbert basis $(\varphi_\ell^\eps)_{\ell \in \N}$ of $L^2(0,L)$ such that 
$$
P_\eps \varphi_\ell^\eps = \lambda_\ell^\eps \varphi_\ell^\eps , \quad \text{ with } \quad \varphi_\ell^\eps(0) = \varphi_\ell^\eps(L)=0 , \quad \lambda_\ell^\eps \leq \lambda_{\ell+1}^\eps ,
$$
as in Theorem~\ref{t:Peps-spect-Allib}.

We will need the following intermediate result, which we state on a time interval $(0,\tau)$ instead of $(0,T)$ for in the proof of Theorem \ref{thmcontrol1D}, this control will be only used in part of the whole time interval $(0,T)$.

\begin{proposition}
\label{propcontrolzerointerm}
Under the assumptions of Theorem \ref{thmcontrol1D}, fix $\delta>0$. Then, there exists $\e_{0}$ and $C>0$ so that for any $0<\e<\e_{0}$, $0<\tau<\delta^{-1}$ and $v_0\in L^{2}(0,L)$, there exists a control $h\in L^{2}(0,\tau)$ to zero of
\begin{equation*}
\left\{
\begin{array}{rl}
(\eps \d_t + P_\eps ) v = 0 ,& (t,x)\in (0,\tau) \times (0,L) , \\
v (t,0) = e^{-\frac{\f(0)}{2\eps}}h(t) , \quad v(t,L) = 0 ,&  t \in (0,\tau) , \\
v(0,x) = v_0 (x),&  x \in (0,L) .
\end{array}
\right.
\end{equation*}
 with the control cost
\begin{align}
\label{estimz}
\nor{h}{L^2(0,\tau)}^2 \leq  \frac{C}{\e^6 \tau^3}e^{\frac{4 T_1^2 +\delta}{\eps \tau}+ \frac{\f(0)}{\eps}}\sum_{n\in \N}\frac{\lambda_n^\eps}{|(\varphi_n^\eps)'(0)|^2}  \left|\int_0^L v_0 (x)\varphi_n^\eps(x) dx \right|^2 .
\end{align}
\end{proposition}
The proof consists in solving the moment problem obtained by testing~\eqref{e:controllability} with $\zeta_*$ ranging in a basis of eigenfunctions of $P_\eps$. It also relies on properties on $P_\eps$ described in Theorem~\ref{t:Peps-spect-Allib} (recall that our assumptions imply $\qf=0$).

\bnp[Proof of Proposition~\ref{propcontrolzerointerm}]
According to~\eqref{e:controllability} controlling $v_0$ to zero in time $\tau$ is equivalent to having existence of $\tilde{z} = h (\tau - \cdot)\in L^2(0,\tau)$ such that for all $\ell \in \N$, 
\bna
0 = e^{-\frac{\tau}{\eps}\lambda_\ell^\eps}  \int_0^L v_0 (x)\varphi_\ell^\eps(x) dx  +\eps e^{-\frac{\f(0)}{2\eps}} (\varphi_\ell^\eps)'(0) \int_0^\tau \tilde{z}(t)e^{-\frac{t}{\eps}\lambda_\ell^\eps} dt .
\ena
The idea of the moment method for finding such $\tilde{z}(t)$ solving
\bnan
\label{defalphaeps}
 \int_0^\tau \tilde{z}(t)e^{-\frac{t}{\eps}\lambda_\ell^\eps} dt = \alpha_\ell^\eps ,
  \quad \text{with} \quad \alpha_\ell^\eps =- \frac{e^{-\frac{\tau}{\eps}\lambda_\ell^\eps} }{\eps e^{-\frac{\f(0)}{2\eps}} (\varphi_\ell^\eps)'(0)}  \int_0^L v_0 (x)\varphi_\ell^\eps(x) dx  ,
\enan
is to construct it as a sum of biorthogonal functions $(\Psi_j^\eps)_{j\in \N} \in L^2(0,\tau)^\N$, namely 
\begin{equation}
\label{e:biorth}
\tilde{z}(t) = \sum_{j\in \N}\alpha_j^\eps  \Psi_j^\eps(t), \quad \text{ with } \int_0^\tau \Psi_j^\eps (t)e^{-\frac{\lambda_\ell^\eps}{\eps} t} dt = \delta_{j\ell} .
\end{equation}
Denoting $z(t):=\eps^{-1}\tilde{z}(t/\eps )$, $\psi_j^\eps (t)=\eps^{-1} \Psi_j^\eps( t/\e)$ and $\beta_\ell^\eps:=\e^{-1}\sqrt{\lambda_\ell^\eps}$, Equation~\eqref{e:biorth} is equivalent to 
\begin{equation}
\label{e:biorthter}
z(t) = \sum_{j\in \N} \alpha_j^\eps  \psi_j^\eps(t), \quad \text{ with } \int_0^{ \eps \tau} \psi_j^\eps (s)e^{-(\beta_\ell^\eps)^{2}s} ds = \delta_{j\ell} .
\end{equation}
For $\delta>0$, Theorem \ref{t:Peps-spect-Allib} yields existence of $\eps_0, \gamma >0$ and $N\in \N$ (all depending on $\delta$) such that for all $0<\eps\leq \eps_0$, the sequence $(\beta_\ell^\eps)_{\ell \in \N}$ satisfies
\bnan
\label{gapbeta}
\beta_{\ell+1}^\eps -\beta_\ell^\eps \geq  \frac{2\pi}{T_1+ \delta/2}, \quad \text{ for all } \ell \geq N, \\
\beta_{\ell+1}^\eps -\beta_\ell^\eps  \geq   \gamma>0 ,\quad \text{for all }  \ell\in\N , \nonumber
\enan
where $T_1$ is defined in~\eqref{e:def-TT1}.
Proposition \ref{propgapheat} with $\gamma_\infty = \frac{2\pi}{T_1+ \delta/2}$ yields existence of $C,\eps_0>0$ such that for all $\e<\eps_0$, setting $S_{\delta} :=\frac12(T_1+ \delta)$, we can find a sequence $(\psi_j^\eps)_{j\in \N}$ satisfying~\eqref{e:biorthter} with 
\begin{align}
\label{e:estim-z}
\nor{z}{L^2(0,\e \tau)}^2\leq \frac{C}{(\e \tau)^3}e^{\frac{(16+\delta)S_{\delta}^{2}}{\eps \tau}}\sum_{\ell\in \N} (\beta_\ell^\eps)^{2} e^{2(\beta_\ell^\eps)^{2}\eps \tau} |\alpha_\ell^\eps|^2
=\frac{C}{(\e \tau)^3}e^{\frac{(16+\delta)S_{\delta}^{2}}{\eps \tau}}\sum_{\ell\in \N} \frac{\lambda_\ell^\eps}{\eps^2} e^{2\tau\frac{\lambda_\ell^\eps }{\e}} |\alpha_\ell^\eps|^2
\end{align}
Choosing now the numbers $\alpha_j^\eps$ as in the second part of~\eqref{defalphaeps}, the function $z$ satisfies the first part of~\eqref{defalphaeps}, and hence is a null-control for $v_0$. Moreover, we can now estimate the last term in~\eqref{e:estim-z} as 
\begin{align*}
\sum_{\ell \in \N} \lambda_\ell^\eps  e^{2\tau\frac{\lambda_\ell^\eps }{\e}} |\alpha_\ell^\eps|^2
&\leq  \frac{1}{\eps^2 e^{-\frac{\f(0)}{\eps}}}\sum_{\ell\in \N}\frac{\lambda_\ell^\eps}{|(\varphi_\ell^\eps)'(0)|^2}  \left|\int_0^L v_0 (x)\varphi_\ell^\eps(x) dx \right|^2
\end{align*}
Combined with~\eqref{e:estim-z}, this yields
\begin{align*}
\nor{z}{L^2(0,\e \tau)}^2\leq \frac{C}{\e^7 \tau^3}e^{\frac{(16+\delta)S_{\delta}^{2}}{\eps \tau}}\frac{1}{e^{-\frac{\f(0)}{\eps}}}\sum_{n\in \N}\frac{\lambda_n^\eps}{|(\varphi_n^\eps)'(0)|^2}  \left|\int_0^L v_0 (x)\varphi_n^\eps(x) dx \right|^2 .
\end{align*}
This gives the result recalling that $S_{\delta}=\frac12(T_1+ \delta)$ and $\nor{h}{L^2(0,\tau)}^2=\nor{\widetilde{z}}{L^2(0,\tau)}^2=\e\nor{z}{L^2(0,\e \tau)}^2$ (up to changing the notation for $\delta$).
\enp

We are now in position to prove Theorem~\ref{thmcontrol1D}. We first take advantage of the natural parabolic dissipation \cite{LR:95,Lissy:12,LL:17hypo}, and then use the control function constructed in Proposition~\ref{propcontrolzerointerm}.

\bnp[Proof of Theorem \ref{thmcontrol1D}]
We construct a control function $h(t)$ for Equation~\eqref{e:heat-control-eps-Witt} under the following form 
\begin{itemize}
\item $h=0$ on $[0,m \TO]$ and use dissipation;
\item $h$ as constructed in Proposition~\ref{propcontrolzerointerm} on the interval $[m\TO,\TO]$ (instead of $[0,\tau]$; this is possible since the equation is invariant by translations in time).
\end{itemize}
At time $m\TO$ the solution of~\eqref{e:heat-control-eps-Witt} is thus given by 
\bna
v(m \TO,x)=\sum_{n\in\N} e^{-\frac{\lambda_n^{\e}}{\eps} m\TO}v_n\varphi_n^{\e}(x) ,
\ena
where $v_n=\int_0^Le^{-\frac{\f(x)}{2\eps}} y_0 (x)\varphi_n^{\e}dx$. Moreover, using the Cauchy-Schwarz estimate, we have
\bna
|v_n|\leq \nor{y_0 }{L^2(0,L)}\nor{e^{-\frac{\f(x)}{2\eps}} \varphi_n^{\e}}{L^2(0,L)}.
\ena
We then take this function $v(m \TO,\cdot)$ as an initial condition for the control problem on $[m\TO,\TO]$.
On this interval, we use the control function $h$ furnished by Proposition \ref{propcontrolzerointerm}. It satisfies Estimate~\eqref{estimz} which reads
\begin{align*}
\nor{h}{L^2(m\TO,\TO)}^2 & \leq  \frac{C}{\e^6 (1-m)^3T^3}e^{\frac{4 T_1^2 +\delta}{\eps (1-m)T}+ \frac{\f(0)}{\eps}}\sum_{n\in \N}\frac{\lambda_n^\eps}{|(\varphi_n^\eps)'(0)|^2}  \left|\int_0^L v(mT,x)\varphi_n^\eps(x) dx \right|^2\\
&\leq   \frac{C_m}{\e^6T^3}e^{\frac{4 T_1^2 +\delta}{\eps (1-m)T}+ \frac{\f(0)}{\eps}}\sum_{n\in \N}\frac{\lambda_n^\eps} {|(\varphi_n^\eps)'(0)|^2}  e^{-2\frac{\lambda_n^{\e}}{\e} m\TO} |v_n|^2\\
&\leq \nor{y_0 }{L^2(0,L)}^2\frac{C_m}{\e^4T^3}e^{\frac{4 T_1^2 +\delta}{\eps (1-m)T}+ \frac{\f(0)}{\eps}}\sum_{n\in \N}\frac{\lambda_n^{\e}} {|\eps(\varphi_n^\eps)'(0)|^2}  e^{-2\frac{\lambda_n^{\e}}{\e} m\TO}\nor{e^{-\frac{\f(x)}{2\eps}} \varphi_n^{\e}}{L^2(0,L)}^2\\
&\leq \nor{y_0 }{L^2(0,L)}^2 \frac{C_m}{\e^4T^3}e^{\frac{4 T_1^2 +\delta}{\eps (1-m)T}+ \frac{\f(0)}{\eps}} A_{\e}B_{\e} ,
\end{align*}
where we have denoted for $0<\theta<1$ small
\begin{align*}
A_{\e}=\sum_{n\in \N}e^{-2\frac{\lambda_n^{\e}}{\e} \theta m\TO}, \quad 
B_{\e} = \sup_{n\in \N}\frac{\lambda_n^{\e}} {|\eps(\varphi_n^\eps)'(0)|^2}  e^{-2\frac{\lambda_n^{\e}}{\e}(1-\theta) m\TO}\nor{e^{-\frac{\f(x)}{2\eps}} \varphi_n^{\e}}{L^2(0,L)}^2.
\end{align*}
We estimate $A_{\e}$ and $B_{\e}$ in Lemma \ref{lmAe} and \ref{lmBe} that we state below. Combined with the previous estimate, it gives (for any $\delta,T_{\max},m>0$, existence of $C_m,\eps_0>0$ such that for all $T\in(0,T_{\max})$, $m\in (0,1)$, $\theta\in(0,1)$ and all $\e\in(0,\e_0)$)
\bnan
\label{e:cost-born-sup-interm}
\nor{h}{L^2(0,\TO)}^2 = \nor{h}{L^2(m\TO,\TO)}^2\leq \frac{C_m}{\theta \eps^4 \TO^3 }e^{\frac{2D(m)}{\e}}\nor{y_0 }{L^2(0,L)}^2 , 
\enan
with (recalling that $W_E = d_{A,E}+\frac{\f}{2}$)
\bna
D(m) = \frac{2 T_1^2}{(1-m)\TO}+\sup_{E\geq E_0}\left[W_E(0)-E (1-\theta) m\TO -\min_{[0,L]} W_E\right]  +\mathsf{C}\delta ,
\ena
for a constant $\mathsf{C}$ depending on $m$, $\TO$, $T_{1}$.  This proves Estimate~\eqref{e:cost-born-sup} in Theorem \ref{thmcontrol1D} after taking $\e_0$ small enough to absorb the polynomial loss, up to changing $\delta$.  Estimate~\eqref{e:Tunif-born-sup} in follow from optimizating in $m$, see Section \ref{calculusI}. We take $\theta=\delta$ and first downgrade the exponential part by using $E\geq E_0 = \frac{| \f'(\xzero)|^2}{4}$ to 
\begin{align*}
D(m)&\leq \frac{2T_{1}^{2}}{ (1-m)\TO}- (1-\theta) \frac{| \f'(\xzero)|^2}{4}m \TO +\sup_{E\geq E_0}\left[W_E(0)-\min_{[0,L]} W_E\right]+\mathsf{C}\delta\\
&\leq \frac{2T_{1}^{2}}{ (1-m)\TO}-  \frac{| \f'(\xzero)|^2}{4}m \TO +\sup_{E\in V(\M)}\left[W_E(0)-\min_{[0,L]} W_E\right]+\widetilde{\mathsf{C}}\delta .
\end{align*}
where we have noticed that $W_{E}=\f/2$ for $E\geq \max_{\M}V$. As a consequence of this together with~\eqref{e:cost-born-sup-interm}, we deduce that  we can infer $T > T_{\unif}$ if 
$$
G(T) := \min_{m\in [0,1)} \left( \frac{2T_{1}^{2}}{ (1-m)\TO}-  \frac{| \f'(\xzero)|^2}{4}m \TO +\sup_{E\in V(\M)}\left[W_E(0)-\min_{[0,L]} W_E\right] \right) <0 .
$$
Lemma \ref{lmcalctmin} then concludes the proof of~\eqref{e:Tunif-born-sup}, and that of Theorem~\ref{thmcontrol1D}.
\enp
It remains to prove the two Lemmata estimating $A_\e$ and $B_\e$.
\begin{lemma}
\label{lmAe}
Under the assumptions of Theorem \ref{thmcontrol1D}, given $T_{\max}>0$, there are $C,\eps_0>0$ such that for all $T\in(0,T_{\max})$, $m\in (0,1)$, $\theta\in(0,1)$ and all $\e\in(0,\e_0)$, 
\bna
\left|A_{\e}\right|\leq \frac{C}{\theta mT} .
\ena
\end{lemma}
\bnp
Item \ref{i:simplicity}, Item \ref{itemallibertgap} and estimate \eqref{e:gap-lambda}, each one to the version adapted to Theorem \ref{t:Peps-spect-Allib} give respectively for $\e$ small enough and for some $\gamma_{2}>0$
\bna
\lambda_0^{\e}\geq  \frac{|\f'(\xzero)|^2}{4}-\delta= : \Lambda_{\delta}>0, \quad \text{ together with }
 \quad \lambda_{n+1}^\eps -\lambda_n^\eps \geq \eps \gamma_{2} >0 .
\ena
In particular, this implies $\lambda_{n}^{\e}\geq \Lambda_{\delta}+n\eps \gamma_{2} $ and we can estimate
\bna
\left|A_{\e}\right|\leq \sum_{n\in \N}e^{-2n \gamma_{2}   \theta m\TO} = \frac{1}{1-e^{-2\gamma_{2}   \theta mT}} , 
\ena
whence the result by the mean value theorem.
\enp
\begin{lemma}
For any $\delta>0$, there exists $\e_{0}$ so that for any $0<\e<\e_{0}$, we have
\label{lmBe}
\bna
\left|B_{\e}\right|\leq C e^{2\frac{F+\delta}{\e}} , \quad \text{with} \quad 
F=\sup_{E\geq E_0}\left(d_{A,E}(0)-E (1-\theta) m\TO -\min_{[0,L]} W_E\right) .
\ena
\end{lemma}

\bnp
For any $n\in \N$ and $\e>0$, we call $E=\lambda_n^{\e}$.
Estimate~\eqref{e:bound-obs} in Theorem \ref{t:allibert-lower-uniform} yields $\frac{\eps}{\sqrt{E}}|(\varphi_n^\eps)'(0)|  \geq e^{-\frac{1}{\eps}(d_{A,E}(0) + \delta)}$, uniformly in $E$ and $\e$ (for $\e$ small enough). For the second term, we simply write $e^{-2\frac{\lambda_n^{\e}}{\e}(1-\theta) m\TO}= e^{-2\frac{E}{\e}(1-\theta) m\TO}$. The last term is estimated thanks to the Agmon type estimate \eqref{e:agmon-H1} of Theorem \ref{t:uniform-agmon} as
\bna
\nor{e^{-\frac{\f(x)}{2\eps}} \varphi_n^{\e}}{L^2(0,L)}= \nor{e^{-\frac{ W_{E}}{\eps}} e^{\frac{d_{A,E}}{\eps} }\varphi_n^{\e}}{L^2(0,L)}\leq e^{-\frac{\min_{[0,L]} W_{E}}{\eps}} \nor{e^{\frac{d_{A,E}}{\eps} }\varphi_n^{\e}}{L^2(0,L)}\leq  e^{-\frac{\min_{[0,L]} W_{E}}{\eps}} e^{\frac{\delta}{\eps}}.
\ena
The combination of these three estimates gives
$$
\frac{\lambda_n^{\e}} {|\eps(\varphi_n^\eps)'(0)|^2}  e^{-2\frac{\lambda_n^{\e}}{\e}(1-\theta) m\TO}\nor{e^{-\frac{\f(x)}{2\eps}} \varphi_n^{\e}}{L^2(0,L)}^2
\leq e^{\frac{2}{\eps}(d_{A,E}(0) + \delta)} e^{-2\frac{E}{\e}(1-\theta) m\TO}e^{-2 \frac{\min_{[0,L]} W_{E}}{\eps}} e^{2\frac{\delta}{\eps}}.
$$
Recalling (see Item \ref{i:simplicity} in Theorem \ref{t:Peps-spect-Allib}) that $E\geq E_0-C\e^2$ and taking the supremum over all $E=\lambda_n^{\e}$ yields the sought result (up to a loss in $\delta$, we can take the supremum in $E\geq E_0$).
\enp

\section{Explicit computations of the various bounds on an example}
\label{s:explicit-comput}

In this section, we provide with some explicit computation of the different bounds we computed in the main part of the paper for concrete examples of functions $\f$. 
For symmetry reasons, we shift the problem and consider the interval $(-L/2,L/2)$ controlled at the point $-L/2$.

For $a>0$, $M>0$, we choose 
\bnan
\label{deff1D} \f_{M,a}^{\pm}(x) & =\pm\int_0^x\sqrt{a^2s^2+M^{2}}ds=\pm \frac{M^{2}}{a} \int_0^{\frac{a x}{M}}\sqrt{y^2+1}dy =\pm\frac{x}{2}\sqrt{a^{2}x^{2} +M^{2}}\pm\frac{ M^{2}}{2a}\arcsinh(\frac{a x}{M}) ,
\enan 
where we have used the identity 
$$
\int_0^y \sqrt{1+t^2} dt = \frac12 \left( y\sqrt{y^2+1}+\arcsinh(y) \right) , \quad y \in \R .
$$
With this choice, we have $\f^{\pm '}_{M,a}(x)=\pm\sqrt{a^2x^2+M^{2}}$ on $(-L/2,L/2)$. The potential 
\begin{align}
\label{e:example-potential}
V(x) = \frac{|\f^{\pm '}_{M,a}(x)|^2}{4}=\frac{a^2x^2+M^{2}}{4}
\end{align}
reaches its minimum at the point $\xzero = 0\in (-L/2,L/2)$. We have chosen this example for the relative simplicity of the computations and because the formal limit when $a\to 0^+$ is the model with constant transport term, well studied in the literature~\cite{CG:05,Gla:10,Lissy:12,Lissy:14,Lissy:15,Darde:17,YM:19,YM:19bis}, and the limit $a\to +\infty$ proves that $T_{\unif}$ can be much larger than the control/flushing time for the transport equation with $T_{\f^{\pm '}_{M,a}} (\{-L/2\})$ (see \cite{LL:18vanish}). Indeed, $V$ is a constant plus a harmonic potential. The parameter $a$ will allow to stress the fact that the convexity is responsible for a concentration of some eigenfunctions close to the minimum, which is not the case for the ``flat potential'' $V=\frac{M^2}{4}$ corresponding to the more studied case $\f(x)=\pm Mx$~\cite{CG:05,Gla:10,Lissy:12,Lissy:14,Lissy:15,Darde:17,YM:19,YM:19bis}.
We now compute explicitly of the quantities involved in the statements of Proposition~\ref{p:transport-seul} and Theorems~\ref{t:estim-Cobs-expo-1D}--\ref{thmlower+-}--\ref{thmcontrol1D}.

\subsection{Computation of $T_{\f^{\pm '}_{M,a}}, T_1, T_{E,B},d_{A,E_{0}}$}
\label{subsectexplcit1D}

\begin{lemma}
\label{l:1D-explicit-limit}
For the function $\f=\f_{M,a}^{\pm}(x)$ defined in~\eqref{deff1D}, the minimal control time (or flushing time, depending on the sign) for the limit equation ($\eps=0$) is given by 
$$
T_{\f^{\pm '}_{M,a}} (\{-L/2\}) =  \frac{2}{a}\arcsinh\left(\frac{aL}{2M}\right).
$$
\end{lemma}
\begin{lemma}
\label{lm:TTEB}
Recalling the definitions~\eqref{e:def-phi}-\eqref{e:def-TT1}-\eqref{TEB}, for the function $\f=\f_{M,a}^{\pm}(x)$ defined in~\eqref{deff1D}, we have
\begin{align}
\label{f:T1}
T_1&=\frac{2\pi}{a}\sqrt{a^2L^2/4+M^2} = \pi \sqrt{L^2+4M^2/a^2} , \\ 
T_{E,B}&=\frac{1}{a}\int_{M^{2}/4}^{a^2L^2/16+M^{2}/4}\log \left|\frac{x+E+2B}{x-E}\right|dx+\frac{2}{a\pi}\int_{a^2L^2/16+M^{2}/4}^{+\infty}\log \left|\frac{x+E+2B}{x-E}\right|\arcsin\left(\frac{aL}{2\sqrt{4x-M^{2}}}\right)dx . \nonumber
\end{align}
Moreover, we have
\begin{equation*}
\begin{array}{ll}
T_1 \underset{a\to 0^{+}}{\rightarrow} +\infty ,  & \quad T_1 \underset{a\to +\infty}{\rightarrow} \pi L, \\
T_{E,B} \underset{a\to 0^{+}}{\rightarrow} \frac{L}{2\pi}\int_{0}^{+\infty}\log \left|\frac{y^2+M^{2}+4E+8B}{y^2+M^{2}-4E}\right|dy
 <+\infty , &  \quad T_{E,B} \underset{a\to +\infty}{\rightarrow} 0,
 \end{array}
\end{equation*}
with in particular $\lim_{a\to 0^{+}} T_{E_0,B}= \frac{L}{2}\sqrt{2M^{2}+8B}$ (case $E=E_{0}=\frac{M^2}{4}$).   
\end{lemma}
Note that the function in second integral in the expression of $T_{E,B}$ behaves like $(x- (a^2L^2/16+M^{2}/4)^{-1/2}$ near $a^2L^2/16+M^{2}/4 $, like $\log(|x-E|)$ near $E$ (if $E$ is in the interval) and like $x^{-3/2}$ at $+\infty$. Therefore, the integral is well defined.
\begin{lemma}
\label{lmAgmonexpl}
For the function $\f=\f_{M,a}^{\pm}(x)$ defined in~\eqref{deff1D}, the Agmon distance to the potential minimum is given by
$d_{A,E_{0}}(x) =\frac{ax^2}{4}.$
\end{lemma}
\bnp[Proof of Lemma~\ref{lmAgmonexpl}]This follows from the expression~\eqref{e:example-potential} of the associated potential and the direct computation
$$
d_{A,E_{0}}(x)= \left| \int_{0}^x \sqrt{\frac{\f'(y)^2}{4}  - \frac{\f'(0)^2}{4}} dy  \right|=\frac{1}{2} \left| \int_{0}^x a|y|~dy  \right|=\frac{ax^2}{4}.
$$
\enp
\bnp[Proof of Lemma~\ref{l:1D-explicit-limit}]
According to Proposition \ref{p:tpstransp}, we have the exact formula
\bna
T_{\f'^{\pm}_{M,a}} = \int_{-L/2}^{L/2} \frac{ds}{|\f'^{\pm}_{M,a}(s)|}=\int_{-L/2}^{L/2} \frac{dy}{\sqrt{ay^2+M^{2}}}=\frac{1}{a}\int_{-\frac{aL}{2M}}^{\frac{aL}{2M}} \frac{dy}{\sqrt{y^2+1}} 
= \frac{2}{a}\arcsinh\left(\frac{aL}{2M}\right),
\ena
where we have used 
\bna
\int_{-x}^{x} \frac{dy}{\sqrt{y^2+1}}= 2\arcsinh(x) , \quad x \in \R .
\ena
\enp

\bnp[Proof of Lemma \ref{lm:TTEB}]
According to~\eqref{deff1D}, we have $4V(x) = |\f_{M,a}^{\pm'}(x)|^2=a^2x^2+M^2$, and $V(x)=\lambda \Leftrightarrow a^2x^2+M^{2}=4\lambda  \Leftrightarrow x=\pm \frac{\sqrt{4\lambda-M^{2}}}{a}$  and it belongs to $[-L/2,L/2]$ if  $4\lambda-M^{2}\leq a^2L^{2}/4$. We rename $\Lambda=\Lambda(\lambda)=4\lambda-M^{2}$ so that \eqref{e:def-TT1} and 
$$
\int_{-t}^t \frac{dy}{\sqrt{1-y^2}}  = 2 \arcsin(t) , \quad t \in [-1,1]
$$
give
\begin{itemize}
\item if $0\leq \Lambda\leq a^2L^{2}/4$, we have $x_\pm(\lambda)=\pm\frac{\sqrt{4\lambda-M^{2}}}{a}=\pm \frac{\sqrt{\Lambda}}{a}$ and (as for the harmonic oscillator)
\bna
T(\lambda)=2\int_{-\frac{\sqrt{\Lambda}}{a}}^{\frac{\sqrt{\Lambda}}{a}}\sqrt{\frac{\Lambda+M^{2}}{\Lambda-a^2x^2}}dx=\frac{2}{a}\int_{-1}^{1}\sqrt{\frac{\Lambda+M^{2}}{1-y^2}}dy=2\pi\frac{\sqrt{\Lambda+M^{2}}}{a} = 4\pi \frac{\sqrt{\lambda}}{a};
\ena
\item if $a^2L^{2}/4 \leq \Lambda$, we have $x_\pm(\lambda)=\pm L/2$ so that 
\begin{align*}
T(\lambda)&=2\int_{-L/2}^{L/2}\sqrt{\frac{\Lambda+M^{2}}{\Lambda-a^2x^2}}dx=\frac{2}{a}\int_{-\frac{aL}{2\sqrt{\Lambda}}}^{\frac{aL}{2\sqrt{\Lambda}}}\sqrt{\frac{\Lambda+M^{2}}{1-y^2}}dy\\
&= \frac{4\sqrt{\Lambda+M^{2}}}{a}\arcsin\left(\frac{aL}{2\sqrt{\Lambda}}\right) = 4\arcsin\left(\frac{aL}{2\sqrt{4\lambda-M^2}}\right)\frac{2\sqrt{\lambda}}{a} .
\end{align*}
\end{itemize}
Coming back to~\eqref{e:def-TT1}, we finally obtain
\bna
T_{1}= \sup_{\lambda \geq E_0} T(\lambda) = \max\left(\sup_{0\leq \Lambda\leq a^2L^{2}/4} 2\pi\frac{\sqrt{\Lambda+M^{2}}}{a}; \sup_{a^2L^{2}/4\leq \Lambda} 4\arcsin\left(\frac{aL}{2\sqrt{\Lambda}}\right)\frac{\sqrt{\Lambda+M^{2}}}{a}\right) .
\ena
The first function is increasing in $\Lambda$.  The second function is decreasing in $\Lambda$ (this is also seen directly from the definition since if  $a^2L^{2}/4 \leq \Lambda$ and $x\in [-L/2,L/2]$, $ \frac{\Lambda+M^{2}}{\Lambda-a^2x^2}=1+\frac{M^{2}+a^2x^2}{\Lambda-a^2x^2}$ is decreasing as a function of $\Lambda$, so that $T(\lambda)$ is actually decreasing). As a consequence, the maximum is reached at $a^2L^{2}/4$. As a consequence, we obtain \eqref{f:T1}.

We now turn to the computation of $\Phi(\lambda)$ and recall \eqref{e:def-phi} 
\bna
\Phi(\lambda)=\int_{x_-(\lambda)}^{x_+(\lambda)}\sqrt{\lambda-V(x)}dx=\int_{x_-(\lambda)}^{x_+(\lambda)}\sqrt{\lambda-\frac{a^2x^2+M^{2}}{4}}dx,
\ena
and use 
\begin{align}
\label{formulacintrac}
\int_{-t}^t \sqrt{1-y^2} dy  = t\sqrt{1-t^2} +  \arcsin(t) , \quad t \in [-1,1] .
\end{align}
Hence, if $\Lambda=4\lambda-M^{2}\leq a^2L^{2}/4$, we obtain (as for the harmonic oscillator)
\begin{align*}
\Phi(\lambda)&=\frac{1}{2}\int_{-\frac{\sqrt{\Lambda}}{a}}^{\frac{\sqrt{\Lambda}}{a}}\sqrt{\Lambda-a^2x^2}dx=\frac{\Lambda}{2a}\int_{-1}^{1}\sqrt{1-y^2}dy=\frac{\pi}{4a}\Lambda=\frac{\pi }{a}(\lambda-\frac{M^{2}}{4}) ,
\end{align*}
while if $\Lambda=4\lambda-M^{2}\geq a^2L^{2}/4$, we get 
\begin{align*}
\Phi(\lambda)&=\frac{1}{2}\int_{-L/2}^{L/2}\sqrt{\Lambda-a^2x^2}dx=\frac{\Lambda}{2a}\int_{-\frac{aL}{2\sqrt{\Lambda}}}^{\frac{aL}{2\sqrt{\Lambda}}}\sqrt{1-y^2}dy=\frac{\Lambda}{2a}\left(\frac{aL}{4\Lambda}\sqrt{4\Lambda-a^{2}L^{2}}+\arcsin(\frac{aL}{2\sqrt{\Lambda}})\right)\\
&=\frac{L}{8}\sqrt{4\Lambda-a^{2}L^{2}}+\frac{\Lambda}{2a}\arcsin(\frac{aL}{2\sqrt{\Lambda}}) ,
\end{align*}
where we have used \eqref{formulacintrac}. The important quantity is mainly the derivative of $\Phi$: 
\begin{align*}
\Phi'(\lambda)&=\frac{\pi}{a} \quad \textnormal{ if } 0\leq 4\lambda-M^{2}\leq a^2L^2/4  ; \\
\Phi'(\lambda)&=4\frac{\d }{\d \Lambda}\Phi(\Lambda)   
=\frac{2}{a} \arcsin \left(\frac{aL}{2\sqrt{\Lambda}}\right) =\frac{2}{a} \arcsin \left(\frac{aL}{2\sqrt{4\lambda-M^2}}\right)  \quad \text{ if } 4\lambda -M^2 \geq a^2L^2/4 .
\end{align*}
Note that this is consistent with Lemma~\ref{l:phiT} stating that $\Phi'(\lambda) = \frac{1}{4\sqrt{\lambda}} T(\lambda)$.
From here, we may now compute, with $E_{0}=V(\xzero)$,
\begin{align*}
T_{E,B}&=\frac{1}{\pi}\int_{0}^{+\infty}\log \left|\frac{\Phi^{-1}(y)+E+2B}{\Phi^{-1}(y)-E}\right|dy\\
&=\frac{1}{\pi}\int_{E_0}^{+\infty}\log \left|\frac{x+E+2B}{x-E}\right|\Phi'(x)dx\\
&=\frac{1}{a}\int_{M^{2}/4}^{a^2L^2/16+M^{2}/4}\log \left|\frac{x+E+2B}{x-E}\right|dx+\frac{2}{a\pi}\int_{a^2L^2/16+M^{2}/4}^{+\infty}\log \left|\frac{x+E+2B}{x-E}\right|\arcsin\left(\frac{aL}{2\sqrt{4x-M^{2}}}\right)dx .
\end{align*}
This concludes the proof of the first part of Lemma~\ref{lm:TTEB}.

To conclude the proof of the lemma, we now analyse the different asymptotic regimes. The limits for $T_{1}$ follow from~\eqref{f:T1}. 
For $T_{E,B}$, the dominated convergence theorem (see the arguments in the proof of Lemma \ref{l:continu-I} for an effective domination) implies that the first term in the previous expression converges to zero as $a\to 0^{+}$, together with
\bna
\lim_{a\to 0^+} T_{E,B}(a) = \frac{L}{\pi}\int_{M^{2}/4}^{+\infty}\log \left|\frac{x+E+2B}{x-E}\right|\frac{1}{\sqrt{4x-M^{2}}}dx=\frac{L}{2\pi}\int_{0}^{+\infty}\log \left|\frac{y^2+M^{2}+4E+8B}{y^2+M^{2}-4E}\right|dy .
\ena
In the case $E=E_{0}=M^{2}/4$, the integral simplifies to 
$$\lim_{a\to 0^+} T_{E_0,B}(a) = \frac{L}{2\pi}\int_{0}^{+\infty}\log \left|\frac{y^2+2M^{2}+8B}{y^2}\right|dy=\frac{L}{2\pi}\sqrt{2M^{2}+8B} \int_{0}^{+\infty}\log \left|\frac{t^2+1}{t^2}\right|dt=\frac{L}{2}\sqrt{2M^{2}+8B},$$
 where we have used $\int_{0}^{+\infty}\log \left|\frac{t^2+1}{t^2}\right|dt=2\int_{0}^{+\infty}\frac{1}{1+t^2}dt=\pi$ by integration by part.
 
 We finally notice that in the limit  $a \to +\infty$, both terms in the expression of $T_{E,B}$ vanish, using $|\arcsin(s)|\leq |s|\pi/2$.
\enp

\subsection{Computation of $G_{\ref{t:estim-Cobs-expo-1D},E}=G_{\ref{thmcontrol1D},E}$ and $G_{\ref{thmlower+-},E_{0}},$}
\label{e:GG-sh}
Recall that the constants $G_{\ref{t:estim-Cobs-expo-1D},E}=G_{\ref{thmcontrol1D},E}$ and $G_{\ref{thmlower+-},E_{0}}$ are defined in~\eqref{e:def-ctes}. 
According to Theorem~\ref{thmlower+-} and Lemma~\ref{l:WE}, we only need to compute the constant $G_{\ref{thmlower+-},E}$ for$E=E_{0}$ (which corresponds to the best estimate). 
Lemma~\ref{l:WE} also implies that in the present setting, $G_{\ref{t:estim-Cobs-expo-1D},E}=G_{\ref{thmcontrol1D},E}$ is independent of $E$ by parity arguments. We are thus left to compute only $G_{\ref{t:estim-Cobs-expo-1D},E_{0}}$ and $G_{\ref{thmlower+-},E_{0}}$. 
Recall also from~\eqref{deff1D} that $\f_{M,a}^{+}$ is increasing and $\f_{M,a}^{-}=-\f_{M,a}^{+}$ is decreasing, which, according to Lemma~\ref{l:WE}, plays a key role in the computations.
\begin{lemma}\label{l:1D-explicit-uniform}
For $E=E_{0}$ and $B=0$, we have:
\begin{itemize}
\item  In case $+$:
 \begin{align*}
G_{\ref{t:estim-Cobs-expo-1D},E_{0}}& =0 , \\
G_{\ref{thmlower+-},E_{0}}& =  \frac{aL^{2}}{8}-\left(\frac{L}{8}\sqrt{a^{2}L^{2}+4M^{2}}+\frac{ M^{2}}{2a}\arcsinh(\frac{a L}{2M})\right) .
   \end{align*}
   In particular,  
     \begin{align*}
G_{\ref{thmlower+-},E_{0}}\underset{a\to 0^{+}}{\rightarrow} & -\frac{ML}{2},  \qquad G_{\ref{thmlower+-},E_{0}} \sim -\frac{M^{2}}{2a}\log(a) \underset{a\to +\infty}{\rightarrow}-\infty .
   \end{align*}
\item In case $-$:  
  \begin{align*}
 G_{\ref{t:estim-Cobs-expo-1D},E_{0}}&  
 =\frac{L}{8}\sqrt{a^{2}L^{2} +4M^{2}}+\frac{ M^{2}}{2a}\arcsinh(\frac{a L}{2M}), \\
 G_{\ref{thmlower+-},E_0}& =\frac{aL^{2}}{8} .
   \end{align*}
   In particular,  
  \begin{align*}
 G_{\ref{t:estim-Cobs-expo-1D},E_{0}}\underset{a\to 0^{+}}{\rightarrow} &\frac{ML}{2} , \qquad  G_{\ref{t:estim-Cobs-expo-1D},E_{0}} \sim \frac{aL^2}{8}\underset{a\to +\infty}{\rightarrow}+\infty ,  \\
G_{\ref{thmlower+-},E_0}\underset{a\to 0^{+}}{\rightarrow} & 0, \qquad G_{\ref{thmlower+-},E_0} \underset{a\to +\infty}{\rightarrow}+\infty
   \end{align*}
   \end{itemize}
\end{lemma}
Recall that in this situation $ T_{\ref{t:estim-Cobs-expo-1D},E_{0}}= G_{\ref{t:estim-Cobs-expo-1D},E_{0}}/E_{0}$ and $E_{0}=M^{2}/4$. Moreover, Theorem~\ref{t:estim-Cobs-expo-1D} formulates
$T_{\unif}(\{-1/2\}) \geq T_{\ref{t:estim-Cobs-expo-1D},E_{0}}= G_{\ref{t:estim-Cobs-expo-1D},E_{0}}/E_{0}$.

\bnp[Proof of Lemma~\ref{l:1D-explicit-uniform}]
Let us begin with the computation of $G_{\ref{t:estim-Cobs-expo-1D},E_{0}}$, following the simplifications of Lemma~\ref{l:WE} using that $\f_{M,a}^{\pm}$ is odd. In case $+$, this lemma yields $G_{\ref{t:estim-Cobs-expo-1D},E}=0$, and in case $-$, we have (recalling~\eqref{deff1D} and that we are working on the translated interval $[-L/2,L/2]$)
\bna
 G_{\ref{t:estim-Cobs-expo-1D},E}&=
-\f_{M,a}^{-}(L/2)=\f_{M,a}^{+}(L/2) 
 =\frac{L}{8}\sqrt{a^{2}L^{2}+4M^{2}}+\frac{ M^{2}}{2a}\arcsinh(\frac{a L}{2M}) .
\ena
We now compute $G_{\ref{thmlower+-},E_{0}}$. 
In case $+$ using that $\f_{M,a}^{+}$ is odd together with Lemmata~\ref{l:WE} and \ref{lmAgmonexpl} gives
\bna
 G_{\ref{thmlower+-},E_0}=&2d_{A, E_{0}}(L/2)-\f_{M,a}^{+}(L/2) = \frac{aL^{2}}{8}-\left(\frac{L}{8}\sqrt{a^{2}L^{2}+4M^{2}}+\frac{ M^{2}}{2a}\arcsinh(\frac{a L}{2M})\right) .
\ena
Now, in the case $-$, using again that $\f_{M,a}^{-}$ is odd together with Lemma~\ref{l:WE}, we obtain 
\bna
 G_{\ref{thmlower+-},E}=&2d_{A, E}(L/2)=\frac{aL^{2}}{8} .
\ena
The asymptotic behaviors follow from the fact that $\arcsinh(s)= s + \grando{s^3}$ near zero and $\arcsinh(s)\sim \log(s)$ near $+\infty$. 
\enp

\subsection{Asymptotics $a\to +\infty$}
\label{s:a=infty}
Recalling that $\arcsinh(t)=\log\left(t+\sqrt{1+t^2}\right)$, we have in this case the following asymptotic behaviors as $a\to +\infty$: 
\begin{itemize}
\item $T_{\f^{\pm '}_{M,a}} (\{-L/2\}) \sim_{a\to +\infty} 2\frac{\log(a)}{a} \underset{a\to +\infty}{\longrightarrow } 0^+$ according to Lemma~\ref{l:1D-explicit-limit}, i.e. the limit transport equation is controllable in small time for large $a$.
\item if we choose the sign $-$ (note that in this case, the control disappears in the limit transport equation and it is only zero on the right), then according to~\eqref{e:def-T14} and Lemma~\ref{l:1D-explicit-uniform}, we have
$$
T_{\unif} \geq T_{\ref{t:estim-Cobs-expo-1D}} \geq \frac{1}{E_0}  G_{\ref{t:estim-Cobs-expo-1D},E_{0}} \sim_{a\to + \infty}\frac{a}{2} \frac{L^2}{M^2}\to + \infty ,
$$
i.e. the minimal uniform control time tends to $+\infty$ for large $a$.

\item  if we choose the sign $+$, this is not useful since in this case $T_{\ref{t:estim-Cobs-expo-1D},E_{0}} =0$ according to Lemmata~\ref{l:1D-explicit-uniform} and~\ref{l:WE}.
\end{itemize}

\subsection{Formal limit $a \to 0^+$: comparison with the Coron-Guerrero case}
\label{s:cascst}
The computations performed and the explicit constants obtained in Sections~\ref{subsectexplcit1D}--\ref{e:GG-sh} do not apply to the situation studied in Coron-Guerrero \cite{CG:05}. The latter would correspond to the function $\f^\pm_{M,a}$ in~\eqref{deff1D} with $a=0$, and thus can be seen as a formal limit $a\to 0$ in Sections~\ref{subsectexplcit1D}--\ref{e:GG-sh}. Even if our results do not apply to the case $a=0$ and our study does not allow to make this limit rigorous, we believe it is worth computing the limit of the different bounds we obtain in this asymptotic regime.

First, we notice that the (formal) limit $a\to 0^+$  in Lemma~\ref{l:1D-explicit-limit} yields the appropriate control/flushing time for the limit equation:
$$
T_{\f^{\pm '}_{M,a}} (\{-L/2\}) \to \frac{L}{M}
$$

Second, we comment on the lower bound $T_{\unif} \geq T_{\ref{thmlower+-}}$ given by Theorem~\ref{thmlower+-}. According to~\eqref{e:def-ctes}--\eqref{e:def-T15}, this rewrites
$$
 T_{\ref{thmlower+-}}= \sup_{E\in V(\M), B\geq 0} \frac{1}{E+B}\left( G_{\ref{thmlower+-},E}+T_{E,B} \right)  \geq \sup_{B\geq 0} \frac{1}{E_0+B}\left( G_{\ref{thmlower+-},E_0}+T_{E_0,B} \right)   .
$$
According to Lemma~\ref{lm:TTEB}
$\lim_{a\to 0^{+}} T_{E_0,B}= \frac{L}{2}\sqrt{2M^{2}+8B}$ (here $E_{0}=\frac{M^2}{4}$) and according to Lemma~\ref{l:1D-explicit-uniform}, we deduce that in the limit $a\to 0^+$,
\begin{align*}
 \frac{1}{E_0+B}\left( G_{\ref{thmlower+-},E_0}+T_{E_0,B} \right)  & \to\frac{1}{M^2/4+B}\left(-\frac{ML}{2} +  \frac{L}{2}\sqrt{2M^{2}+8B} \right)  \\
  & \quad  = \frac{2L}{M^2+4B}\left( -M+\sqrt{2}\sqrt{M^2+4B}\right) ,  \quad   \text{ (Case $+$)} ,\\
  \frac{1}{E_0+B}\left( G_{\ref{thmlower+-},E_0}+T_{E_0,B} \right)  & \to\frac{1}{M^2/4+B}\left(0 + \frac{L}{2}\sqrt{2M^{2}+8B} \right)  = \frac{2\sqrt{2}L}{\sqrt{M^2+4B}} ,\quad  \text{ (Case $-$)} .
\end{align*}
Theorem~\ref{thmlower+-} gives us that 
\begin{align*}
\liminf_{a \to 0^+}T_{\unif,a} \geq \liminf_{a \to 0^+}T_{\ref{thmlower+-},a}\geq &\frac{2L}{M^{2}+4B}\left(-M+\sqrt{2}\sqrt{M^{2}+4B}\right) ,\quad \text{ (Case $+$)}  , \\
\liminf_{a \to 0^+}T_{\unif,a} \geq \liminf_{a \to 0^+}T_{\ref{thmlower+-},a}\geq&\frac{2\sqrt{2}L}{\sqrt{M^{2}+4B}},\quad \text{ (Case $-$)} .
\end{align*}

The maximum of $x\mapsto -\frac{M}{x}+\frac{\sqrt{2}}{\sqrt{x}}$ (for $x>0$) is reached for $x=2M^{2}$, so the maximum of the first expression is reached when $B=M^{2}/4$. The second case is better when $B=0$, so we get~\eqref{e:ato0+}--\eqref{e:ato0-}.

\bigskip
Let us now comment on the upper bound of Theorem \ref{thmcontrol1D} when $a \to 0^+$. The fact that $T_1 \underset{a\to 0^{+}}{\rightarrow} +\infty$ as stated in Lemma \ref{lm:TTEB} suggests that the quantity $T_1$ (which appears as the spectral gap of the $\beta_\ell$ in \eqref{gapbeta}) is not the appropriate one (at least in this regime). 
Indeed, in the case $a=0$, the operator is $ P_\eps := - \eps^2 \d_x^2 + \frac{M^2}{4} $ and the associated eigenfunctions are $\psi_k^{\e}=\sin\left(\frac{k\pi x}{L}\right)$, $k\in \N^{*}$, with the eigenvalues $\lambda_k^{\e}=\left(\frac{\e k\pi }{L}\right)^{2} + \frac{M^2}{4} $. In particular (compare with Theorem~\ref{t:Peps-spect-Allib}), one can check that the family $\e^{-1}\sqrt{\lambda_\ell^\eps}$ does not have a uniform (in $\eps$) gap.
However, in this particular setting, this issue is solved by making a translation of the spectrum, replacing $\lambda_k^{\e}$ by $\lambda_k^{\e}- \frac{M^2}{4}$. From the control point of view, this only consists in making the change of unknown $u=e^\frac{M^2 t}{4\e}v$ and new control $h_u(t)=e^\frac{M^2 t}{4\e}h(t)$ inside Proposition \ref{propcontrolzerointerm}.

The new family $\eps^{-1}\sqrt{\lambda_k^{\e}- \frac{M^2}{4}}$ then enjoys a uniform gap as in~\eqref{e:gap-rac-asympt} with $T_{1}=2L$. Hence,~\eqref{e:gap-rac-asympt} (or equivalently~\eqref{gapbeta}) is fulfilled and our proof of Theorem~\ref{thmcontrol1D} then adapts to this problem. The constants involved are however slightly less accurate than those available in the literature~\cite{CG:05,Gla:10,Lissy:12,Lissy:14,Lissy:15,Darde:17}, and we therefore do not pursue in this direction.

\appendix

\section{Some results from~\cite{Allibert:98}}
\label{s:Allibert}
In this section, we extract and translate in our context some of the results of Allibert~\cite{Allibert:98}. 
All along the section, to match the setting of~\cite{Allibert:98}, we assume that $V \in C^\infty([0,L])$, and that Items~\ref{A1}--\ref{A4} in Assumption~\ref{assumptions} are satisfied. We apply as a blackbox the results of~\cite{Allibert:98} and hence also have to assume that $\qf=0$.
We recall that the function $\Phi (\lambda)$ and $T_1$ are defined respectively in~\eqref{e:def-phi} and~\eqref{e:def-TT1}.

The goal of this appendix is to deduce a proof of the following result from the results of~\cite{Allibert:98}.
\begin{theorem}
\label{t:Peps-spect-Allib}
Consider the operator $P_\eps := - \eps^2 \d_x^2 + V(x)+ \eps  \qf + \eps^2 W$, with $V \in C^\infty([0,L];\R)$, $W \in L^\infty((0,L);\R)$, acting on the space $L^2((0,L), dx)$, with domain $H^2 \cap H^1_0(0,L)$. Assume further that $V$ satisfies Items~\ref{A1}--\ref{A4} in Assumption~\ref{assumptions} and that $\qf=0$. 
Denote by $(\lambda_k^{\e})_{k \in \N}$, the sequence of eigenvalues of the operator $P_\eps$, sorted so that $\lambda_k^{\e} \leq \lambda_{k+1}^{\e}$.
The following properties hold true:
\begin{enumerate} 
\item \label{i:simplicity} There is $C>0$ such that we have $V(\xzero)- C\eps^{2} \leq \lambda_0^{\e}$ and $\lambda_k^{\e} < \lambda_{k+1}^{\e}$ for all $k \in \N$ and $\eps \in (0,1)$.
\item  \label{i:Phi-Allibert} There exist $D,C_0, \eps_0>0$ such that for all $\eps\in (0,\eps_0)$ and $k \in \N$, we have 
\bnan
\label{stimlambda}
\Phi^{-1}\big(\e (\pi k-D) \big)-C_0\eps^{3/2}\leq \lambda_k^{\e}\leq \Phi^{-1}\big(\e (\pi k+D) \big)+C_0\eps^{3/2}.
\enan
\item \label{itemallibertgap} 
For all $\delta>0$, there are $\eps_0, N >0$ such that for all $\eps\leq \eps_0$, 
\bnan
\label{e:gap-rac-asympt}
\sqrt{\lambda_{\ell+1}^\eps} -\sqrt{\lambda_\ell^\eps} \geq \eps \frac{2\pi}{T_1+\delta}, \quad \text{ for all } \ell \geq N .
\enan
Moreover, there are $\eps_0, \gamma, \gamma_{2} >0$ such that for all $\eps\leq \eps_0$, 
\begin{align}
\sqrt{\lambda_{\ell+1}^\eps} -\sqrt{\lambda_\ell^\eps}  \geq  \eps \gamma>0 ,\quad \text{for all }  \ell\in\N \label{e:gap-rac}\\
\lambda_{\ell+1}^\eps -\lambda_\ell^\eps\geq \eps \gamma_{2} >0 ,\quad \text{for all }  \ell\in\N .  \label{e:gap-lambda}
\end{align}
\end{enumerate}
\end{theorem}
To prove this result, let us now recall the setting of~\cite{Allibert:98}.
Allibert~\cite{Allibert:98} considers on $(0,b)$ for $b>0$, the operator
$$
P_h^{All} u = - \frac{h^2}{R(z)\sqrt{1+R'(z)^2}} \d_z \left(\frac{R(z)}{\sqrt{1+R'(z)^2}} \d_z u \right) + \frac{1}{R^2(z)} u .
$$
This operator is selfadjoint on the space $L^2\big((0,b) ,R(z)\sqrt{1+R'(z)^2} dz \big)$ with domain $H^2\cap H^1_0(0,b)$, with the assumption that $z\mapsto \frac{1}{R^2(z)}$ admits at $c \in (0,b)$ a strict nondegenerate minimum, and that $\frac{1}{R^2(c)} < \frac{1}{R^2(b)} < \frac{1}{R^2(0)}$.

The geometric quantities entering into the discussion are
\begin{align*}
\Phi^{All} (\lambda) & = \int_{z_-(\lambda)}^{z_+(\lambda)}\sqrt{1+R'(z)^2}  \sqrt{\lambda - \frac{1}{R^2(z)}} dz ,\\
T_1^{All} & =  \sup_{\lambda \geq \frac{1}{R^2(c)}} T^{All}(\lambda) , \quad T^{All}(\lambda)  = 2 \int_{z_-(\lambda)}^{z_+(\lambda)} \sqrt{1+R'(z)^2} \frac{\sqrt{\lambda}}{\sqrt{\lambda - \frac{1}{R^2(z)}}} dz .
\end{align*}
In this expression, $z_-(\lambda)$ is the solution to $\frac{1}{R^2(z_-(\lambda))} = \lambda$ with $z_-(\lambda)\leq c$ for $\lambda \leq \frac{1}{R^2(0)}$, and $z_-(\lambda) = 0$ for $\lambda \geq \frac{1}{R^2(0)}$.
Similarly, $z_+(\lambda)$ is the solution to $\frac{1}{R^2(z_+(\lambda))} = \lambda$ with $z_+(\lambda)\geq c$ for $\lambda \leq \frac{1}{R^2(b)}$, and $z_+(\lambda) = b$ for $\lambda \geq \frac{1}{R^2(b)}$.

\medskip
We want to compare this setting to the one considered here. Namely, we study the operator $P_\eps := - \eps^2 \d_x^2 + V(x)+ \eps  \qf $, with $V(x) = \frac{|\f'(x)|^2}{4}$ and $\qf= \frac{\f''}{2}-\q$, acting on the space $L^2((0,L), dx)$, with domain $H^2\times H^1_0(0,L)$.

\medskip
We define the increasing diffeomorphism \begin{equation*}
\begin{array}{rcl}
x:  [0,b] & \to & [0,L] \\
z&\mapsto & x(z) = \int_0^z\sqrt{1+R'(t)^2} dt
\end{array}
\end{equation*}
where we set $L := \int_0^b \sqrt{1+R'(t)^2} dt$. Given a potential $V \in C^k([0,L])$ satisfying Items~\ref{A1}--\ref{A4} in Assumption~\ref{assumptions}, 
the function  $\frac{1}{R^2(z)}:=V(x(z))$ satisfies the assumptions of Allibert~\cite{Allibert:98} with $c$ given by $x(c) = \xzero$.

We obtain, with this change of variables, that $\Phi^{All} (\lambda) = \Phi(\lambda),  T^{All}(\lambda)=T(\lambda)$ and $T_1^{All}=T_1$, where $\Phi(\lambda), T(\lambda)$ and $T_1$ are defined in~\eqref{e:def-phi}-\eqref{e:def-TT1}. Note that $x_-(\lambda)= x(z_-(\lambda))$ and $x_+(\lambda)= x(z_+(\lambda))$ in these definitions.

\medskip
Moreover, under this change of variable, we have $\frac{\d}{\d s} = \frac{1}{\sqrt{1+R'(z)^2}}\frac{\d}{\d z}$ so that the operator $P_h^{All}$ becomes 
$$
 - \frac{h^2}{\mathsf{R}(s)} \d_s \left(\mathsf{R}(s)\d_s  \cdot \right) +  V(s) ,
$$
where we have written $\mathsf{R}(x(z)) = R(z)$ (and hence $V(s) = \frac{1}{\mathsf{R}(s)^2}$). This operator acts on the space $L^2\big((0,L) ,\mathsf{R}(s)ds \big)$ with domain $H^2\cap H^1_0(0,L)$. 

Now, observe that the map
\begin{equation*}
\begin{array}{rcl}
T: L^2((0,L),\mathsf{R}(s)ds) & \to & L^2((0,L), ds) \\
u  & \mapsto & Tu, \quad \text{with} \quad (Tu)(s) = \mathsf{R}(s)^{\frac12} u(s)
\end{array}
\end{equation*}
is an isometry and the conjugated operator of $\frac{1}{\mathsf{R}(s)} \d_s \left(\mathsf{R}(s)\d_s  \cdot \right)$ is
\begin{align*}
T\left( \frac{1}{\mathsf{R}(s)} \d_s \left(\mathsf{R}(s)\d_s  \cdot \right) \right) T^{-1}u = \mathsf{R}^{1/2}\left( \frac{1}{\mathsf{R}(s)} \d_s \left(\mathsf{R}(s)\d_s (\mathsf{R}^{-1/2}u) \right) \right)= \partial_s^2 u- V_1 u  ,
\end{align*}
where $V_1(s)=-\frac{1}{4}\frac{\mathsf{R}'(s)^2}{\mathsf{R}(s)^2}+\frac{1}{2}\frac{\mathsf{R}''(s)}{\mathsf{R}(s)^2}$. 

We thus obtain 
\begin{equation}
\label{e:Pall}
\mathsf{P}_h := T P_h^{All} T^{-1}  = - h^2 \partial_s^2 + V(s) + h^2 V_1(s) .
\end{equation}
This is almost the operator we consider, except for lower order terms.

Then, Allibert~\cite{Allibert:98} describes the spectrum of the operator $P_h^{All}$:
\begin{enumerate}
\item he constructs in \cite[Lemmata~6-7 and Section~3.1.2]{Allibert:98} approximate eigenvalues and eigenfunctions. The approximate eigenvalues are $O(h^{3/2})$ close to real eigenvalues;
\item he proves in \cite[Section~3.1.3]{Allibert:98} that the sequence of real eigenvalues constructed in the first point actually contains {\em all} eigenvalues (using a Sturm-Liouville argument); in particular, the spectrum is simple;
\item he computes in \cite[Section~3.1.4]{Allibert:98} the spectral gap (using the explicit expression of the approximate eigenvalues).
\end{enumerate}
We first collect the following properties of $\mathsf{P}_\eps$ from~\cite{Allibert:98}.

\begin{theorem}
\label{lmstimlambda}
Consider the operator $\mathsf{P}_\eps$ acting on the space $L^2((0,L), dx)$, with domain $H^2\cap H^1_0(0,L)$. 
Denote by $(\lambda_k^{\e})_{k \in \N}$, the sequence of eigenvalues of the operator $\mathsf{P}_\eps$, sorted so that $\lambda_k^{\e} \leq \lambda_{k+1}^{\e}$.
Assuming that $V \in C^\infty([0,L])$ satisfies Items~\ref{A1}--\ref{A4} in Assumption~\ref{assumptions}, then Items~\ref{i:simplicity},~\ref{i:Phi-Allibert},~\ref{itemallibertgap} of Theorem~\ref{t:Peps-spect-Allib} hold for the eigenvalues $(\lambda_k^{\e})_{k \in \N}$ of $\mathsf{P}_\eps$.
\end{theorem}
Note that a much finer property than~\eqref{stimlambda} is actually proved in~\cite{Allibert:98}, but this weaker form is sufficient for our needs.
\bnp[Proof of Theorem~\ref{lmstimlambda} from~\cite{Allibert:98}]
The first lower bound in Item~\ref{i:simplicity} comes from $$(\mathsf{P}_\eps u,u)_{L^2} \geq \left(V(\xzero)- \|V_1\|_{L^\infty} \eps^{2}\right) \|u\|_{L^2}^2,$$ and the simplicity of the spectrum from the fact that we consider Dirichlet boundary conditions (hence, the space of solutions to the ODE eigenvalue equation has dimension one).

Let us now explain how Item~\ref{i:Phi-Allibert} is deduced from~\cite[Lemme~6 and Lemme~7]{Allibert:98}. Firstly, note that these properties concern the eigenvalues of the operator $P_\eps^{All}$, which, according to \eqref{e:Pall}, are exactly those of $\mathsf{P}_\eps$. 
Secondly~\cite[Lemme~6 and Lemme~7]{Allibert:98} prove the existence of a sequence $\mu_{k,\eps}$ such that 
\begin{equation}
\label{e:Phi-mu-k}
\Phi(\mu_{k,\eps}) = \eps k \pi + \eps \Theta (\eps, \mu_{k,\eps}) ,
\end{equation}
where $\Theta : [0,1]\times \R^+ \to \R$ is a uniformly bounded function, and an eigenvalue of $\mathsf{P}_\eps$, $\tilde{\lambda}_k^{\e} \in \{\lambda_\ell^{\e}, \ell \in \N\}$, such that $|\tilde{\lambda}_k^{\e} - \mu_{k,\eps}| \leq C_0 \eps^{3/2}$.
Then, in~\cite[Section~3.1.3]{Allibert:98}, he proves that the set $\{\tilde\lambda_\ell^{\e}, \ell \in \N\}$ constructed that way coincides with the spectrum, that is $\tilde{\lambda}_k^{\e} =\lambda_k^{\e}$ for all $k \in \N$. This implies 
\begin{equation}
\label{e:mu-k-lambda-k}
\mu_{k,\eps} - C_0 \eps^{3/2} \leq \lambda_k^{\e} \leq \mu_{k,\eps} + C_0 \eps^{3/2}
\end{equation}
We set $D := \sup_{(\eps,\mu) \in [0,1]\times \R^+} |\Theta(\eps, \mu)|$. As a consequence of~\eqref{e:Phi-mu-k},~\eqref{e:mu-k-lambda-k}, together with the fact that $\Phi$ is increasing, we obtain
\begin{align*}
\Phi( \lambda_k^{\e} - C_0 \eps^{3/2}) & \leq \Phi( \mu_{k,\eps} ) = \eps k \pi + \eps \Theta (\eps, \mu_{k,\eps}) \leq \eps (k \pi + D)\\
\Phi( \lambda_k^{\e} + C_0 \eps^{3/2}) & \geq \Phi( \mu_{k,\eps} ) = \eps k \pi + \eps \Theta (\eps, \mu_{k,\eps}) \geq \eps (k \pi - D) ,
\end{align*}
which proves~\eqref{stimlambda}.

Finally, concerning Item \ref{itemallibertgap}, \eqref{e:gap-rac-asympt} and~\eqref{e:gap-rac} are proved in~\cite[Proposition~2 p 1511 and Section~3.1.4]{Allibert:98}.
The last estimate~\eqref{e:gap-lambda} comes from 
$$
\lambda_{\ell+1}^\eps -\lambda_\ell^\eps =\left(\sqrt{\lambda_{\ell+1}^\eps} -\sqrt{\lambda_\ell^\eps} \right)\left(\sqrt{\lambda_{\ell+1}^\eps} +\sqrt{\lambda_\ell^\eps} \right) 
\geq \eps \gamma 2 \sqrt{\lambda_0^\eps} \geq   2 \eps \gamma\left(V(\xzero)- C\eps^{2} \right) , 
$$
where we have used~\eqref{e:gap-rac} together with Item~\ref{i:simplicity}.
\enp

We now deduce from this result for $\mathsf{P}_\eps$ a proof of Theorem~\ref{t:Peps-spect-Allib}, that is, prove that the same properties hold for $P_\eps$.

The proof of Theorem~\ref{t:Peps-spect-Allib} from Theorem~\ref{lmstimlambda} consists in a classical perturbative (deformation) argument, and relies of the following lemma.
\begin{lemma}
\label{l:rank-proj-cont}
Let $H$ be a Banach space, and $(P(t))_{t\in [0,1]} \in \L(H)^{[0,1]}$ be a family of projectors (in the sense that $P(t)^2=P(t)$) having finite rank $r(t) \in \N$, and such that the map $t\mapsto P(t)$ is continuous $[0,1]\to \L(H)$. Then, all projectors have the same rank, i.e. $r(t)=r(0)$ for all $t\in [0,1]$.
\end{lemma}

\bnp[Proof of Theorem~\ref{t:Peps-spect-Allib} from Theorem~\ref{lmstimlambda}]
We write $\tilde{W} = W-V_1$ and set 
\begin{align*}
A_\eps(t) : =(1-t)\mathsf{P}_\eps+ tP_\eps = \mathsf{P}_\eps+ t \eps^2 \tilde{W} .
\end{align*}
We denote by $(\lambda_k^\eps)_{k\in \N}$ the spectrum of $\mathsf{P}_\eps$, which satisfies Items~\ref{i:simplicity},~\ref{i:Phi-Allibert},~\ref{itemallibertgap} of Theorem~\ref{t:Peps-spect-Allib}.
We have, for $z \notin \Sp(\mathsf{P}_\eps)$ 
\begin{align}
\label{e:link-res}
z - A_\eps(t) = (z-\mathsf{P}_\eps)\big(\id - (z-\mathsf{P}_\eps)^{-1}t\eps^2 \tilde{W}\big) .
\end{align}
Next, we remark that for all $t,\eps,z$ such that $|t|\eps^2 \nor{(z-\mathsf{P}_\eps)^{-1}}{\L}\nor{\tilde{W}}{\infty} < 1$, the operator $\big(\id - (z-\mathsf{P}_\eps)^{-1}t\eps^2 \tilde{W}\big)$ is invertible with 
\begin{align*}
\big(\id - (z-\mathsf{P}_\eps)^{-1}t\eps^2 \tilde{W}\big)^{-1} & = \sum_{n=0}^\infty \left(t\eps^2 (z-\mathsf{P}_\eps)^{-1}\tilde{W}\right)^n \\
\nor{\big(\id - (z-\mathsf{P}_\eps)^{-1}t\eps^2 \tilde{W}\big)^{-1}}{\L} & \leq \sum_{n=0}^\infty \left( |t|\eps^2 \nor{(z-\mathsf{P}_\eps)^{-1}}{\L}\nor{\tilde{W}}{\infty}\right)^n\\
&  \leq \frac{1}{1- |t|\eps^2 \nor{(z-\mathsf{P}_\eps)^{-1}}{\L}\nor{\tilde{W}}{\infty}}.
\end{align*}
Recalling that $\nor{(z-A_\eps(t))^{-1}}{\L} = \frac{1}{\dist(z, \Sp(A_\eps(t)))}$ for $z \notin \Sp(A_\eps(t))$ (since $A_\eps(t)$ is selfadjoint) together with~\eqref{e:link-res}, we deduce
\begin{align}
\label{e:info-res-A-eps}
 |t|\eps^2 \nor{(z-\mathsf{P}_\eps)^{-1}}{\L}\nor{\tilde{W}}{\infty} <1 \implies 
\left\{
\begin{array}{c}
z \notin \Sp(A_\eps(t)) , \\ 
(z- A_\eps(t))^{-1} =\big(\id - (z-\mathsf{P}_\eps)^{-1}t\eps^2 \tilde{W}\big)^{-1}  (z-\mathsf{P}_\eps)^{-1} \\
\nor{(z- A_\eps(t))^{-1}}{\L} \leq \frac{\nor{(z-\mathsf{P}_\eps)^{-1}}{\L} }{1- |t|\eps^2 \nor{(z-\mathsf{P}_\eps)^{-1}}{\L}\nor{\tilde{W}}{\infty}},
\end{array}
\right.
\end{align}
and hence
\begin{align*}
|t|\eps^2\nor{\tilde{W}}{\infty}< \dist(z, \Sp(\mathsf{P}_\eps)) \implies  
\dist(z, \Sp(\mathsf{P}_\eps)) \leq \dist(z, \Sp(A_\eps(t))) + |t|\eps^2\nor{\tilde{W}}{\infty}.
\end{align*}
Now taking $z= \nu_\eps(t) \in \Sp(A_\eps(t))$ implies that
\begin{align}
\label{e:estim-dist-spect-t}
\dist(\nu_\eps(t), \Sp(\mathsf{P}_\eps))\leq |t|\eps^2\nor{\tilde{W}}{\infty} .
\end{align}
We now recall the gap property~\eqref{e:gap-lambda} of the spectrum $(\lambda_k^\eps)_{k\in \N}$ of the operator $\mathsf{P}_\eps$, and define the contour (oriented counterclockwise)
$$
\Gamma_k^\eps = \d B(\lambda_k^\eps , \frac{\gamma_2 \eps}{3}) . 
$$
According to~\eqref{e:gap-lambda}, these sets are disjoint and each contains exactly one eigenvalue of $\mathsf{P}_\eps$. We define the associated orthogonal projector onto $\ker(\mathsf{P}_\eps - \lambda_k^\eps)$ by
$$
\Pi_k^\eps =  \int_{\Gamma_k^\eps} (z-\mathsf{P}_\eps)^{-1} dz
$$
As a consequence of~\eqref{e:estim-dist-spect-t}, we obtain that for all $t\in [0,1]$ and all $\eps\in(0,\eps_0)$ with $\eps_0$ such that $\eps_0^2\nor{\tilde{W}}{\infty}< \frac{\gamma_2\eps_0}{3}$, we have
$$
\Gamma_k^\eps \cap \Sp(A_\eps(t)) = \emptyset , \quad \Sp(A_\eps(t)) \subset \bigcup_{k\in \N} B(\lambda_k^\eps , \frac{\gamma_2 \eps}{3}) .
$$
In particular, we can define the orthogonal projector onto the spectral subspace of $A_\eps(t)$ associated to its eigenvalues inside $\Gamma_k^\eps$, namely,
$$
\Pi_k^\eps(t) :=  \int_{\Gamma_k^\eps} (z-A_\eps(t))^{-1} dz
$$
(hence $\Pi_k^\eps(0)=\Pi_k^\eps$). According to~\eqref{e:info-res-A-eps}, we have the uniform bound for $z\in \Cup_{k\in \N}\Gamma_k^\eps$,
$$\nor{(z- A_\eps(t))^{-1}}{\L} \leq \frac{\nor{(z-\mathsf{P}_\eps)^{-1}}{\L} }{1-\eps^2 \nor{(z-\mathsf{P}_\eps)^{-1}}{\L}\nor{\tilde{W}}{\infty}}$$ for $t\in [0,1]$ so that $t\mapsto\Pi_k^\eps(t)$ is continuous $[0,1] \to \L(L^2)$.
According to Lemma~\ref{l:rank-proj-cont}, we obtain $\Rk (\Pi_k^\eps(t))= \Rk (\Pi_k^\eps) = 1$ for all $t\in [0,1]$. 
As a consequence, for all $t\in [0,1]$, $A_\eps(t)$ has a single eigenvalue $\lambda_k^\eps(t)$ inside $\Gamma_k^\eps$, having multiplicity one. Moreover, from~\eqref{e:estim-dist-spect-t} we have $|\lambda_k^\eps(t)-\lambda_k^\eps| \leq |t|\eps^2\nor{\tilde{W}}{\infty}$. 

Items~\ref{i:simplicity},~\ref{i:Phi-Allibert},~\ref{itemallibertgap} of Theorem~\ref{lmstimlambda} being satisfied by $\lambda_k^\eps$, they are thus satisfied as well by $\lambda_k^\eps(t)$ for all $t\in [0,1]$ and $\eps \in (0,\eps_0)$ for $\eps_0$ sufficiently small. Note that we use that $s\mapsto \sqrt{s}$ is uniformly Lipschitz on $[V(\xzero)-C\e_0^2,+\infty)$ thanks to Item~\ref{A1} in Assumption~\ref{assumptions} and an appropriate choice of $\e_0$. This yields the sought result in case $t=1$.
\enp

We finally prove Lemma~\ref{l:rank-proj-cont}, which is a consequence of the following remark.
\begin{lemma}
\label{l:proj-cte-rank}
Let $P_1$ and $P_2$ two continuous projections with finite respective rank $r_1>r_2$ in a Banach $H$. Then, $\nor{P_1-P_2}{H \rightarrow H}\geq 1$.
\end{lemma}
\bnp[Proof of Lemma~\ref{l:proj-cte-rank}]
We define $H_1$ (resp. $H_2$) the range of $P_1$ (resp. $P_2$) which are spaces of finite dimension. We define the application $F:H_1 \rightarrow H_2$, defined by $F(x_1)=P_2(x_1)$. By the rank-nullity theorem and the assumption $r_1>r_2$, we have $\dim \ker (F)>0$ and there exists $x_1\in H_1$ with $\nor{x_1}{H}=1$ so that $P_2(x_1)=0$. But since $x_1\in H_1$, $\nor{P_1(x_1)}{H}=\nor{x_1}{H}=1$. This gives the result.
\enp
\bnp[Proof of Lemma~\ref{l:rank-proj-cont}]
given $t \in [0,1]$, there exist $\delta>0$ such that for all $t'\in [0,1]$, $|t'-t| \leq \delta \implies \nor{P(t')-P(t)}{H\to H} \leq 1/2$. This implies that $r(t') = r(t)$ for all $t' \in [t-\delta,t+\delta]\cap [0,1]$ (this would otherwise contradict Lemma~\ref{l:proj-cte-rank} since we assume that all projectors have finite rank). A connectedness argument concludes that $r$ is globally constant on $[0,1]$. 
\enp

\section{A moment result}
\label{s:moment}
The purpose of this Section is the proof of Proposition \ref{propgapheat} below which may not be new, but for which we did not find any reference, especially for the uniform dependence of the constants. The study of of biorthogonal sequences and their application to controllability of parabolic equations is classical and dates back to Fattorini-Russell \cite{FR:71,FR:74}. We also refer to Hansen \cite{H:91}, and Ammar Khodja--Benabdallah--Gonz\'alez Burgos--de Teresa \cite{ABGT:11b}. 
At the time of writing this article, Cannarsa, Martinez and Vancostenoble~\cite{CMV:20} obtained results close to the one we obtain in this section. We have chosen to keep this section since our method seems simpler, with a slightly more explicit constant. Our proof relies on an Ingham inequality, see e.g.~\cite{H:89-2,KL:05}, together with a transmutation argument
due to Ervedoza-Zuazua~\cite{EZ:11s}.

The main result of this section is the following proposition.
\begin{proposition}
\label{propgapheat}
For any $\gamma_{\infty}>0$, $\gamma>0$, $N \in \N$ and any $S>\frac{\pi}{\gamma_{\infty}}$ and $\eps>0$, we can find a constant $C=C(\gamma_{\infty},\gamma,N,S, \eps)>0$ so that for any sequence $(\beta_{n})_{n\in \N^*}$ satisfying 
\begin{enumerate}
\item \label{e:item1-biorth} $\beta_{n+1}-\beta_{n}\geq \gamma$ for all $n\in \N^*$ and $\beta_1 \geq \gamma$,
\item \label{e:item2-biorth}$\beta_{n+1}-\beta_{n}\geq \gamma_{\infty}$ for any $n\in \N^*$ with $n \geq N$,
\end{enumerate}
 and for any $0<T\leq 1$, there exists a sequence of functions $(u_{n})_{n\in \N^*}\in L^{2}(0,T)^{\N^*}$ so that 
\begin{enumerate}
\item \label{i:bibiorth} for any $l,n\in \N^*$, we have $\int_{0}^{T} u_{n}(t)e^{-\beta_{l}^{2}t}dt=\delta_{n,l}$,
\item For any $(a_n)_{n\in \N^*}$ so that $\beta_n^2 e^{\beta_{n}^{2}T}a_n\in\ell^2(\N^*)$, we have
\bna
\nor{z}{L^2(0,T)}^2\leq \frac{C}{T^3}e^{\frac{(16+\eps)S^{2}}{T}}\sum_{n\in \N^*}\beta_n^2e^{2\beta_{n}^{2}T} |a_n|^2, \quad \text{ with } z(t)=\sum_{n\in \N^*}a_n u_n(t) .
\ena
\end{enumerate}
\end{proposition}

Our proof relies on the following classical inequality due to Ingham-Haraux~\cite{H:89-2,KL:05}.
\begin{theorem}[Ingham-Haraux]
\label{thmIngham}
For any $\gamma_{\infty}>0$, $\gamma>0$, and $N \in \N$ and any $S>\frac{2\pi}{\gamma_{\infty}}$, we can find a constant $C_0=C_0(\gamma_{\infty},\gamma,N,S)>0$ so that for any sequence $(\mu_{k})_{k\in \Z}$ satisfying:
\begin{enumerate}
\item $\mu_{k+1}-\mu_{k}\geq \gamma$ for all $k\in \Z$,
\item $\mu_{k+1}-\mu_{k}\geq \gamma_{\infty}$ for any $k \in \Z$ with $|k| \geq N$,
\end{enumerate}
then, we have
\bna
C_0^{-1}\int_{0}^{S} \left| \sum_{k\in \Z}a_{k}e^{i\mu_{k}s}\right|^{2}ds \leq \sum_{k\in \Z}|a_{k}|^{2}\leq C_0 \int_{0}^{S} \left| \sum_{k\in \Z}a_{k}e^{i\mu_{k}s}\right|^{2}ds
\ena
for all $(a_{k})_{k\in \Z}\in \ell^{2}(\Z)$ with finite support.
\end{theorem}
Note that in these estimates, only the length of the time interval $(0,S)$ is relevant; under the assumption that $S>\frac{\pi}{\gamma_{\infty}}$ the conclusion holds with the integrals over $(0,S)$ replaced by integrals over $(-S,S)$.
\begin{corollary}
\label{corbiorthIngham}
For any $\gamma_{\infty}>0$, $\gamma>0$, $N \in \N$ and any $S>\frac{\pi}{\gamma_{\infty}}$, we can find a constant $C_0=C_0(\gamma_{\infty},\gamma,N,S)>0$ so that for any sequence $(\beta_{n})_{n\in \N^*}$ satisfying Item \ref{e:item1-biorth}-\ref{e:item2-biorth} of Proposition \ref{propgapheat},
there exists a sequence of functions $(v_{n})_{n\in \N}\in L^{2}((-S,S))^{\N^*}$ so that 
\begin{align}
\int_{-S}^{S} v_{n}(s)\sin(\beta_{l}s)~ds=\delta_{n,l}, \quad \text{for all }l,n\in \N^*,\quad \text{and} \nonumber\\
\label{orthodual}
C_0^{-1}\int_{-S}^{S} \left| \sum_{n\in \N^*}b_{n}v_n(s)~ds\right|^{2}\leq \sum_{n\in \N^*}|b_{n}|^{2}\leq C_0 \int_{-S}^{S} \left| \sum_{n\in \N^*}b_{n}v_n(s)~ds\right|^{2}
\end{align}
for all $(b_{n})_{n\in \N^*}\in \ell^{2}(\N^*)$.
\end{corollary}
\bnp
For $k\in \Z$, we set $\mu_{k}:=\beta_{k}$ if $k> 0$, $\mu_{k}:=-\beta_{-k}$ if $k< 0$ and $\mu_0=0$. That the sequence $(\mu_k)_{k\in \Z}$ satisfies the assumptions of Theorem \ref{thmIngham} readily follows from the assumptions.

Given $(b_{n})_{n\in \N^*}\in \ell^{2}(\N^*)$, we define for $k\in \Z$, $a_{k}:=\frac{b_{k}}{2i}$ if $k> 0$, $a_{k}:=-\frac{b_{-k}}{2i}$ if $k<0$ and $a_0:=0$. We have
\bna
\sum_{k\in \Z}a_{k}e^{i\mu_{k}s}= \sum_{k\in \N^*} \frac{b_{k}}{2i}e^{i\beta_k s} -\frac{b_{-k}}{2i}e^{-i\beta_k s}= \sum_{n\in \N^*}b_{n}\sin(\beta_{n}s).
\ena
Theorem \ref{thmIngham} (applied on the time interval $(-S,S)$, of length $>\frac{2\pi}{\gamma_{\infty}}$) gives 
\bna
C_0^{-1}\int_{-S}^{S} \left| \sum_{n\in \N^*}b_{n}\sin(\beta_{n}s)\right|^{2}ds \leq \sum_{n\in \N^*}|b_{n}|^{2}\leq C_0\int_{-S}^{S} \left| \sum_{n\in \N^*}b_{n}\sin(\beta_{n}s)\right|^{2}ds.\ena
In particular, the family $\left(\sin(\beta_{n}s)\right)_{n\in \N}$ forms a Riesz basis of the space it spans in $L^{2}(-S,S)$. Lemma~\ref{lmbiorth} below in $H=L^2(-S,S)$ yields the existence of a biorthogonal family $(v_n)_{n\in \N^*}$ to $(\sin(\beta_n s))_{n\in \N^*}$ satisfying~\eqref{orthodual} for the same constant $C_0$.
\enp

To deduce a proof of Proposition~\ref{propgapheat}, we now construct from the  sequence biorthogonal to $(\sin(\beta_n s))_{n\in \N^*}$ in $L^{2}(-S,S)$, a sequence biorthogonal to $(e^{-\beta_n^2 t})_{n\in \N^*}$ in $L^2(0,T)$ satisfying precise bounds. To this aim, we use ideas coming from transposition from heat to waves, see \cite{Miller:06b,EZ:11,EZ:11s}, and more precisely a kernel constructed in~\cite{EZ:11s}.
\bnp[Proof of Proposition~\ref{propgapheat}]
According to~\cite[Section~3.1]{EZ:11s}, given $\alpha > 2 S^2$, there exists a kernel function $k_T(t,s) \in C^\infty(\R^2)$ solution to 
\bneqn
\label{heatkernel}
\partial_t k_T(t,s) + \partial_s^2k_T(t,s) &= &0, \quad \text{ for } s \in (-S,S) ,~ t\in (0,T),\\
\left(k_T(t,s),\partial_{s}k_T(t,s)\right)|_{s=0}&=& \left(0,e^{-\alpha \left(\frac{1}{t}+\frac{1}{T-t}\right)}\right),\\
k_T(t,s)|_{t=0}=k_T(t,s)|_{t=T}&=&0, \quad \supp(k_T) \subset [0,T] \times \R.
\eneqn
and such that~\cite[Proposition~3.1]{EZ:11s} for all $\delta \in (0,1)$ and all $(t, s) \in (0,T) \times (-S,S)$, $k_T$ satisfies
\bnan
\label{e:estim-kT}
|k_T(t,s)|\leq |s|\exp\left( \frac{1}{\min \left\{ t,T-t\right\}}\left(\frac{s^2}{\delta}-\frac{\alpha}{(1+\delta)}\right)\right) .
\enan
Let $(v_n)_{n\in \N^*}$ the sequence given by Corollary~\ref{corbiorthIngham}. We define $w_n\in C^\infty(\R)$ by
\bna
w_{n}(t):= \int_{-S}^{S}k_{T}(t,s)v_{n}(s)ds, \quad \supp w_n \subset [0,T]. 
\ena
We compute
\bnan
\label{e:fl}
\int_{0}^{T} w_{n}(t)e^{-\beta_{l}^{2}t}dt=\int_{0}^{T} \int_{-S}^{S}k_{T}(t,s)v_{n}(s) e^{-\beta_{l}^{2}t}~dtds=\int_{-S}^{S} v_{n}(s)f_{l}(s)ds
\enan
where we have set $f_{l}(s)= \int_{0}^{T} k_{T}(t,s) e^{-\beta_{l}^{2}t}~dt$. Using~\eqref{heatkernel}, we have for $s \in (-S,S)$
\bna
\frac{d^{2}}{ds^{2}}f_{l}(s)=\int_{0}^{T} \frac{\d^{2}}{\d s^{2}}k_{T}(t,s) e^{-\beta_{l}^{2}t}~dt=-\int_{0}^{T} \frac{\d}{\d t}[k_{T}(t,s)] e^{-\beta_{l}^{2}t}~dt=- \beta_{l}^{2}\int_{0}^{T} k_{T}(t,s) e^{-\beta_{l}^{2}t}dt=- \beta_{l}^{2}f_{l}(s),
\ena
where we have performed an integration by parts using the zero boundary conditions of $k_{T}$ at $t=0$ and $t=T$. Noticing that $f_{n}(0)=0$ and $f_n'(0)= \int_0^Te^{-\alpha \left(\frac{1}{t}+\frac{1}{T-t}\right)}e^{-\beta_{l}^{2}t}dt$, using~\eqref{heatkernel}, we obtain
\bnan
\label{e:def-cl}
f_{l}(s)=c_{l} \sin(\beta_{l}s) , \quad \text{with} \quad c_{l} =\frac{1}{\beta_l}\int_{0}^{T}e^{-\alpha \left(\frac{1}{t}+\frac{1}{T-t}\right)}e^{-\beta_{l}^{2}t}~dt. 
\enan
In particular, using the definition of $v_{n}$ in Corollary~\ref{corbiorthIngham} together with~\eqref{e:fl}, we have 
\bna
\int_{0}^{T} w_{n}(t)e^{-\beta_{l}^{2}t}dt=c_{l }\int_{-S}^{S} v_{n}(s)\sin(\beta_{l}s)ds=c_{l}\delta_{n,l} , \quad n,l \in \N^*.
\ena
Therefore, defining $u_{n}(t):=c_{n}^{-1}w_{n}(t)$ for any $n\in \N^*$, the sequence $(u_{n})_{n\in \N^*}$ forms a family biorthogonal to $(e^{-\beta_{l}^{2}t})_{l\in\N^*}$ in $L^2(0,T)$, which proves Item~\ref{i:bibiorth} of the proposition. It only remains to estimate $u_{n}$ to conclude the proof, that is, estimate $c_{l}^{-1}$. We have, performing the change of variable $\sigma=\frac{2t}{T}-1$, for all $\nu \in (0,1)$,
$$
\int_{0}^{T}e^{-\alpha \left(\frac{1}{t}+\frac{1}{T-t}\right)}dt  =\frac{T}{2} \int_{-1}^1 e^{-\frac{4\alpha}{T}\left(\frac{1}{1-\sigma^2}\right)}d\sigma
\geq \frac{T}{2}  \int_{-\nu}^\nu e^{-\frac{4\alpha}{T}\left(\frac{1}{1-\sigma^2}\right)}d\sigma
\geq T\nu e^{-\frac{4\alpha}{T}\left(\frac{1}{1-\nu^2}\right)}
$$
and thus, with $\nu = \left(1+ \frac{\alpha}{T} \right)^{-1/2}$, we have $\frac{1}{1-\nu^2}=1+\frac{T}{\alpha}$ and this lower bound reads
$$
\int_{0}^{T}e^{-\alpha \left(\frac{1}{t}+\frac{1}{T-t}\right)}dt \geq T\left(1+ \frac{\alpha}{T} \right)^{-1/2} e^{-\frac{4\alpha}{T}-4} \geq C T^{3/2} e^{-\frac{4\alpha}{T}}, 
$$
for $T \in [0,1]$. As a consequence, with $c_l$ defined in~\eqref{e:def-cl}, we have the rough esimate
\bnan
\label{estimcl}
c_{l} \geq \frac{1}{\beta_l}e^{-\beta_{l}^{2}T}\int_{0}^{T}e^{-\alpha \left(\frac{1}{t}+\frac{1}{T-t}\right)}~dt\geq C\frac{ T^{3/2}}{\beta_l} e^{-\beta_{l}^{2}T}e^{-\frac{4\alpha }{T}}.
\enan
Finally, for a finite sequence $(a_n)_{n\in \N^*}$ (see below for a definition) and $b_n=c_n^{-1}a_n$, we write 
$$z(t)=\sum_{n\in \N^*}a_n u_n(t)=\sum_{n\in \N^*}b_n w_n(t)= \int_{-S}^{S}\sum_{n\in \N^*}k_{T}(t,s)b_nv_{n}(s)ds.$$
By Cauchy-Schwarz inequality in $L^2(-S,S)$, we deduce
\bna
\nor{z}{L^2(0,T)}^2&=&\int_0^T \left(\int_{-S}^{S}\sum_{n\in \N^*}k_{T}(t,s)b_nv_{n}(s)ds\right)^2dt\leq \int_0^T \nor{k_{T}(t,\cdot)}{L^{2}(-S,S)}^2\int_{-S}^{S}\left|\sum_{n\in \N^*}b_nv_{n}(s)\right|^2dsdt\\
&\leq &C_0 \int_0^T \nor{k_{T}(t,\cdot)}{L^{2}(-S,S)}^2dt\sum_{n\in \N^*}|b_n|^2 ,
\ena
after having used \eqref{orthodual}. Then, we fix $\alpha := 2S^2(1+\eps)$ for $\eps>0$, and next fix $\delta \in (0,1)$ close to $1$ so that $\delta/(1+\delta)<2$ and $S^2/\delta<\alpha/(1+\delta)$, and by~\eqref{e:estim-kT}, $|k_T(t,s)| \leq S e^{\frac{2}{T}\left(\frac{S^2}{\delta}- \frac{\alpha}{1+\delta} \right)} \leq S$ uniformly in $T\in [0,1]$. Therefore, we obtain 
\bna
\nor{z}{L^2(0,T)}^2 
\leq C \sum_{n\in \N^*}c_n^{-2} |a_n|^2\leq \frac{C}{T^3} e^{\frac{8\alpha }{T}} \sum_{n\in \N^*} \beta_n^2 e^{2\beta_{n}^{2}T} |a_n|^2 
= \frac{C}{T^3} e^{\frac{16 S^2(1+\eps) }{T}} \sum_{n\in \N^*} \beta_n^2 e^{2\beta_{n}^{2}T} |a_n|^2 
\ena
after having used~\eqref{estimcl}.
\enp
We have used the following classical lemma that we state and prove only because we did not find any reference precising the constants involved. The proof we present is taken from Gohberg-Krein \cite[Theorem 2.1 p310]{GK:69}.

In the following, we shall say that a sequence $(a_k)_{k\in \N}$ is finite if $a_k \neq 0$ for only a {\em finite number} of indices $k \in \N$.
\begin{lemma}[Biorthogonal family with explicit constants]
\label{lmbiorth}
Let $H$ be a Hilbert space with norm $\nor{\cdot}{H}$, $C_1,C_2>0$ two constants, and $(\varphi_k)_{k\in \N} \in H^\N$ a sequence so that
\bnan
\label{assumorth}
C_1\nor{\sum_{k\in \N}a_{k}\varphi_k}{H}^2\leq \sum_{k\in \N}|a_{k}|^{2}\leq C_2\nor{\sum_{k\in \N}a_{k}\varphi_k}{H}^2
\enan
for any finite sequence $(a_k)_{k\in \N}$. Then, there exists a sequence $(\psi_k)_{k\in \N}$ in $\overline{\vect_{k\in \N} \varphi_k}$ so that 
\bna
\left(\varphi_k,\psi_n \right)_H=\delta_{k,n} , \quad \text{ for all } k,n \in \N,
\ena
and 
\bna
C_2^{-1}\nor{\sum_{k\in \N}a_{k}\psi_k}{H}^2\leq \sum_{k\in \N}|a_{k}|^{2}\leq C_1^{-1}\nor{\sum_{k\in \N}a_{k}\psi_k}{H}^2
\ena
for any finite sequence $(a_k)_{k\in \N}$. 
\end{lemma}
\bnp
Let $(e_k)_{k\in \N}$ be an arbitrary orthonormal basis of the Hilbert space $\widetilde{H}=\overline{\vect_{k\in \N} \varphi_k}$ endowed with the norm $\nor{\cdot}{\widetilde{H}} = \nor{\cdot}{H}$. We define two linear operators, $A$ on $\vect_{k\in \N} e_k$ and $A_1$ on $\vect_{k\in \N} \varphi_k$,  by 
\bna
A\left(\sum_{k\in\N}a_k e_k\right)=\left(\sum_{k\in\N}a_k \varphi_k\right);\quad A_1\left(\sum_{k\in\N}a_k \varphi_k\right)=\left(\sum_{k\in\N}a_k e_k\right)
\ena
 for finite sequences $(a_k)_{k\in \N}$. Note that it is uniquely defined thanks to the orthogonality of the family $(e_k)_{k\in \N}$ and Assumption~\eqref{assumorth}. Assumption~\eqref{assumorth} actually gives more precisely
\bna
\nor{A\left(\sum_{k\in\N}a_k e_k\right)}{H}^2\leq C_1^{-1}\sum_{k\in \N}|a_{k}|^{2}=C_1^{-1} \nor{\sum_{k\in\N}a_k e_k}{H}^2 , \\
\nor{A_1\left(\sum_{k\in\N}a_k \varphi_k\right)}{H}^2=\nor{\sum_{k\in\N}a_k e_k}{H}^2=\sum_{k\in \N}|a_{k}|^{2}\leq C_2\nor{\sum_{k\in\N}a_k \varphi_k}{H}^2 .
\ena
In particular, $A$ and $A_1$ can be extended uniquely by uniform continuity to $\widetilde{H}$ (recall that $\ovl{\vect_{k\in \N} e_k} = \widetilde{H} =\ovl{\vect_{k\in \N} \varphi_k}$ by definition) with $\nor{A}{\widetilde{H}\rightarrow\widetilde{H}} \leq C_1^{-1/2}$ and $\nor{A_1}{\widetilde{H}\rightarrow\widetilde{H}} \leq C_2^{1/2}$. Moreover, they satisfy $AA_1=A_1A=\id_{\widetilde{H}}$. Then, we define $\psi_n:=A_1^*e_n$. With this definition, we have
\bna
\left(\varphi_k,\psi_n \right)_H=\left(\varphi_k,A_1^* e_n \right)_H =\left(A_1\varphi_k,e_n \right)_H=\left(e_k,e_n \right)_H=\delta_{k,n}.
\ena
Moreover,
\bna
\nor{\sum_{k\in \N}a_{k}\psi_k}{H}^2=\nor{A_1^*\left(\sum_{k\in \N}a_{k}e_k\right)}{H}^2\leq \nor{A_1^*}{\widetilde{H}\rightarrow\widetilde{H}}^2\nor{\sum_{k\in \N}a_{k}e_k}{H}^2\leq  C_2\sum_{k\in \N}|a_{k}|^{2},
\ena
\bna
 \sum_{k\in \N}|a_{k}|^{2}=\nor{A^*A_1^*\left(\sum_{k\in \N}a_{k}e_k\right)}{H}^2\leq\nor{A^*}{\widetilde{H}\rightarrow\widetilde{H}}^2\nor{\sum_{k\in \N}a_{k}\psi_k}{H}^2\leq  C_1^{-1}\nor{\sum_{k\in \N}a_{k}\psi_k}{H}^2,
\ena
which concludes the proof of the lemma.
\enp

\section{Proofs of technical results}
\label{A:technical}
In this section, we provide with proofs of some technical results stated in the introduction.
\subsection{Proof of Lemma~\ref{l:phiT}}
\bnp[Proof of Lemma~\ref{l:phiT}]
Note that for $\EE>\min V=V(\xzero)$, $\Phi$ is differentiable at all points where $x_-$ and $x_+$ are, that is for $\EE \in \R\setminus \{V(L), V(0)\}$, with 
\begin{align*}
\Phi'(\EE) & = x_+'(\EE) \sqrt{\EE-V(x_+(\EE))} - x_-' (\EE) \sqrt{\EE-V(x_-(\EE))} + \int_{x_-(\EE)}^{x_+(\EE)} \frac{1}{2\sqrt{\EE -V(s)}} ds \\
 & =  \int_{x_-(\EE)}^{x_+(\EE)} \frac{1}{2 \sqrt{\EE -V(s)}} ds = \frac{1}{4\sqrt{\EE}} T(\EE) .
\end{align*}
As a consequence, we have $\big(\Phi(\EE^2) \big)' = 2 \EE \Phi'(\EE^2) = \frac{1}{2} T(\EE^2)$.
\enp

\subsection{Proof of Lemma~\ref{l:TEB-T}}
\bnp[Proof of Lemma~\ref{l:TEB-T}]
Hence, $T_1$ (defined in~\eqref{e:def-TT1}) and $T_{E,B}$ (defined in~\eqref{TEB}) are linked by: (recall $E_0=V(\xzero) = \frac{|\f'(\xzero)|^2}{4}$)
\begin{align*}
T_{E,B} = \frac{1}{\pi}\int_{V(\xzero)}^{+\infty}\log \left|\frac{x+E+2B}{x-E}\right|\Phi'(x)dx = \frac{1}{\pi}\int_{V(\xzero)}^{+\infty}\log \left|\frac{x+E+2B}{x-E}\right| \frac{T(x)}{4 \sqrt{x}}dx  \leq T_1 \Gamma(E,B ,V(\xzero)) ,
\end{align*}
with 
$$ 
\Gamma(E,B,E_0)  = \frac{1}{\pi}\int_{E_0}^{+\infty}\log \left|\frac{x+E+2B}{x-E}\right|  \frac{1}{4 \sqrt{x}} dx =  \frac{1}{2\pi}\int_{\sqrt{E_0}}^{+\infty}\log \left|\frac{y^2+E+2B}{y^2-E}\right|  dy.
$$
Changing variables in this last integral, we obtain 
$$ 
\Gamma(E,B,E_0)  = \frac{\sqrt{E_0}}{2\pi} \Gamma_0\left(\sqrt{\frac{E+2B}{E_0}} , \sqrt{\frac{E}{E_0}}\right), \quad \text{ with }
\quad \Gamma_0(\alpha, \beta) = \int_{1}^{+\infty}\log \left|\frac{y^2 +\alpha^2}{y^2-\beta^2}\right|  dy
$$
We have obtained 
$$
T_{E,B}  \leq \frac{T_1 \sqrt{E_0}}{2\pi} \Gamma_0\left(\sqrt{\frac{E+2B}{E_0}} , \sqrt{\frac{E}{E_0}}\right),
$$
and we now compute $\Gamma_0(\alpha, \beta)$ for $\alpha,\beta \geq 1$ (since $E\geq E_0$).
We set, for $\alpha,\beta \geq 1$
\begin{align*}
F_\alpha(y) & := y \left( \log(y^2+\alpha^2) - 2\right) + 2 \alpha \arctan \left( \frac{y}{\alpha}\right) , \quad y \in \R , \\
 G_\beta^+(y)& := y \left( \log(y^2-\beta^2) - 2\right) + \beta \log \left( \frac{y+\beta}{y-\beta}\right) ,\quad y > \beta ,\\
  G_\beta^-(y)& := y \left( \log(\beta^2-y^2) - 2\right) + \beta \log \left( \frac{\beta+y}{\beta-y}\right) ,\quad 1\leq y < \beta ,
\end{align*}
and notice that $F_\alpha'(y) =\log(y^2+\alpha^2)$ for all $y\in \R$, $(G_\beta^+)'(y) =\log(y^2-\beta^2)$ for all $y > \beta$ and $(G_\beta^-)'(y) =\log(\beta^2-y^2)$ for all $1\leq y<\beta$. Moreover, we notice that $\lim_{y \to \beta^+} G_\beta^+(y) = 2\beta (\log(2\beta)-1) = \lim_{y \to \beta^-} G_\beta^-(y)$, and thus $G_\beta : = \mathds{1}_{(-\infty,\beta)}G_\beta^-+\mathds{1}_{[\beta,+\infty)}G_\beta^+$ is continuous.
As a consequence, we can compute explicitly for $\alpha,\beta \geq 1$
\begin{align*} 
\Gamma_0(\alpha,\beta) & = \left[ F_\alpha(y) - G_\beta(y)\right]_{1}^\infty 
= \lim_{y \to + \infty} (F_\alpha(y) - G_\beta^+(y)) - (F_\alpha(1) - G_\beta^-(1) ) = \pi \alpha - F_\alpha(1) + G_\beta^-(1) \\
& =\pi \alpha - \log(1+\alpha^2)  - 2 \alpha \arctan \left( \frac{1}{\alpha}\right) 
+ \log(\beta^2-1)  + \beta \log \left( \frac{\beta+1}{\beta-1}\right) ,
\end{align*}
which is the sought result.
\enp

\subsection{Proof of Lemma~\ref{l:compare}}
\bnp[Proof of Lemma~\ref{l:compare}]
 Note first that, recalling that $W_E = \frac{\f}{2}+d_{A,E}$ and $\widetilde{W_E} = \frac{\f}{2}-d_{A,E}$, and that $d_{A,E}=0$ on $K_E$, we obtain
  $$\min_{\M} W_E\leq \min_{K_E} W_E=\min_{K_E} \widetilde{W_E}\leq \sup_{[0,L]}\widetilde{W_E} ,$$
  and thus~\eqref{e:compG} holds true.
Next, according to Lemma~\ref{l:TEB-T}, the quantities $S_{\ref{thmlower+-},E,B}$ and $S_{\ref{thmcontrol1D}}$ are linked by  
\begin{align*}
0= S_{\ref{t:estim-Cobs-expo-1D}}\leq S_{\ref{thmlower+-},E,B} = T_{E,B}  \leq \frac{S_{\ref{thmcontrol1D}}}{4\pi \sqrt{2}} \Gamma_0\left(\sqrt{\frac{E+2B}{E_0}} , \sqrt{\frac{E}{E_0}}\right),
\end{align*}
As a consequence, using that $2 \frac{E+B}{E_0} = \left(\sqrt{\frac{E+2B}{E_0}}\right)^2 + \left(\sqrt{\frac{E}{E_0}}\right)^2$
$$
\frac{S_{\ref{thmlower+-},E,B}}{E+B} \leq  \frac{E_0}{E+B} \frac{S_{\ref{thmcontrol1D}}}{E_0} \frac{1}{4\pi \sqrt{2}} \Gamma_0\left(\sqrt{\frac{E+2B}{E_0}},\sqrt{\frac{E}{E_0}}\right) \leq  \frac{S_{\ref{thmcontrol1D}}}{E_0} \frac{1}{2\pi \sqrt{2}} \sup_{\alpha\geq \beta\geq1}  \frac{1}{\alpha^2+\beta^2} \Gamma_0(\alpha,\beta) .
$$
 Note that the supremum is actually a maximum according to Lemma~\ref{l:TEB-T}, whence~\eqref{e:linky}.
\enp

\subsection{Proof of Lemma~\ref{l:WE}}

\bnp[Proof of Lemma~\ref{l:WE}]
We write $\f=\pm g$ with $g$ strictly increasing $[0,L]$; the case $\f$ increasing (resp. decreasing) will be denoted the case $+$ (resp. $-$) and in both cases we have $g'\geq0$.
 
  Note that we only need to prove the result in the case $E \in V([0,L])$, for if $E> \max V$, we have $d_{A,E}=0$ identically on $[0,L]$ and thus $W_E=\widetilde{W_E}=\frac{\f}{2}$ and the result follows.

For $E\geq E_{0}$, we recall that $x_\pm(E)$ are defined just after~\eqref{e:def-TT1}.
Outside of $K_{E}=[x_-(E),x_+(E)]$, we have $d_{A,E}'(x) = \sqrt{\frac{|g'(x)|^{2}}{4} - E  }$ for $x\geq x_+(E)$, and $d_{A,E}'(x) = -\sqrt{\frac{|g'(x)|^{2}}{4} - E  }$  for $x\leq x_-(E)$. As a consequence, recalling the definition of $W_E$ in~\eqref{e:def-WE}, we have
\begin{align*}
\begin{array}{ll}
W_E' (x)=-\sqrt{\frac{|g'(x)|^{2}}{4} - E}\pm \frac{g'(x)}{2}=\frac{g'(x)}{2} \left(-\sqrt{ 1- \frac{4E}{|g'(x)|^{2}}}\pm 1\right)& \textnormal{ for }x\leq x_-(E), \\
W_E' (x)=\pm \frac{g'(x)}{2}& \textnormal{ for }x\in [x_-(E),x_+(E)] ,\\
W_E' (x)=\sqrt{\frac{|g'(x)|^{2}}{4} - E}\pm \frac{g'(x)}{2}=\frac{g'(x)}{2} \left(\sqrt{ 1- \frac{4E}{|g'(x)|^{2}}}\pm 1\right)& \textnormal{ for }x\geq x_+(E).
\end{array}
\end{align*}
Outside of $K_E = [x_-(E),x_+(E)]$, we always have $0\leq \sqrt{ 1- \frac{4E}{|g'(x)|^{2}}}\leq 1$, so that for $x\in [0,L]$, $W_{E}$ is increasing in the case $+$  and decreasing in the case $-$.

Concerning $\widetilde{W_E}(s)=\frac{\f}{2}(s)-d_{A,E}(s)$, we compute similarly
\begin{align*}
\begin{array}{ll}
\widetilde{W_E}'(x)=\pm \frac{g'(x)}{2}+\sqrt{\frac{|g'(x)|^{2}}{4} - E}=\frac{g'(x)}{2} \left(\pm 1+\sqrt{ 1- \frac{4E}{|g'(x)|^{2}}}\right) & \textnormal{ for }x\leq x_-(E) , \\
\widetilde{W_E}'(x)=\pm \frac{g'(x)}{2}& \textnormal{ for }x\in [x_-(E),x_+(E)] , \\
\widetilde{W_E}'(x)=\pm \frac{g'(x)}{2}-\sqrt{\frac{|g'(x)|^{2}}{4} - E}=\frac{g'(x)}{2} \left(\pm 1-\sqrt{ 1- \frac{4E}{|g'(x)|^{2}}}\right) & \textnormal{ for }x\geq x_+(E) .
\end{array}
\end{align*}
So, as for $W_E$, the function $\widetilde{W_E}(s)$ is increasing on $[0,L]$ in the case $+$  and decreasing in the case $-$. 

To summarize, in the case $+$, we have
\bna
\min_{[0,L]} W_E=W_E(0); \quad \sup_{[0,L]}\widetilde{W_E}=\widetilde{W_E}(L);  
\ena
while in the case $-$, we have 
\bna
\min_{[0,L]} W_E=W_E(L) ; \quad \sup_{[0,L]}\widetilde{W_E}=\widetilde{W_E}(0); \ena
The statements concerning $G_{\ref{t:estim-Cobs-expo-1D},E}=G_{\ref{thmcontrol1D},E}=W_E(0) - \min_{\M} W_E =0$ and 
$G_{\ref{thmlower+-},E}=W_E(0)-\sup_{[0,L]}\widetilde{W_E}$ are then direct consequences of the above results.

Concerning the last properties of these functions, we notice that $d_{A, E}(0)$ and $d_{A, E}(L)$ are non-increasing functions of $E$. This proves that  $G_{\ref{thmlower+-},E}$ is non increasing in both cases. 

Finally, if $g$ is odd, then $g'=|g'|$ is even and $d_{A,E}$ is even. All sought simplifications follow.
\enp

\subsection{Elementary computations}
\label{calculusI}
We collect here two elementary 
 lemmata, that are used in the proof of Theorem~\ref{thmcontrol1D}.
\begin{lemma}
\label{lmcalcI}
Let $a,b>0$, and set $F(m) := \frac{a}{(1-m)}-bm$ for $m\in [0,1)$. Then,  
\begin{itemize}
\item if $a\leq b$ then, $\min_{m\in [0,1)} F(m) = F\left( 1-\sqrt{\frac{a}{b}} \right) = 2\sqrt{ab}-b$;
\item if $a\geq b$ then,  $\min_{m\in [0,1)} F(m) = F(0) = a$.
\end{itemize}
\end{lemma}
\bnp We simply write
$F'(m)=\frac{a}{(1-m)^2}-b\geq 0 \Leftrightarrow (1-m)^2\leq \frac{a}{b}\Leftrightarrow 1-m\leq \sqrt{\frac{a}{b}}\Leftrightarrow m\geq 1-\sqrt{\frac{a}{b}}$.
\enp
\begin{lemma}
\label{lmcalctmin}
Let $a,b, c>0$. Let $G(\TO)=\min_{m\in[0,1)} \frac{a}{(1-m)\TO}-bm\TO+c$. Then, 
\bna
G(T)<0 \quad \text{ if and only if } \quad \TO> 2\sqrt{\frac{a}{b}}+\frac{c}{b} .
\ena
 \end{lemma}
\bnp
In the case $a\geq b\TO^2$, it follows from Lemma~\ref{lmcalcI} that $G(T)\geq 0$. In case $a\leq b\TO^2$, we have from Lemma~\ref{lmcalcI} that
$G(T)=\min_{m\in[0,1)} \frac{a}{(1-m)\TO}-bm\TO+c=2\sqrt{ab}-b\TO+c$ which gives the result.
\enp

\bigskip
\noindent
{\em Acknowledgements.} 
The authors would like to thank Franck Sueur for interesting discussions on the subject of this article and the anonymous referee for his/her useful suggestions, which helped to improve the exposition of the article.
The second author is partially supported by the Agence Nationale de la Recherche under grants SALVE (ANR-19-CE40-0004) and ADYCT (ANR-20-CE40-0017).

\small
\bibliographystyle{alpha}
\bibliography{bibli}
\end{document}

%% file: mymacros.tex
\newtheorem{lemma}{Lemma}[section]
\newtheorem{theorem}[lemma]{Theorem}
\newtheorem{proposition}[lemma]{Proposition}
\newtheorem{corollary}[lemma]{Corollary}
\newtheorem{assumption}[lemma]{Assumption}

\theoremstyle{definition}

\newtheorem{definition}[lemma]{Definition}
\newtheorem{remark}[lemma]{Remark}

\makeatletter
\def\keywords{
    \vspace{1ex}
    \noindent
    \if@twocolumn
      \small{\bf  Keywords}\/---$\!$    \else
      \begin{center}\small\ {\bf Keywords}\end{center}\quotation\small
    \fi}
\def\endkeywords{\vspace{0.6em}\par\if@twocolumn\else\endquotation\fi
    \normalsize\rm}
\makeatother

\renewcommand{\O}{\ensuremath{\mathcal O}}

\renewcommand{\L}{\ensuremath{\mathcal L}}

\newcommand{\x}{\ensuremath{\bf x}}

\DeclareMathOperator{\Sp}{Sp}
\DeclareMathOperator{\loc}{loc}

\DeclareMathOperator{\Vol}{Vol}

\DeclareMathOperator{\Rk}{Rk}

\newcommand{\mb}[1]{\ensuremath{\mathbb{#1}}}
\newcommand{\N}{{\mb{N}}}

\newcommand{\R}{{\mb{R}}}
\newcommand{\C}{{\mb{C}}}

\newcommand{\f}{\ensuremath{\mathfrak{f}}}

\newcommand{\eps}{\varepsilon}

\newcommand{\M}{\ensuremath{\mathcal M}}

\let \Im \relax
\DeclareMathOperator{\Im}{Im}

\newcommand{\ovl}[1]{\overline{#1}}

\DeclareMathOperator{\supp}{supp}

\DeclareMathOperator{\dist}{dist}

\DeclareMathOperator{\vect}{span}

\DeclareMathOperator{\id}{Id}

\renewcommand{\d}{\ensuremath{\partial}}

\def\x{x}

\newcommand{\Z}{\mathbb Z}
\newcommand{\F}{\mathscr F}





%% file: Vanishing1D.bbl
\begin{thebibliography}{AKBGBdT11}

\bibitem[AKBGBdT11]{ABGT:11b}
Farid Ammar-Khodja, Assia Benabdallah, Manuel Gonz\'{a}lez-Burgos, and Luz
  de~Teresa.
\newblock The {K}alman condition for the boundary controllability of coupled
  parabolic systems. {B}ounds on biorthogonal families to complex matrix
  exponentials.
\newblock {\em J. Math. Pures Appl. (9)}, 96(6):555--590, 2011.

\bibitem[All98]{Allibert:98}
Brice Allibert.
\newblock Contr\^ole analytique de l'\'equation des ondes et de l'\'equation de
  {S}chr\"odinger sur des surfaces de r\'evolution.
\newblock {\em Comm. Partial Differential Equations}, 23(9-10):1493--1556,
  1998.

\bibitem[AM19a]{YM:19bis}
Youcef Amirat and Arnaud M\"{u}nch.
\newblock Asymptotic analysis of an advection-diffusion equation and
  application to boundary controllability.
\newblock {\em Asymptot. Anal.}, 112(1-2):59--106, 2019.

\bibitem[AM19b]{YM:19}
Youcef Amirat and Arnaud M\"{u}nch.
\newblock On the controllability of an advection-diffusion equation with
  respect to the diffusion parameter: asymptotic analysis and numerical
  simulations.
\newblock {\em Acta Math. Appl. Sin. Engl. Ser.}, 35(1):54--110, 2019.

\bibitem[BDE20]{BDE:20}
Karine Beauchard, J\'er\'emi Dard\'e, and Sylvain Ervedoza.
\newblock {Minimal time issues for the observability of Grushin-type
  equations}.
\newblock {\em Annales de l'Institut Fourier}, 70(1):247--312, 2020.

\bibitem[CF96]{CF:96}
{Jean-Michel} Coron and Andrei~V. Fursikov.
\newblock Global exact controllability of the $2$-{D} {N}avier-{S}tokes
  equations on a manifold without boundary.
\newblock {\em Russian J. Math. Phys.}, 4:429--448, 1996.

\bibitem[CG05]{CG:05}
{Jean-Michel} Coron and Sergio Guerrero.
\newblock Singular optimal control: A linear $1$-{D} parabolic-hyperbolic
  example.
\newblock {\em Asympt. Anal.}, 44:237--257, 2005.

\bibitem[Cha43]{Chan:43}
Subrahmanyan Chandresekhar.
\newblock Stochastic problems in physics and astronomy.
\newblock {\em Rev. Modern Phys.}, 15:1--89, 1943.

\bibitem[Cha09]{ChaNavier:09}
Marianne Chapouly.
\newblock On the global null controllability of a {N}avier-{S}tokes system with
  {N}avier slip boundary conditions.
\newblock {\em J. Differential Equations}, 247:2094--2123, 2009.

\bibitem[CMS20]{CMS:19}
Jean-Michel Coron, Fr\'{e}d\'{e}ric Marbach, and Franck Sueur.
\newblock Small-time global exact controllability of the {N}avier-{S}tokes
  equation with {N}avier slip-with-friction boundary conditions.
\newblock {\em J. Eur. Math. Soc. (JEMS)}, 22(5):1625--1673, 2020.

\bibitem[CMV20]{CMV:20}
Piermarco Cannarsa, Patrick Martinez, and Judith Vancostenoble.
\newblock Precise estimates for biorthogonal families under asymptotic gap
  conditions.
\newblock {\em Discrete Contin. Dyn. Syst. Ser. S}, 13(5):1441--1472, 2020.

\bibitem[Cor96]{Cor:96}
{Jean-Michel} Coron.
\newblock On the controllability of the $2$-{D} incompressible
  {N}avier-{S}tokes equations with the {N}avier slip boundary conditions.
\newblock {\em ESAIM Control Optim. Calc. Var.}, 1:35--75, 1996.

\bibitem[Cor07]{Cor:book}
Jean-Michel Coron.
\newblock {\em Control and nonlinearity}, volume 136 of {\em Mathematical
  Surveys and Monographs}.
\newblock American Mathematical Society, Providence, RI, 2007.

\bibitem[Daf00]{Daf:00}
Constantine~M. Dafermos.
\newblock {\em Hyperbolic conservation laws in continuum physics}.
\newblock Springer-Verlag, Berlin, 2000.

\bibitem[DE19]{Darde:17}
J\'{e}r\'{e}mi Dard\'{e} and Sylvain Ervedoza.
\newblock On the cost of observability in small times for the one-dimensional
  heat equation.
\newblock {\em Anal. PDE}, 12(6):1455--1488, 2019.

\bibitem[DR77]{DR:77}
Szymon Dolecki and David~L. Russell.
\newblock A general theory of observation and control.
\newblock {\em SIAM J. Control Optim.}, 15(2):185--220, 1977.

\bibitem[DS99]{DS:book}
Mouez Dimassi and Johannes Sj\"ostrand.
\newblock {\em Spectral asymptotics in the semi-classical limit}, volume 268 of
  {\em London Mathematical Society Lecture Note Series}.
\newblock Cambridge University Press, Cambridge, 1999.

\bibitem[EZ11a]{EZ:11}
Sylvain Ervedoza and Enrique Zuazua.
\newblock Observability of heat processes by transmutation without geometric
  restrictions.
\newblock {\em Math. Control Relat. Fields}, 1(2):177--187, 2011.

\bibitem[EZ11b]{EZ:11s}
Sylvain Ervedoza and Enrique Zuazua.
\newblock Sharp observability estimates for heat equations.
\newblock {\em Arch. Ration. Mech. Anal.}, 202(3):975--1017, 2011.

\bibitem[FI96]{FI:96}
Andrei~V. Fursikov and Oleg~Yu. Imanuvilov.
\newblock {\em Controllability of evolution equations}, volume~34 of {\em
  Lecture Notes Series}.
\newblock Seoul National University Research Institute of Mathematics Global
  Analysis Research Center, Seoul, 1996.

\bibitem[FR71]{FR:71}
Hector~O. Fattorini and David~L. Russell.
\newblock Exact controllability theorems for linear parabolic equations in one
  space dimension.
\newblock {\em Arch. Rational Mech. Anal.}, 43:272--292, 1971.

\bibitem[FR75]{FR:74}
Hector~O. Fattorini and David~L. Russell.
\newblock Uniform bounds on biorthogonal functions for real exponentials with
  an application to the control theory of parabolic equations.
\newblock {\em Quart. Appl. Math.}, 32:45--69, 1974/75.

\bibitem[GG07]{GG:07}
Olivier Glass and Sergio Guerrero.
\newblock On the uniform controllability of the {B}urgers equation.
\newblock {\em SIAM J. Control Optim.}, 46:1211--1238, 2007.

\bibitem[GG08]{GG:08}
Olivier Glass and Sergio Guerrero.
\newblock Some exact controllability results for the linear {K}d{V} equation
  and uniform controllability in the zero-dispersion limit.
\newblock {\em Asymptot. Anal.}, 60:61--100, 2008.

\bibitem[GG09]{GG:09}
Olivier Glass and Sergio Guerrero.
\newblock Uniform controllability of a transport equation in zero
  diffusion-dispersion limit.
\newblock {\em Math. Models Methods Appl. Sci.}, 19:1567--1601, 2009.

\bibitem[GK69]{GK:69}
Israel~C. Gohberg and Mark~G. Krein.
\newblock {\em Introduction to the theory of linear non-selfadjoint operators},
  volume~18 of {\em Transl. Math. Monogr.}
\newblock Amer. Math. Soc., Providence, R.I., 1969.

\bibitem[GL07]{GL:07}
Sergio Guerrero and Gilles Lebeau.
\newblock Singular optimal control for a transport-diffusion equation.
\newblock {\em Comm. Partial Differential Equations}, 32:1813--1836, 2007.

\bibitem[Gla10]{Gla:10}
Olivier Glass.
\newblock A complex-analytic approach to the problem of uniform controllability
  of a transport equation in the vanishing viscosity limit.
\newblock {\em J. Funct. Anal.}, 258:852--868, 2010.

\bibitem[Han91]{H:91}
Scott~W. Hansen.
\newblock Bounds on functions biorthogonal to sets of complex exponentials;
  control of damped elastic systems.
\newblock {\em J. Math. Anal. Appl.}, 158(2):487--508, 1991.

\bibitem[Har89]{H:89-2}
Alain Haraux.
\newblock S\'{e}ries lacunaires et contr\^{o}le semi-interne des vibrations
  d'une plaque rectangulaire.
\newblock {\em J. Math. Pures Appl. (9)}, 68(4):457--465 (1990), 1989.

\bibitem[Hel88]{Helffer:booksemiclassic}
Bernard Helffer.
\newblock {\em Semi-classical analysis for the {S}chr\"odinger operator and
  applications}, volume 1336 of {\em Lecture Notes in Mathematics}.
\newblock Springer-Verlag, Berlin, 1988.

\bibitem[HR84]{HR:84}
Bernard Helffer and Didier Robert.
\newblock Puits de potentiel g\'{e}n\'{e}ralis\'{e}s et asymptotique
  semi-classique.
\newblock {\em Ann. Inst. H. Poincar\'{e} Phys. Th\'{e}or.}, 41(3):291--331,
  1984.

\bibitem[HS84]{HS:84}
Bernard Helffer and Johannes Sj\"ostrand.
\newblock Multiple wells in the semiclassical limit. {I}.
\newblock {\em Comm. Partial Differential Equations}, 9(4):337--408, 1984.

\bibitem[HS85]{HS:85}
Bernard Helffer and Johannes Sj\"ostrand.
\newblock Puits multiples en m\'{e}canique semi-classique. {IV}. \'{E}tude du
  complexe de {W}itten.
\newblock {\em Comm. Partial Differential Equations}, 10(3):245--340, 1985.

\bibitem[KL05]{KL:05}
Vilmos Komornik and Paola Loreti.
\newblock {\em Fourier series in control theory}.
\newblock Springer Monographs in Mathematics. Springer-Verlag, New York, 2005.

\bibitem[Koo88]{Koosis:bookI}
Paul Koosis.
\newblock {\em The logarithmic integral. {I}}, volume~12 of {\em Cambridge
  Studies in Advanced Mathematics}.
\newblock Cambridge University Press, Cambridge, 1988.

\bibitem[Kru70]{Kru:70}
Stanislav~N. Kru\v{z}kov.
\newblock First order quasilinear equations with several independent variables.
  (russian).
\newblock {\em Mat. Sb. (N.S.)}, 81:228--255, 1970.

\bibitem[L{\'e}a10]{Lea:10}
Matthieu L{\'e}autaud.
\newblock Spectral inequalities for non-selfadjoint elliptic operators and
  application to the null-controllability of parabolic systems.
\newblock {\em J. Funct. Anal.}, 258:2739--2778, 2010.

\bibitem[L{\'e}a12]{Lea:11}
Matthieu L{\'e}autaud.
\newblock Uniform controllability of scalar conservation laws in the vanishing
  viscosity limit.
\newblock {\em SIAM J. Control Optim.}, 50(3):1661--1699, 2012.

\bibitem[Lis12]{Lissy:12}
Pierre Lissy.
\newblock A link between the cost of fast controls for the 1-{D} heat equation
  and the uniform controllability of a 1-{D} transport-diffusion equation.
\newblock {\em C. R. Math. Acad. Sci. Paris}, 350(11-12):591--595, 2012.

\bibitem[Lis14]{Lissy:14}
Pierre Lissy.
\newblock An application of a conjecture due to {E}rvedoza and {Z}uazua
  concerning the observability of the heat equation in small time to a
  conjecture due to {C}oron and {G}uerrero concerning the uniform
  controllability of a convection-diffusion equation in the vanishing viscosity
  limit.
\newblock {\em Systems Control Lett.}, 69:98--102, 2014.

\bibitem[Lis15]{Lissy:15}
Pierre Lissy.
\newblock Explicit lower bounds for the cost of fast controls for some 1-{D}
  parabolic or dispersive equations, and a new lower bound concerning the
  uniform controllability of the 1-{D} transport-diffusion equation.
\newblock {\em J. Differential Equations}, 259(10):5331--5352, 2015.

\bibitem[LL21]{LL:18vanish}
Camille Laurent and Matthieu L\'{e}autaud.
\newblock On uniform observability of gradient flows in the vanishing viscosity
  limit.
\newblock {\em J. \'{E}c. polytech. Math.}, 8:439--506, 2021.

\bibitem[LL22a]{LL:17hypo}
Camille Laurent and Matthieu L\'{e}autaud.
\newblock Tunneling estimates and approximate controllability for hypoelliptic
  equations.
\newblock {\em Mem. Amer. Math. Soc.}, 276(1357):vi+95, 2022.

\bibitem[LL22b]{LL:22}
Camille Laurent and Matthieu L\'eautaud.
\newblock Uniform observation of semiclassical {S}chr\"odinger eigenfunctions
  on an interval.
\newblock {\em submitted, arXiv:2203.03271}, 2022.

\bibitem[LR95]{LR:95}
Gilles Lebeau and Luc Robbiano.
\newblock Contr\^ole exact de l'\'equation de la chaleur.
\newblock {\em Comm. Partial Differential Equations}, 20:335--356, 1995.

\bibitem[Mil06]{Miller:06b}
Luc Miller.
\newblock The control transmutation method and the cost of fast controls.
\newblock {\em SIAM J. Control Optim.}, 45(2):762--772, 2006.

\bibitem[SM79]{SM:79}
Zeev Schuss and Bernard~J. Matkowsky.
\newblock The exit problem: a new approach to diffusion across potential
  barriers.
\newblock {\em SIAM J. Appl. Math.}, 36(3):604--623, 1979.

\bibitem[Wit82]{Witten:82}
Edward Witten.
\newblock Supersymmetry and {M}orse theory.
\newblock {\em J. Differential Geometry}, 17(4):661--692 (1983), 1982.

\end{thebibliography}
